\documentclass[a4paper]{amsart} 
\usepackage{amsfonts,amsmath,amssymb,latexsym,epsfig,graphicx,color} 
\usepackage[matrix,arrow,curve]{xy} 
\newtheorem{theorem}{Theorem} 
\newtheorem{corollary}{Corollary} 
\newtheorem{lemma}{Lemma} 
\newtheorem{proposition}{Proposition} 
\newtheorem{definition}{Definition} 
 
\theoremstyle{definition} 
\newtheorem{remark}{Remark} 
\newtheorem{example}{Example} 
\newcommand{\para}{\begin{PARA}\rm} 
\newcommand{\arap}{\end{PARA}\rm} 
\newcommand{\dfn}{\begin{definition}\rm} 
\newcommand{\nfd}{\end{definition}\rm} 
\newcommand{\rmk}{\begin{remark}\rm} 
\newcommand{\kmr}{\end{remark}\rm} 
\newcommand{\xmpl}{\begin{example}\rm} 
\newcommand{\lpmx}{\end{example}\rm} 
\newcommand{\cA}{\mathcal{A}}

\newcommand{\cC}{\mathcal{C}} 
 
\newcommand{\cF}{\mathcal{F}} 
\newcommand{\cE}{\mathcal{E}} 
 
\newcommand{\cH}{\mathcal{H}} 
 
\newcommand{\cJ}{\mathcal{J}} 
 
\newcommand{\cL}{\mathcal{L}} 
\newcommand{\cM}{\mathcal{M}}

\newcommand{\cP}{\mathcal{P}}

\newcommand{\cU}{\mathcal{U}} 
 
\newcommand{\cV}{\mathcal{V}}

\newcommand{\oeps}{{\overline{\epsilon}}} 
\newcommand{\og}{{\overline{\gamma}}} 
\newcommand{\ug}{{\underline{\gamma}}} 
\newcommand{\oev}{\overline{\mathrm{ev}}} 
\newcommand{\uev}{\underline{\mathrm{ev}}} 
\newcommand{\oL}{{\overline{L}}} 
\newcommand{\uL}{{\underline{L}}}

\newcommand{\oT}{{\overline{T}}} 
\newcommand{\uT}{{\underline{T}}}

\newcommand{\oSigma}{\overline{\Sigma}} 
\newcommand{\uSigma}{\underline{\Sigma}} 
\newcommand{\oz}{\overline{z}} 
\newcommand{\uz}{\underline{z}}

\newcommand{\utheta}{\underline{\theta}} 
\newcommand{\one} 
{{{\mathchoice \mathrm{ 1\mskip-4mu l} \mathrm{ 1\mskip-4mu l} 
\mathrm{ 1\mskip-4.5mu l} \mathrm{ 1\mskip-5mu l}}}} 

\newcommand{\C}{{\mathbb{C}}}

\newcommand{\Q}{{\mathbb{Q}}} 
\newcommand{\R}{{\mathbb{R}}}

\renewcommand{\u}{{\mathbf{u}}}

\newcommand{\Z}{{\mathbb{Z}}} 
\newcommand{\Arg}{\mathrm{Arg}}  
\newcommand{\im}{\mathrm{ im }}        
 
\newcommand{\ind}{\mathrm{ind}} 
\newcommand{\Jreg}{\cJ_{\mathrm{reg}}}   
 
\newcommand{\Freg}{\cF_{\mathrm{reg}}} 
 
\newcommand{\ev}{\mathrm{ev}} 

\newcommand{\eps}{{\varepsilon}} 
\newcommand{\om}{{\omega}}

\newcommand{\ta}{{\widetilde{a}}}  
 
\newcommand{\tC}{{\widetilde{C}}} 
 
\newcommand{\tgamma}{{\widetilde{\gamma}}} 
\newcommand{\tp}{{\widetilde{p}}} 
\newcommand{\tq}{{\widetilde{q}}}

\newcommand{\tf}{{\widetilde{f}}} 
\newcommand{\tu}{{\widetilde{u}}}

\def\NABLA#1{{\mathop{\nabla\kern-.5ex\lower1ex\hbox{$#1$}}}} 
\def\Nabla#1{\nabla\kern-.5ex{}_{#1}} 
\def\Tabla#1{\Tilde\nabla\kern-.5ex{}_{#1}} 
\renewcommand{\Tilde}{\widetilde}

\newcommand{\p}{{\partial}} 
\newcommand{\dbar}{{\bar\partial}}

\newenvironment{enum}{\begin{enumerate}}
{\end{enumerate}}

\begin{document} 
 
\title{An exact sequence for contact- 
  and symplectic homology} 
 
\author[Bourgeois and Oancea]{Fr\'ed\'eric {\sc Bourgeois}, \
Alexandru \sc{Oancea} 
           \\ \quad \\
         {\it Universit\'e Libre de Bruxelles, B-1050 Bruxelles,
Belgium} \\
         {\it Universit\'e Louis Pasteur, F-67084 Strasbourg, France}
\\ 
{\tt fbourgeo@ulb.ac.be,\qquad oancea@math.u-strasbg.fr}
}

\date{15 October 2008} 
 
\maketitle

 
\begin{abstract}  
A symplectic manifold $W$ with contact type boundary 
$M = \partial W$ induces a linearization of the contact homology of  
$M$ with corresponding  linearized contact homology $HC(M)$. 
We establish a Gysin-type exact sequence in which the symplectic homology 
$SH(W)$ of $W$ maps to $HC(M)$, which in turn maps to $HC(M)$, by a map 
of degree $-2$, which then maps to $SH(W)$.     
Furthermore, we give a description of the degree $-2$ map in terms of rational 
holomorphic curves with constrained asymptotic markers, in 
the symplectization of $M$. 
\end{abstract} 
 
\tableofcontents 
 
 
\section{Introduction} 
 
Let $(W,\om)$ be a compact symplectic manifold with contact type\break 
boundary $M:=\p W$. This means that there exists a vector field $X$ 
defined in a neighbourhood of $M$, transverse and pointing outwards 
along $M$, and such that 
$$ 
\cL _X \om = \om. 
$$ 
Such an $X$ is called a {\bf Liouville vector field}. The $1$-form 
$\lambda:=(\iota_X\om)|_M$ is a contact form on $M$. We denote by 
$\xi$ the contact distribution defined by $\lambda$ and we call 
$(W,\om)$ a {\bf filling} of $(M,\xi)$.  
 
We assume throughout the paper that $(W,\om)$ satisfies the condition 
\begin{equation} \label{eq:asph} 
\int_{T^2} f^*\om =0 \quad \mbox{for all smooth } f:T^2\to W, 
\end{equation} 
where $T^2$ is the $2$-torus. This condition guarantees that the 
energy of a Floer trajectory (for a definition, see for example 
\cite[Section 2]{BOauto})  
does not depend on its homology class, but only on its endpoints.  
Our main class of examples is provided by exact symplectic forms.  
  
Theorem~\ref{thm:intro} ties together the symplectic homology 
groups of   
$(W,\om)$ and the linearized contact homology groups of 
$(M,\xi)$. Both these invariants encode   
algebraically the dynamics 
of the same vector field, the {\bf Reeb vector field} $R_\lambda$  
defined by $\ker \, \om|_M = \langle R_\lambda \rangle$ and 
$\lambda(R_\lambda)=1$. But their natures are quite different: the 
former belongs to the realm of Floer theory~\cite{FH1,V}, whereas the 
latter belongs to the realm of symplectic field theory 
(SFT)~\cite{EGH}. Our result can be read as a way to make 
symplectic homology fit into SFT.   
 
Let us introduce some relevant notation. Given 
a free homotopy class $a$ of loops in $W$ we denote by $SH^a_*(W,\om)$ 
the symplectic homology groups of $(W,\om)$ in the homotopy class 
$a$. The free homotopy class of the constant loop will be denoted by $0$. 
We also denote by $SH_*^+(W,\om)$ the symplectic homology groups 
in the trivial homotopy class truncated at a small positive value 
of the action functional. We refer to Section~\ref{sec:symplectic} for the  
definitions.    
 
Let $i:M\hookrightarrow W$ be the inclusion. Given a free homotopy 
class $a$ of loops in $W$ we denote by $i^{-1}(a)$ the set of free 
homotopy classes in $M$ which are mapped to $a$ via $i$, and we use 
the convention $i^{-1}(+):=i^{-1}(0)$.   
We denote by $HC_*^{i^{-1}(a)}(M,\xi)$ the linearized    
contact homology groups of $(M,\xi)$ based on closed Reeb orbits whose 
free homotopy class belongs to $i^{-1}(a)$. We  
refer to Section~\ref{sec:contact} for the definition. 
 
Both $SH_*^a(W,\om)$ and $HC_*^{i^{-1}(a)}(M,\xi)$ are 
defined over the Novikov ring $\Lambda_\om$ with $\Q$-coefficients 
consisting of formal combinations $\lambda:=\sum_{A\in H_2(W;\Z)} 
\lambda_A e^A$, $\lambda_A\in\Q$ such that  
$$ 
\#\{ A | \lambda_A\neq 0, \om(A) \le c\} <\infty 
$$ 
for all $c>0$. The multiplication in $\Lambda_\om$ is given by the 
convolution product.  
 
We assume the existence of an almost complex structure $J$ such that 
linearized contact homology is defined. This means that $J$ needs to 
be regular for rigid holomorphic planes in the symplectic completion 
of $W$, as well as for  rational holomorphic curves with one  
positive puncture in the symplectization of $M$ satisfying the following  
property. These curves are asymptotic, at all negative 
punctures except at most one,  to Reeb orbits  
which can be capped with rigid holomorphic planes in  
the symplectic completion of $W$.   
We refer to Section~\ref{sec:contact}, Remark~\ref{rmk:transv_exples} 
for a discussion of these regularity 
assumptions. We expect this technical assumption to be completely 
removed using the new ongoing approach to transversality by Cieliebak 
and Mohnke (see~\cite{CM} for the symplectic case), or using the polyfold  
theory developed by Hofer Wysocki and Zehnder \cite{H,HWZ2}.

\begin{theorem} \label{thm:intro} 
If $a\neq 0$ or $a=+$ there exists a long exact sequence  
\begin{eqnarray} \label{eq:intro} 
\ldots \to SH_{k-(n-3)}^a(W,\om) \to HC_{k}^{i^{-1}(a)}(M,\xi) 
\stackrel D \to 
HC_{k-2}^{i^{-1}(a)}(M,\xi) \to \\ 
\hspace{4.5cm}  \to SH_{k-1-(n-3)}^a(W,\om) \to 
\ldots  \nonumber & 
\end{eqnarray} 
Moreover, the map $D$ can be described exclusively in terms of 
rational holomorphic curves with constrained asymptotic markers 
in the symplectization of $(M,\xi)$, and of    
rigid holomorphic planes in the symplectic completion of $W$.   
\end{theorem}  
 
The description of the map $D$ is given in Proposition~\ref{prop:D} of 
Section~\ref{sec:D}.   
We emphasize the fact that, in general, the linearized contact 
homology groups and the map $D$ depend on the filling $(W,\om)$.  
The example of Riemann surfaces in Section~\ref{sec:Riemann} shows that this 
is the case even in the simple situation $M=S^1$. On the other hand,  
there are classes of contact manifolds for which the linearized contact 
homology groups and the map $D$ only depend on $(M,\xi)$. 
This is illustrated in Section~\ref{sec:subcritical} by subcritically Stein 
fillable contact manifolds $(M,\xi)$ of dimension $\ge 3$ with $c_1(\xi) = 0$.

\medskip  
 
Besides the Floer/SFT distinction mentioned above, the groups $SH_*$ and 
$HC_*$ differ in a more subtle way, related    
to the fact that the loop space naturally carries  
an $S^1$-action. The construction of contact homology groups is 
intrinsically $S^1$-equivariant in the sense that the generators of 
the complex are 
\emph{unparametrized} Reeb orbits and $S^1$ acts on the relevant  
spaces of solutions, whereas the construction of symplectic homology groups is 
non-equivariant, i.e. the Hamiltonian is time-dependent and the 
generators of the complex are \emph{parametrized} Hamiltonian orbits.  
 
One should recall at this point the Gysin exact sequence relating ordinary 
and $S^1$-equivariant homology of an $S^1$-space $X$, which reads  
\begin{equation} \label{eq:Gysin-intro}  
\ldots \longrightarrow H_*(X) \longrightarrow H_*^{S^1}(X) \stackrel 
{\cap e} \longrightarrow H_{*-2}^{S^1}(X) \longrightarrow H_{*-1}(X)  
\longrightarrow \ldots 
\end{equation}    
The analogy -- modulo shifts in the grading -- 
between~\eqref{eq:intro} and~\eqref{eq:Gysin-intro} is by no means 
formal. We prove in~\cite{BO} that an $S^1$-equivariant version of 
symplectic homology is isomorphic to linearized contact homology and 
that~\eqref{eq:intro} is the corresponding Gysin exact sequence,  
whereas the paper~\cite{CO} constructs a non-equivariant version of 
(linearized) contact homology fitting into a Gysin exact sequence with 
the usual contact homology groups.  
 
>From this point of view, contact- and symplectic 
homology are closely related complementary theories, linked via a Gysin exact 
sequence. This perspective on symplectic homology also relates to a 
recent 
conjecture of Seidel~\cite{Se} predicting that symplectic homology is 
isomorphic to the Hochschild homology $HH_*(\cC)$ of a suitable 
$A_\infty$-category $\cC$. Then $S^1$-equivariant symplectic homology 
should be isomorphic to the cyclic homology $HC_*(\cC)$  and the 
Gysin exact sequences mentioned above should be isomorphic to the 
standard Connes exact sequence connecting Hochschild and cyclic 
homology  
$$ 
\ldots \longrightarrow HH_*(\cC) \longrightarrow HC_*(\cC) \longrightarrow 
HC_{*-2}(\cC) \longrightarrow HH_{*-1}(\cC) \longrightarrow \ldots 
$$ 
 
As far as the connecting map $D$ is concerned, the analogy 
with the finite dimensional case is again fertile. It can be thought of as 
a cap product with an Euler class, just as the map 
$H_*^{S^1}(X)\to H_{*-2}^{S^1}(X)$ in~\eqref{eq:Gysin-intro} is the cap 
product with the Euler class of the $S^1$-bundle over the homotopy 
quotient $X \times_{S^1} ES^1$.  
 
\medskip 
 
This paper essentially consists of a proof of Theorem~\ref{thm:intro} 
and we now give an overview of the proof.      
We draw the reader's attention to Section~\ref{sec:MBcomplex}, where we have 
concentrated the key statements in rigorous  
form. The preliminary constructions are given  
in Sections~\ref{sec:symplectic} and~\ref{sec:bigcontact}. 
 
The main technical tool for our proof is the Morse-Bott complex developed 
in~\cite{BOauto} following ideas from~\cite{B}. The construction, 
which is summarized  
in Section~\ref{sec:MBsymp}, gives a recipe to compute the symplectic homology 
groups in terms of the moduli spaces of Floer trajectories for a 
\emph{time-independent} Hamiltonian $H$ under the assumption -- 
generic for autonomous Hamiltonians -- that the 
$1$-periodic orbits of the latter are either constant and 
nondegenerate, or nonconstant and transversally nondegenerate. The 
Morse-Bott complex mimicks a time-dependent perturbation of $H$ via the 
choice of a perfect Morse function $f_\gamma$ along the geometric image of 
each nonconstant orbit $\gamma$, which is a circle. Each 
(unparametrized) nonconstant orbit $\gamma$ gives rise to two 
generators in the Morse-Bott complex, one for the minimum and one for 
the maximum of $f_\gamma$.  
 
For the special type of autonomous Hamiltonians used to define symplectic 
homology the nonconstant orbits $\gamma$ are in one-to-one correspondence with 
closed Reeb orbits $\gamma'$ on $M$. A nice feature of the     
Morse-Bott complex is that it is naturally filtered by the Maslov 
index $\mu(\gamma')$, and this filtration gives rise to a spectral 
sequence supported in two lines. As seen in Section~\ref{sec:proofleq}, any 
such spectral sequence gives rise to a long exact sequence of the type 
$$ 
\ldots \longrightarrow E^\infty \longrightarrow E^2_{k,0} 
\stackrel {d^2} \longrightarrow E^2_{k-2,1} \longrightarrow 
E^\infty \longrightarrow \ldots 
$$ 
By definition we have $E^\infty\simeq SH$. In order to establish  
Theorem~\ref{thm:intro} we prove that $E^2_{k,i} \simeq HC_{k+(n-3)}$,  
$i=0,1$ and identify the differential $d^2$ in the following way.    
In Section~\ref{sec:nonequiv} we define a 
non-equivariant version of contact homology, inspired by~\cite{CO},   
by means of a construction of a Morse-Bott complex which we 
call the \emph{$S^1$-parametrized contact complex},    
analogous to the one in Section~\ref{sec:MBsymp}. This complex is also  
filtered by $\mu(\gamma')$ and the $E^2$-term of the associated 
spectral sequence is trivially identified with $HC_*$. We obtain thus 
two filtered complexes, one for symplectic homology and another one 
for contact homology, which we prove in Section~\ref{sec:filtered} to be 
isomorphic. This automatically implies that the $E^2$-terms of the 
corresponding spectral sequences are isomorphic, and also that the 
corresponding $d^2$-differentials coincide.  
 
The isomorphism between the two filtered complexes is established by 
considering ``mixed'' moduli spaces consisting of punctured curves  
defined on the cylinder $\R\times S^1$ and taking values in the 
symplectization $M\times \R$. Near $-\infty$, these curves are holomorphic  
and asymptotic to a Reeb orbit; near $+\infty$, they satisfy Floer's 
equation and are asymptotic to a $1$-periodic orbit of $H$.      
Since the contact action decreases along such curves the 
resulting chain map has upper triangular form, and we show, by  
constructing solutions to the mixed problem described above and showing  
they are unique, that      
the entries on the diagonal are $\pm 1$. This method of establishing 
an isomorphism at the chain level    
by using mixed moduli spaces is reminiscent of~\cite{AbboSch}.  
 
We note at this point the fact that the Morse-Bott construction  
of the $S^1$-parametrized contact complex is necessary only in order 
to identify the differential $d^2$ in the exact sequence. The 
isomorphism between the $E^1$-term of the symplectic homology spectral 
sequence and the linearized contact complex can be established 
directly by using the mixed moduli spaces described above.   
 
Two more remarks are in order. The first one concerns the fact that 
the Floer trajectories for $H$ might wander deep inside the filling 
$W$, whereas the isomorphism between the two filtered complexes is 
constructed in the  
symplectization $M\times \R$. The basic technique 
in Section~\ref{sec:stretch} is to stretch the 
neck near the boundary of $W$ and show that, when the stretching 
parameter is large enough, the Floer trajectories in $W$ are in bijective 
correspondence with punctured Floer trajectories in $M\times \R$, 
capped at the punctures with rigid holomorphic planes in the 
symplectization of $W$. The maximum principle plays a crucial 
role for showing that the limit building contains no curve with more 
than one positive puncture.  
 
The second remark concerns the problem of good and bad orbits 
(see Section~\ref{sec:contact} for the definition). One of the most pleasant 
features of the spectral sequences described above is that, although 
the starting complex contains two generators  
for each Reeb orbit, only the generators corresponding to good Reeb 
orbits survive to $E^1$. At this point it is crucial to use 
$\Q$-coefficients for $\Lambda_\om$ rather than 
$\Z$-coefficients. Besides the analysis of signs borrowed  
from~\cite{BOauto}, we need to show that, for a suitable choice of the  
Hamiltonian $H$, certain rigid Floer 
trajectories do not appear and hence do not    
contribute to the expression of $d^0$ in the case of symplectic 
homology (see Section~\ref{sec:MBcomplex} for details). This is done 
in Section~\ref{sec:slow} by slowing down the rate of variation of $H$ and 
using the regularity assumptions on the time-independent almost 
complex structure $J$ for contact homology. 
 
The paper ends with Section~\ref{sec:examples} in which we treat four 
examples: Riemann surfaces with one boundary component, subcritical 
Stein domains, negative disc bundles and unit cotangent bundles. 
 
\medskip  
 
\begin{figure}[hpt] 
\input{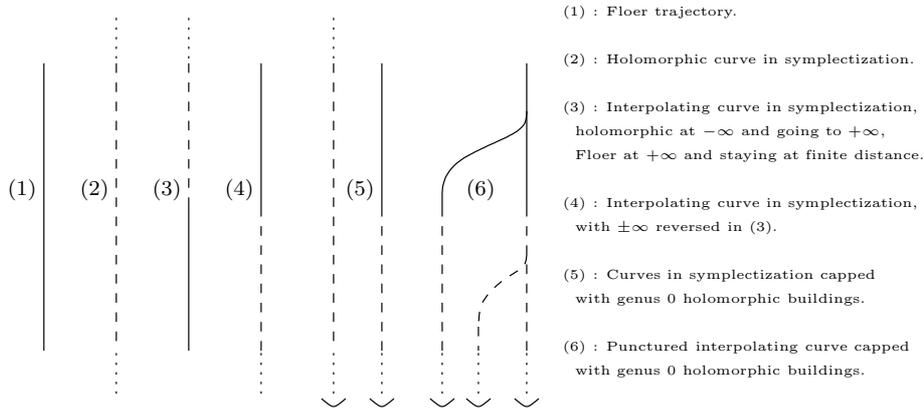} 
\caption{Pictorial conventions \label{fig:legend}} 
\end{figure}

\noindent {\bf Note on pictorial conventions.} 
  We use several different types of moduli spaces in the paper, 
  the most important of which are shown in Figure~\ref{fig:Phi}  
  on page~\pageref{fig:Phi}. There and throughout the paper we  
  use the following conventions (cf. Figure~\ref{fig:legend}):   
  gradient trajectories of Morse functions are represented by 
  horizontal lines, solutions of Floer's equation by vertical lines, 
  and holomorphic curves in a symplectization by dashed vertical 
  lines. Vertical dots stand for holomorphic curves in 
  a symplectization going to $\pm\infty$ at a puncture.  
 
\medskip  
 
\noindent {\it Acknowledgements.} We thank Otto van 
Koert and Paul Seidel for useful suggestions, and also the  
anonymous  
referees for their careful reading of the manuscript and their 
precise comments and remarks, which helped us clarify and improve the 
presentation.   
The non-equivariant construction 
of Section~\ref{sec:nonequiv} and the subsequent treatment of the 
isomorphism between the $E^1$-term of the spectral sequence and the 
contact complex  
have been influenced by discussions with Kai Cieliebak.  
Both authors acknowledge financial support from Fonds National de 
la Recherche Scientifique (Belgium), Institut 
de Recherche Ma\-th\'e\-ma\-tique Avanc\'ee (Strasbourg), Forschungsinstitut 
f\"ur Mathematik (Z\"urich), and the Mathematisches 
Forschungsinstitut (Oberwolfach).

 
\section{Symplectic homology} \label{sec:symplectic} 
 
\subsection{Construction and basic properties of Symplectic homology}  
\label{sec:recollections}  
 
Let $(W,\om)$ be a compact symplectic manifold with   
boundary $M:=\p W$ of contact type, satisfying~\eqref{eq:asph}. The 
Liouville vector field is denoted by $X$, the induced contact form on 
$M$ is $\lambda$, the contact distribution is $\xi$ and 
the Reeb vector field is $R_\lambda$.   
 
We now construct the symplectic homology groups of $(W,\omega)$. 
For more details, we refer to \cite[Section~2]{BOauto}.    
Let $\phi$ be the flow of $X$. We parametrize a neighbourhood $U$ of $M$ by 
$G: M \times [-\delta, 0] \to U$, $(p,t) \mapsto \phi^t(p)$. 
Then $d(e^t\lambda)$ is a symplectic form on $M\times \R^+$ and 
$G^*\om = d(e^t \lambda)$. 
The {\bf symplectic completion} $(\widehat W, \widehat \omega)$ is 
defined by 
$$ 
\widehat W : = W \ \bigcup _{G} \ M\times \R^+,\qquad  
\widehat \om : = 
\left\{\begin{array}{ll} 
\om, & \textrm{ on } W, \\ 
d(e^t \lambda), & \textrm{ on } M\times \R^+. 
\end{array} \right. 
$$ 
 
Given a time-dependent Hamiltonian $H :S^1\times \widehat W \to \R$, 
we define the 
{\bf Hamiltonian vector field} $X^\theta_H$ by 
$$ 
\widehat \om (X^\theta_H,\cdot) = d H_\theta, \qquad \theta\in S^1 = \R/\Z, 
$$ 
where $H_\theta:=H(\theta,\cdot)$. We denote by $\phi_H$ the flow of 
$X_H^\theta$, defined by $\phi_H^0=\textrm{Id}$ and 
$$ 
  \frac d {d\theta} \phi_H^\theta (x) = X^\theta_H(\phi_H^\theta(x)),  
\qquad \theta\in \R. 
$$ 
 
We denote by $\cP(H)$ the set of $1$-periodic orbits of $X^\theta_H$.   
Given a free homotopy class $a$ of loops in $W$, we denote by 
$\cP^a(H)$ the set of $1$-periodic orbits of $X^\theta_H$ in the class 
$a$. The {\bf symplectic homology groups} of $(W,\om)$ are defined as 
the direct limit 
\begin{equation*} 
   SH_*^a(W,\om) := \lim _{\substack{ \to \\ H\in \cH} 
    } SH_*^a(H). 
\end{equation*} 
Here $\cH$ is a suitable class of Hamiltonians and $SH_*^a(H)$ are 
the Floer homology groups of $H$ in the class $a$. The underlying 
complex $SC_*^a(H,J)$ is generated by the elements of $\cP^a(H)$, and 
the differential $d$ is defined using a time-dependent almost 
complex structure $J$ on $\widehat W$ which is regular for $H$ and has 
a special behaviour at infinity. The resulting homology groups do 
not depend on $J$ and we omit it from the notation.  
 
Let now $a=0$ be the trivial homotopy class. In this case we denote 
the symplectic homology groups by $SH_*(W,\om)$. We define the {\bf 
reduced Hamiltonian action functional} 
$$ 
\cA_H^0 : C^\infty_{\textrm{contr}}(S^1,\widehat W) \to \R 
$$ 
by 
$$ 
\cA_H^0(\gamma) := -\int_{D^2} \sigma^*\widehat \om - \int_{S^1} 
H(\theta,\gamma(\theta)) \, d\theta. 
$$ 
Here $C^\infty_{\textrm{contr}}(S^1,\widehat W)$ denotes the space of 
smooth contractible loops in $\widehat W$ and $\sigma:D^2\to \widehat 
W$ is a smooth extension of $\gamma$. Note that $\cA_H^0$ is 
well-defined thanks to condition~(\ref{eq:asph}) and that   
$\cA_H^0$ is  
decreasing along Floer trajectories. We define the {\bf action spectrum} of  
$(M,\lambda)$ by     
$$ 
\textrm{Spec}(M,\lambda) := \{ T \in \R^+\, | \, \textrm{ there is a 
   closed } R_\lambda\textrm{-orbit of period } T\}. 
$$ 
We fix $\epsilon>0$ such that 
$\epsilon < T$ for all $T\in \textrm{Spec}(M,\lambda)$. For a regular 
almost complex structure $J$ we define the chain complexes 
\begin{equation} \label{eq:SC-} 
SC_*^-(H,J) := \bigoplus _{\substack{ \gamma \in 
     \cP^0(H) \\ \cA_H^0(\gamma) \le \epsilon }} \Lambda_\om 
\langle \gamma \rangle \ \subset SC_*(H,J) 
\end{equation} 
and 
$$ 
SC_*^+(H,J) := SC_*(H) / SC_*^-(H,J). 
$$ 
The differential on $SC_*^\pm(H,J)$ is induced by $d$. The groups 
$$ 
SH_*^\pm(H) := H_*(SC_*^\pm(H,J),d) 
$$ 
neither depend on $J$ nor on $\epsilon$.     
We define 
$$ 
SH_*^\pm(W,\om) := \lim _{\substack{ \to \\ H\in \cH}} SH_*^\pm(H). 
$$ 
We call $SH_*^+(W,\om)$ the {\bf positive symplectic homology group} 
of $(W,\om)$. 
   
\begin{remark} {\rm 
For contractible orbits condition~\eqref{eq:asph} can be  
replaced by the weaker {\bf symplectic asphericity} condition  
$\langle \om,\pi_2(W)\rangle =0$. } 
\end{remark} 
 
Let us assume that $M$ is of {\bf positive contact 
   type}~\cite[\S4.3]{O1}. This means that every  
closed $R_\lambda$-orbit $\gamma$ on $M$ which is contractible in $W$  
has positive action $\cA_\om(\gamma)$ bounded away from zero, where 
$$ 
\cA_\om(\gamma) := \int_{D^2} \sigma^*\om 
$$ 
for some extension $\sigma:D^2\to W$ of $\gamma$. This condition is 
automatically satisfied if the boundary $M$ is of restricted contact 
type, i.e. the vector field $X$ is globally defined on $W$. Under the 
positive contact type assumption we 
have~\cite[Proposition~1.4]{V}~(see also~\cite{O1})  
$$ 
SH_*^-(W,\om) = H_{*+n}(W,\p W;\Lambda_\om), \qquad n=\frac 1 2 \dim 
\, W. 
$$ 
 
Moreover, the short exact sequence $SC_*^-(H) \to SC_*(H) 
\to SC_*^+(H)$ induces the long exact sequence~\cite{V} 
\begin{equation}{\scriptstyle 
   \label{eq:Vseq} 
\ldots \to SH_{*+1}^+(W,\om) \!\to H_{*+n}(W,\p W;\Lambda_\om) \!\to 
SH_*(W,\om) \!\to SH_*^+(W,\om) \to \ldots} 
\end{equation} 
 
\subsection{Morse-Bott description of symplectic homology} 
\label{sec:MBsymp} 
 
In this section we recall the Morse-Bott formalism  
of~\cite{BOauto}.  We assume in this section that the closed 
$R_\lambda$-orbits are transversally nondegenerate in $M$. 
We denote by $\phi_\lambda$ the flow of $R_\lambda$.  
 
In~\cite[\S3]{BOauto} we used a class $\cH'$ of admissible 
Hamiltonians consisting of elements $H:\widehat W\to \R$ such that 
\begin{enum} 
\item $H|_W$ is a $C^2$-small Morse function and $H<0$ on $W$; 
\item $H(p,t)=h(t)$ outside $W$, where $h(t)$ is a strictly 
   increasing function with $h(t)=\alpha e^t+\beta$, 
   $\alpha,\beta\in\R$, $\alpha \notin \textrm{Spec}(M,\lambda)$ for 
   $t$ bigger than some $t_0$, and such that $h''-h'>0$ on $[0,t_0[$. 
\end{enum} 
 
Note that the $1$-periodic orbits of $X_H$ in $W$ are constant and 
nondegenerate by assumption~(i). A direct computation shows that 
\begin{equation} \label{eq:XHR} 
X_H(p,t) = -e^{-t}h'(t) R_\lambda, \qquad \textrm{for }  
(p,t) \in M \times [0,\infty). 
\end{equation}   
 
The $1$-periodic orbits of $X_H$ fall in two classes: 
\begin{enumerate} 
\item[(1)] critical points of $H$ in $W$; 
\item[(2)] nonconstant $1$-periodic orbits of $X_h$, located on levels 
   $M\times \{t\}$, $t\in ]0,t_0[$, which are in one-to-one 
   correspondence with closed $-R_\lambda$-orbits of period 
   $e^{-t}h'(t)$. 
\end{enumerate} 
 
Let $\alpha := \lim_{t \to \infty} e^{-t} H(p,t)$. Let 
$\cP_\lambda$ be the set of closed unparametrized $R_\lambda$-orbits in $M$.  
We denote by $\cP_\lambda^{\le \alpha}$ the set of all $\gamma' \in 
\cP_\lambda$  
such that $\cA_\lambda(\gamma') \le \alpha$. 
Because $H$ is independent of $\theta$, every 
   orbit $\gamma'\in\cP_\lambda^{\le \alpha}$ gives rise to a whole 
circle of nonconstant $1$-periodic orbits $\gamma$ of $X_H$, which are  
   transversally nondegenerate and whose parametrizations differ by a 
shift $\theta\in S^1$. We denote by $S_\gamma$ the set of such orbits, 
so that $S_\gamma=S_{\gamma(\cdot + \theta)}$ for all $\theta\in S^1$. 
 
Let $\cJ=\cJ(\widehat W,\widehat \omega)$ be the space of 
$\theta$-dependent almost complex structures $J$ such that 
\begin{enum} 
\item $J$ is compatible with $\widehat \omega$; 
\item for $t$ large enough, $J$ is independent of $\theta$; 
\item $J$ preserves the contact distribution $\xi$; 
\item  $J\frac \partial {\partial t}=R_\lambda$.  
\end{enum} 
Given $\og,\ug\in\cP(H)$, $\tq\in \textrm{Crit}(H)$ and $J \in \cJ$, 
we denote by  
$$ 
\widehat \cM^A(S_{\og},S_{\ug};H,J), \qquad 
\widehat \cM^A(S_{\og},\tq;H,J) 
$$ 
the spaces of Floer 
trajectories for $(H,J)$ starting at $S_{\og}$ and ending  
at $S_{\ug}$ or $\tq$, respectively.  
Such a Floer trajectory is a map $u : \R \times S^1 \to \widehat W$   
satisfying 
\begin{equation} \label{eq:Floereq} 
\p_s u + J (\p_\theta u - X_H) = 0 \ \ \mbox{for all} \ \  
(s,\theta)\in\R\times S^1, 
\end{equation} 
as well as the conditions   
\begin{equation} \label{eq:asyog} 
\lim_{s\to-\infty} u(s,\cdot) \in S_\og 
\end{equation} 
and, respectively,  
\begin{equation} \label{eq:asyug} 
\lim_{s\to\infty}u(s,\cdot)\in S_\ug \ \ \mbox{or} \ \  
\lim_{s\to\infty}u(s,\cdot)=\tq. 
\end{equation} 
 
The {\bf Morse-Bott moduli spaces of Floer trajectories} 
are defined by  
$$ 
\cM^A(S_{\og},S_{\ug};H,J):=\widehat \cM^A(S_{\og},S_{\ug};H,J)/\R 
$$ 
and 
$$ 
\cM^A(S_{\og},\tq;H,J):=\widehat \cM^A(S_{\og},\tq;H,J)/\R. 
$$ 
Let $\Jreg(H)\subset \cJ$ be the set of those almost complex 
structures for which the linearization of the  
equation~\eqref{eq:Floereq} at its solutions is surjective, so that 
$\Jreg(H)$ is dense in $\cJ$~\cite[Proposition 3.5]{BOauto}. Given 
$J\in\Jreg(H)$ the Morse-Bott moduli spaces of Floer trajectories are 
smooth manifolds, and their respective dimensions are 
\begin{eqnarray} \label{eq:dim} 
\dim \, \cM^A(S_{\og},S_{\ug};H,J) & = & \mu(\og) - \mu(\ug) + 2\langle 
c_1(TW),A\rangle,\\ 
\dim \, \cM^A(S_{\og},\tq;H,J) & = & \mu(\og) - \mu(\tq) + 2\langle 
c_1(TW),A\rangle. \nonumber 
\end{eqnarray} 
Here $\mu(\tq)=\ind(\tq;-H)-n$ is the Conley-Zehnder index of the 
constant orbit $\tq$, whereas $\mu(\og),\mu(\ug)$ denote the 
Conley-Zehnder indexes of the linearized Hamiltonian flows restricted 
to $\xi$. 
  
We have natural evaluation maps 
$$ 
\oev : \cM^A(S_{\og},S_{\ug};H,J) \to S_{\og}, \qquad 
\uev : \cM^A(S_{\og},S_{\ug};H,J) \to S_{\ug} 
$$ 
and 
$$ 
\oev : \cM^A(S_{\og},\tq;H,J) \to S_{\og} 
$$ 
defined by 
$$ 
\oev([u]):= \lim_{s\to-\infty} u(s,\cdot), \qquad \uev([u]):= 
\lim_{s\to \infty} u(s,\cdot). 
$$ 
 
For each $S_\gamma$, $\gamma\in \cP(H)$ we choose a Morse function 
$f_{S_\gamma}:S_\gamma\to \R$ with exactly one maximum $M$ and one 
minimum $m$. To simplify notation, we shall write in the sequel 
$f_\gamma$ instead of $f_{S_\gamma}$, so that 
$f_\gamma=f_{\gamma(\cdot+\theta)}$ for all $\theta\in S^1$.  
We denote by $\gamma_m$, $\gamma_M$ the orbits in 
$S_\gamma$ starting at the minimum and the maximum of 
$f_\gamma$ respectively. 
For a generic choice of these Morse functions \cite[Lemma 3.6]{BOauto}, 
all the maps $\oev$ are transverse to the unstable manifolds $W^u(p)$, $p\in 
\textrm{Crit}(f_\gamma)$, all the maps $\uev$ are transverse to the 
stable manifolds $W^s(p)$, $p\in \textrm{Crit}(f_\gamma)$ and all 
pairs 
\begin{equation} \label{eq:transvf} 
(\oev,\uev) : \cM^A(S_{\og},S_{\ug};H,J)   \to  S_{\og} \times 
S_{\ug},  
\end{equation}  
$$ 
(\oev,\uev) : \cM^{A_1}(S_{\og},S_{\gamma_1};H,J) \ _{\uev}\times_{\oev}  
\cM^{A_2}(S_{\gamma_1},S_{\ug};H,J)   \to  S_{\og} \times 
S_{\ug} 
$$ 
are transverse to products $W^u(p)\times W^s(q)$, $p\in 
\textrm{Crit}(f_{\og})$, $q\in \textrm{Crit}(f_{\ug})$. 
The unstable and stable manifolds are 
understood with respect to $\nabla f_\gamma$, so that 
$W^u(M) = \{ M \}$, $W^s(M) = S_\gamma \setminus \{ m \}$,  
$W^u(m) = S_\gamma \setminus \{ M \}$ and $W^s(m) = \{ m \}$.   
We denote by 
$\Freg(H,J)$ the set consisting of collections $\{f_\gamma\}$ of Morse 
functions which satisfy the above transversality conditions.  
 
Let now $J\in \Jreg(H)$ and $\{f_\gamma\}\in \Freg(H,J)$. For $p\in 
\textrm{Crit}(f_\gamma)$ we denote the Morse index by  
$$\ind(p):=\dim \, W^u(p;\nabla f_\gamma). 
\index{$\ind(p)$}   
$$ 
Let $\og,\ug\in\cP(H)$ and $p\in 
\mathrm{Crit}(f_{\og})$, $q\in \mathrm{Crit}(f_{\ug})$. For $m\ge 0$ 
we denote by  
\begin{equation} \label{eq:bigM} 
\cM_m^A(p,q;H,\{f_\gamma\},J) 
\index{$\cM_m^A(p,q;H,\{f_\gamma\},J)$, $\cM^A(p,q;H,\{f_\gamma\},J)$|(} 
\end{equation} 
the union for $\tgamma_1,\ldots,\tgamma_{m-1}\in\cP(H)$ and 
$A_1+\ldots +A_m=A$ of the fibered products  
\begin{eqnarray*} 
&& 
\hspace{-.7cm}W^u(p)  
\times_{\oev} 
(\cM^{A_1}(S_{\og}\,,S_{\tgamma_1})\!\times\!\R^+) 
{_{\varphi_{f_{\tgamma_1}}\!\circ\uev}}\!\times   
_{\oev} 
(\cM^{A_2}(S_{\tgamma_1},S_{\tgamma_2})\!\times\!\R^+) \\ 
&&  
{_{\varphi_{f_{\tgamma_2}}\!\circ\uev}\times_{\oev}} \ldots\, 
{_{\varphi_{f_{\tgamma_{m-1}}}\!\!\circ\uev}}\!\!\times 
_{\oev}   
\cM^{A_m}(S_{\tgamma_{m-1}},\!S_{\ug})  
{_{\uev}\times} W^s(q). 
\end{eqnarray*} 
This is a smooth manifold of dimension 
\begin{eqnarray} \label{eq:dimension} 
\lefteqn{\dim \, \cM_m^A(p,q;H,\{f_\gamma\},J)} \\ 
& = & \ \mu(\og) + \ind(p) - \mu(\ug) - \ind(q) + 2\langle 
c_1(TW),A\rangle - 1 , 
\end{eqnarray} 
as shown in~\cite{BOauto}.     
Note that the moduli space $\cM_0^A(p,q;H,\{f_\gamma\},J)$ is a
submanifold of $\cM^A(S_{\og},S_{\ug};H,J)$. We denote  
$$ 
\cM^A(p,q;H,\{f_\gamma\},J):=\bigcup_{m\ge 0} 
\cM_m^A(p,q;H,\{f_\gamma\},J) 
\index{$\cM_m^A(p,q;H,\{f_\gamma\},J)$, $\cM^A(p,q;H,\{f_\gamma\},J)$|)} 
$$ 
and we call this {\bf the moduli space of Morse-Bott broken 
trajectories}, whereas $\cM_m^A(p,q;H,\{f_\gamma\},J)$ is called {\bf 
the moduli space of Morse-Bott broken trajectories with $m$ 
sublevels}. We refer to Figure~\ref{fig:Phi}.(a) on 
page~\pageref{fig:Phi} for a visual 
representation of the elements of these moduli spaces.

In the sequel we use only $0$-dimensional moduli spaces 
$\cM^A(p,q;H,\{f_\gamma\},J)$. Nevertheless, we give now a brief 
description of the topology of the compactification 
$\overline\cM^A(p,q;H,\{f_\gamma\},J)$ in the general case.   
We first start with the compactification $\overline 
\cM^A_m(p,q;H,\{f_\gamma\},J)$. There are three  
types of codimension $1$ degeneracies for sequences $\u_n$ 
of elements in $\cM^A_m(p,q;H,\{f_\gamma\},J)$. Firstly, one of the $m$ Floer 
trajectories composing $\u_n$ can break in two Floer 
trajectories as $n\to\infty$. Secondly, the flow time of one of the 
$m-1$ finite  
gradient trajectories can shrink to $0$. Thirdly, one of the $m+1$ 
gradient trajectories can break in two gradient trajectories. Higher 
codimension degeneracies of sequences of elements of 
$\cM^A_m(p,q;H,\{f_\gamma\},J)$ are obtained by combining the above 
three types of degeneracies. The space $\overline \cM^A(p,q;H,\{f_\gamma\},J)$ 
is obtained by gluing the moduli spaces 
$\overline \cM^A_m(p,q;H,\{f_\gamma\},J)$ along their common boundary strata 
involving only degenerations of the first two types. For simplicity, we 
just describe the case of codimension $1$ boundary strata:  
these are common boundary strata for 
$\overline \cM^A_{m-1}(p,q;H,\{f_\gamma\},J)$ and  
$\overline \cM^A_m(p,q;H,\{f_\gamma\},J)$ and correspond to degenerations of 
the first and second type respectively. The boundary (and the corners) 
of $\overline \cM^A(p,q;H,\{f_\gamma\},J)$ correspond to at least one degeneracy 
of the third type.  
 
Given $\og\in\cP(H)$, 
$p\in \mathrm{Crit}(f_{\og})$, $\tq\in \mathrm{Crit}(H)$, we define 
the moduli spaces $\cM_m^A(p,\tq;H,\{f_\gamma\},J)$, $m\ge 0$  
of Morse-Bott broken trajectories by replacing the last  
term $\cM^{A_m}(S_{\tgamma_{m-1}},S_\ug) {_{\uev}\times} W^s(q)$ 
in~\eqref{eq:bigM} with $\cM^{A_m}(S_{\tgamma_{m-1}},\tq;H,J)$. The 
union over $m\ge 0$ of these spaces is denoted by 
$\cM^A(p,\tq;H,\{f_\gamma\},J)$.  This is 
again well defined as a smooth manifold of dimension  
\begin{eqnarray*} 
\lefteqn{\dim \, \cM^A(p,\tq;H,\{f_\gamma\},J)} \\ 
& = & \ \mu(\og) + \ind(p) - \ind(\tq;-H)  
+ n + 2 \langle c_1(TW),A\rangle -1. 
\end{eqnarray*} 
Again, $\cM_0^A(p,\tq;H,\{f_\gamma\},J)$ is a submanifold of 
the space  
$\cM^A(S_{\og},\tq;H,J)$.  
 
\begin{remark} \label{rmk:grad_traj} \rm 
Since $H$ is $C^2$-small, the moduli spaces 
$\cM^A(\tp,\tq;H,J)$, $\tp,\tq\in \mathrm{Crit}(H)$ of expected 
dimension  
$$ 
\ind(\tp;-H)-\ind(\tq;-H) + 2\langle c_1(TW),A)\rangle -1=0 
$$ 
consist exclusively of gradient trajectories of $H$ in 
$W$~\cite[Theorem~6.1]{HS}(see also~\cite[Theorem~7.3]{SZ}). As a 
consequence, these moduli spaces are empty whenever $A\neq 0$. 
\end{remark} 
 
For each  
$[u]\in \cM^A(p,q;H,\{f_\gamma\},J)$ or  
$[u]\in \cM^A(p,\tq;H,\{f_\gamma\},J)$ we have defined 
in~\cite{BOauto} a sign $\bar \epsilon(u)$. Let $a$ be a free homotopy 
class of loops in $W$.  
We define the {\bf Morse-Bott chain groups} by 
\begin{eqnarray}  \label{eq:MBchain_a} 
   BC_*^a(H) & := & 
\bigoplus _{  \gamma\in \cP^a(H)} 
\Lambda_\om \langle \gamma_m,\gamma_M\rangle, \qquad a\neq 0, \\ 
BC_*^0(H) & := & \bigoplus_{\scriptstyle \tp\in \mathrm{Crit}(H)} 
   \Lambda_\om \langle \tp \rangle \ \oplus 
\bigoplus _{ \gamma\in \cP^0(H) } 
\Lambda_\om \langle \gamma_m,\gamma_M\rangle. 
\end{eqnarray} 
The grading is defined by 
\begin{eqnarray*} 
|e^A \tp| & := & \ind(\tp;-H) - n - 2\langle c_1(TW),A\rangle, \\ 
|e^A \gamma_m | & := &  \mu(\gamma) + 1 - 2\langle c_1(TW),A\rangle, 
\\ 
|e^A \gamma_M | & := &  \mu(\gamma) - 2\langle c_1(TW),A\rangle. 
\end{eqnarray*} 
 
We define the {\bf Morse-Bott differential} 
$$ 
d : BC_*^a(H) \to BC_{*-1}^a(H) 
$$ 
by 
\begin{eqnarray} \label{eq:MBdiff_p} 
   d \tp & := & \sum_{\substack{ \tq\in \mathrm{Crit}(H) \\ 
      |\tp|-|\tq| =1}} 
\ \sum_{\scriptstyle [u]\in \cM^0(\tp,\tq;H,\{f_\gamma\},J)} \bar 
   \epsilon(u)\tq, \\  
\label{eq:MBdiff_gamma} 
d \gamma_p & := & \sum_{\substack{ 
  \tq\in \mathrm{Crit}(H) \\ 
|\gamma_p| - |e^A \tq|=1}} 
\ \sum_{\scriptstyle [u]\in \cM^A(p,\tq;H,\{f_\gamma\},J)} 
\bar \epsilon(u)e^A \tq \\ 
& & + \sum_{\substack{ 
  \ug\in \cP(H), q\in \mathrm{Crit}(f_{\ug}) \\ 
|\gamma_p| - |e^A \ug_q|=1}} 
\ \sum_{\scriptstyle [u]\in \cM^A(p,q;H,\{f_\gamma\},J)} 
\bar \epsilon(u)e^A \ug_q, \ p \in {\rm Crit}(f_\gamma). \nonumber 
\end{eqnarray} 
The sums~(\ref{eq:MBdiff_p}) and~(\ref{eq:MBdiff_gamma}) clearly 
involve only periodic orbits in the same free homotopy class as that 
of $\tp$ or $\gamma_p$ respectively. 
 
The Correspondence Theorem 3.7 in~\cite{BOauto} 
implies $d\circ d=0$ and 
$$ 
\lim_{\substack{ \to \\ H\in \cH'}} 
H_*(BC_*^a(H),d) = SH_*^a(W,\om). 
$$ 
We shall denote in the sequel  
$$ 
SH_*^a(H,J):= H_*(BC_*^a(H),d). 
$$ 
Moreover, if we define the subcomplex 
$$ 
BC_*^-(H) := \bigoplus_{\scriptstyle \tp\in \mathrm{Crit}(H)} 
   \Lambda_\om \langle \tp \rangle 
$$ 
and the quotient 
$$ 
BC_*^+(H):= BC_*^0(H)/BC_*^-(H), 
$$ 
we have 
$$ 
\lim_{\substack{ \to \\ H\in \cH'}} H_*(BC_*^+(H),d) = SH_*^+(W,\om). 
$$ 
 
\noindent {\bf Notation.} From now on the letter $a$ will denote 
either a free homotopy class in $W$ or one of the symbols $\pm$. The 
notation $i^{-1}(+)$ stands for $i^{-1}(0)$. 
 
\medskip  
 
The previous description 
of the Floer differential can be 
generalized to the case of an $s$-dependent family of autonomous 
Hamiltonians. More precisely, let $H_s$, $s\in\R$ be 
a homotopy of autonomous Hamiltonians satisfying the following 
conditions:   
\begin{enum} 
\item $H_s$ is increasing with respect to $s$, and constant for $|s|$ 
  large enough;  
\item $H_s(p,t)=h_s(t)$ outside $W$, where $h_s$ is a strictly 
increasing function with $h_s(t)=\alpha(s) e^t+\beta(s)$ for $t$ 
bigger than some $t_0$; 
\item $H_\pm:=\lim_{s\to \pm\infty}H_s$ belong to the class $\cH'$ of 
admissible Hamiltonians. 
\end{enum} 
We call this an {\bf admissible homotopy of Hamiltonians}.  
Similarly, we define an 
{\bf admissible homotopy of almost complex structures} to be a family 
$J_s$, $s\in\R$ of elements of $\cJ$, which is constant for $|s|$ 
large enough.  
 
Given $A\in H_2(W,\Z)$, $\og\in\cP(H_-)$, $\ug\in\cP(H_+)$, and $\tq\in 
\mathrm{Crit}(H_+)$, we denote by  
$$ 
\cM^A(S_{\og},S_{\ug};H_s,J_s), \qquad 
\cM^A(S_{\og},\tq;H_s,J_s) 
$$ 
the spaces of Floer 
trajectories $u:\R\times S^1\to \widehat W$ satisfying  
\begin{equation} \label{eq:sdep} 
\partial _s u + J_s\partial _\theta u = X_{H_s}(u) \ \ \mbox{for all} \ \  
(s,\theta)\in\R\times S^1, 
\end{equation} 
as well as conditions~\eqref{eq:asyog} and~\eqref{eq:asyug}. Since 
equation~\eqref{eq:sdep} is $s$-dependent, the additive group $\R$ 
does not act on the spaces of solutions, and we shall refer to 
these as the {\bf Morse-Bott moduli spaces of $s$-dependent Floer 
trajectories}. One shows as in~\cite[Section~4.1]{BOauto} that, for a 
generic choice of the homotopy $J_s$, these are smooth manifolds whose 
respective dimensions are given by~\eqref{eq:dim}.  
 
Let us now choose for each $S_\gamma$, $\gamma\in\cP(H_\pm)$ a 
perfect Morse function $f_\gamma^\pm:S_\gamma\to \R$. Given 
$\og\in\cP(H_-)$, $\ug\in\cP(H_+)$, $p\in 
\mathrm{Crit}(f^-_\og)$,  
$q\in\mathrm{Crit}(f^+_\ug)$, $m\ge 1$, and $i\in\{1,\dots,m\}$, we define  
$$ 
\cM^A_m(p,q;H_s,\{f_\gamma^\pm\},J_s) 
$$ 
as the (disjoint) union over $\tgamma_1,\dots,\tgamma_{i-1}\in \cP(H_-)$, 
$\tgamma_i,\dots,\tgamma_{m-1}\in \cP(H_+)$, and $A_1+\dots+A_m=A$ 
of the following fibered products (with the convention 
$\tgamma_0:=\og$, $\tgamma_m:=\ug$)  
\begin{eqnarray*} 
&& 
W^u(p)  
\times_{\oev} 
(\cM^{A_1}(S_{\tgamma_0},S_{\tgamma_1};H_-,J_-)\!\times\!\R^+) \\ 
&& {_{\varphi_{f^-_{\tgamma_1}}\!\circ\uev}}\!\times _{\oev}  
\ldots {_{\varphi_{f^-_{\tgamma_{i-2}}}\!\circ\uev}}\!\times _{\oev}  
(\cM^{A_{i-1}}(S_{\tgamma_{i-2}},S_{\tgamma_{i-1}};H_-,J_-)\!\times\!\R^+) \\ 
&&  
{_{\varphi_{f^-_{\tgamma_{i-1}}}\!\circ\uev}\times_{\oev}}  
(\cM^{A_i}(S_{\tgamma_{i-1}},S_{\tgamma_i};H_s,J_s)\!\times\!\R^+) \\ 
&&{_{\varphi_{f^+_{\tgamma_i}}\!\circ\uev}\times_{\oev}}  
(\cM^{A_{i+1}}(S_{\tgamma_i},S_{\tgamma_{i+1}};H_+,J_+)\!\times\!\R^+) \\ 
&&  
{_{\varphi_{f^+_{\tgamma_{i+1}}}\!\circ\uev}\times_{\oev}} \ldots\, 
{_{\varphi_{f^+_{\tgamma_{m-1}}}\!\!\circ\uev}}\!\!\times 
_{\oev}   
\cM^{A_m}(S_{\tgamma_{m-1}},\!S_{\tgamma_m};H_+,J_+)  
{_{\uev}\times} W^s(q). 
\end{eqnarray*} 
We note the similarity to the fibered product 
defining~\eqref{eq:bigM}, the difference being that the term  
$\cM^{A_i}(S_{\tgamma_{i-1}},S_{\tgamma_i};H,J)$ being replaced by  
$\cM^{A_i}(S_{\tgamma_{i-1}},S_{\tgamma_i};H_s,J_s)$.  
One shows as in~\cite[Lemma~3.6]{BOauto} that, for a generic choice of 
the collection of Morse functions $\{f_\gamma^\pm\}$, each space  
$\cM^A_m(p,q;H_s,\{f_\gamma^\pm\},J)$ is a smooth manifold of 
dimension  
\begin{eqnarray*} 
\lefteqn{\dim \, \cM^A_m(p,q;H_s,\{f_\gamma^\pm\},J_s)} \\ 
& = & \mu(\og) +\ind(p) - \mu(\ug) - 
\ind(q) +2\langle c_1(TW),A\rangle. 
\end{eqnarray*} 
The union over $m\ge 1$ of the moduli spaces 
$\cM^A_m(p,q;H_s,\{f_\gamma^\pm\},J_s)$ is denoted by  
$\cM^A(p,q;H_s,\{f_\gamma^\pm\},J_s)$ and is called  
the {\bf moduli space of 
$s$-dependent Morse-Bott broken trajectories}.  
The topology of the compactification\break  
$\overline \cM^A(p,q;H_s,\{f_\gamma^\pm\},J_s)$ 
is described similarly to the $s$-independent case.  
 
An increasing homotopy as above defines a {\bf continuation morphism} 
$$ 
\sigma_{H_+,H_-}:BC_*^a(H_-)\to BC_*^a(H_+), 
$$ 
which preserves the degree and is obtained by a count of rigid 
configurations in $\cM^A(p,q;H_s,\{f_\gamma^\pm\},J_s)$. Via the 
identification proved in~\cite[Theorem~3.7]{BOauto} between the 
Morse-Bott complexes of $H_\pm$ and the Floer  
complexes of suitable time-dependent perturbation of these 
Hamiltonians, 
the continuation morphism $\sigma_{H_+,H_-}$ coincides with the usual 
continuation morphism in Floer homology~\cite{FH1}. This can be seen 
for example by repeating the gluing arguments of~\cite{BOauto} in the 
context of $s$-dependent families.

 
\section{Contact homology} \label{sec:bigcontact} 
 
\subsection{Linearized contact homology} \label{sec:contact} 
 
In this section we define linearized contact homology of the contact 
manifold $(M,\xi)$ with the symplectic filling $(W,\om)$ 
following~\cite{BEE}.  
 
For each free homotopy class of loops $b$ 
in $M$ we denote by $\cP_\lambda^b$ the set of all ${\gamma'}\in 
\cP_\lambda$ in the homotopy class $b$. The inclusion 
$i:M\hookrightarrow W$ induces a map (still denoted by $i$) between 
the sets of free homotopy classes of loops in $M$ and $W$ 
respectively. For each 
free homotopy class $a$ in $W$ we denote 
$$ 
\cP_\lambda^{i^{-1}(a)} := \bigcup _{b\in i^{-1}(a)} \cP_\lambda^b. 
$$ 
We assume in this section that all the closed Reeb orbits are 
transversally nondegenerate in $M$. This means that, for every orbit ${\gamma'}$ 
of period $T>0$, we have 
$$ 
\det \left( \one - d\phi_\lambda^T({\gamma'}(0))|_\xi \right) \neq 0. 
$$ 
This can always be achieved by an arbitrarily small perturbation of 
$\lambda$ or, equivalently, of $X$, which does not change the 
symplectic homology groups. In this situation one can assign to each 
${\gamma'}\in \cP_\lambda$ a Conley-Zehnder index $\mu_{CZ}({\gamma'})$ 
according to the following recipe. 
 
We fix a reference loop $l_a:S^1\to \widehat W$ for each free homotopy 
class $a$ in $\widehat W$ such that $[l_a]=a$. If $a$ is the trivial 
homotopy class we choose $l_a$ to be a constant loop, and we require that 
$l_{a^{-1}}$ coincides with the loop $l_a$ with the opposite 
orientation. We also choose symplectic trivializations  
$$ 
\Phi_a:S^1\times \R^{2n} \to l_a^*T\widehat W 
$$ 
for each class $a$. If $a$ is the trivial homotopy class we choose the 
trivialization to be constant, and we require that 
$\Phi_{a^{-1}}(\theta,\cdot)=\Phi_a(-\theta,\cdot)$, $\theta\in 
S^1=\R/\Z$.  
 
We fix a reference loop $l_b:S^1\to M$ for each free 
homotopy class $b$ in $M$ such that $[l_b]=b$. If $b$ is the 
trivial homotopy class we choose $l_b$ to be a constant loop and we 
require that $l_{b^{-1}}$ coincides with $l_b$ with the opposite 
orientation. We define symplectic trivializations 
$$ 
\Phi_b:S^1\times \R^{2n-2} \to l_b^*\xi 
$$ 
as follows. For each class $b$ we choose a homotopy $h_{ab}:S^1\times 
[0,1] \to W$ from $l_a$, $a=i(b)$ to $l_b$ such that 
\begin{equation} \label{eq:CRtrivhom} 
h_{a^{-1}b^{-1}}(\tau,\cdot) = h_{ab}(-\tau,\cdot). 
\end{equation} 
We extend the trivialization $\Phi_a:S^1\times \R^{2n} \to 
l_a^*T\widehat W$ over the homotopy $h_{ab}$ to get a trivialization 
$\Phi'_b:S^1\times \R^{2n} \to l_b^*T\widehat W$. This trivialization 
is homotopic to another one, still denoted $\Phi'_b$, such that 
\begin{eqnarray} 
   \Phi'_b(S^1\times \R^{2n-2} \times \{0\} \times \{0\}) & = & 
   l_b^*\xi, \nonumber \\ 
   \Phi'_b(S^1\times \{0\} \times \R \times \{0\}) & = & l_b^*\langle 
   \frac \p {\p t} \rangle, \label{eq:homotopic_triv} \\ 
   \Phi'_b(S^1\times \{0\} \times \{0\} \times \R) & = & l_b^*\langle 
   R_\lambda \rangle. \nonumber 
\end{eqnarray} 
We define $\Phi_b:=\Phi'_b|_{S^1\times \R^{2n-2} \times \{0 \} \times 
   \{0\}}$. If $b$ is the trivial homotopy 
class we choose $h_{ab}$ to be a path of 
constant loops, so that $\Phi_b$ is constant. 
 
We fix for each ${\gamma'}\in 
\cP_\lambda$ a map $\sigma_{\gamma'}:\Sigma \to M$, with $\Sigma$ 
a Riemann surface with two boundary components 
$\p_0\Sigma$ (with the opposite boundary orientation) and 
$\p_1\Sigma$ (with the boundary orientation), satisfying 
\begin{equation} \label{eq:sigmaM} 
\sigma|_{\p_0\Sigma} = l_{[{\gamma'}]}, \qquad \sigma|_{\p_1\Sigma}={\gamma'}. 
\end{equation} 
For each ${\gamma'}\in \cP_\lambda$ there exists a unique (up to 
homotopy) trivialization 
$$ 
\Phi: \Sigma \times \R^{2n-2} \to \sigma_{\gamma'}^*\xi 
$$ 
such that $\Phi=\Phi_{[{\gamma'}]}$ on $\p_0\Sigma \times \R^{2n-2}$. Let 
\begin{equation} \label{eq:CRtriv} 
\Psi:[0,T]\to \textrm{Sp}(2n-2),  
\Psi(\tau):= \Phi^{-1} \circ d\phi_\lambda^\tau(p) \circ \Phi,  
p\in {\gamma'}([0,T]). 
\end{equation} 
Because ${\gamma'}$ is nondegenerate we can define the 
{\bf Conley-Zehnder index} $\mu({\gamma'})$ by 
\begin{equation} 
   \label{eq:CRmu} 
   \mu({\gamma'}):=\mu({\gamma'},\sigma_{\gamma'}) := \mu_{CZ}(\Psi), 
\end{equation} 
where $\mu_{CZ}(\Psi)$ is the Conley-Zehnder index of a path of 
symplectic matrices~\cite{RS}. 
 
\begin{remark} \rm  
Given  $B\in H_2(M;\Z)$, we define a map $\sigma_{\gamma'} \# B$ up to homology 
as the connected sum of $\sigma_{\gamma'}$ with a surface  
representing $B$.     
If, in the previous construction, we replace 
   $\sigma_{\gamma'}$ with $\sigma_{\gamma'} \# B$, then the 
   resulting index will be 
  \begin{equation} \label{eq:CRc1} 
   \mu({\gamma'},\sigma_{\gamma'} \# B) = \mu({\gamma'},\sigma_{\gamma'}) + 
   2\langle c_1(\xi),B\rangle. 
  \end{equation} 
Note that $c_1(\xi)=i^*c_1(TW)$ because $i^*TW=\xi \oplus \langle 
\frac \p {\p t},R_\lambda \rangle$. Moreover, the parity of 
$\mu({\gamma'})$ is well-defined independently of the trivialization of 
$\xi$ along ${\gamma'}$. 
\end{remark} 
 
For each simple orbit ${\gamma'}\in \cP_\lambda$ we denote by ${\gamma'}^k$, 
$k\in \Z^+$ its positive iterates. The parity of the Conley-Zehnder 
index of all the odd, respectively even iterates is the same. If 
these two parities differ we say that all even iterates ${\gamma'}^{2k}$, 
$k\in \Z^+$ are {\bf bad orbits}. It can be seen that 
the even iterates of a simple orbit ${\gamma'}$ of period $T$ 
are bad if and only if $d\phi_\lambda^T(p)|_\xi$, $p\in {\gamma'}([0,T])$ 
has an odd number of real negative eigenvalues strictly smaller than 
$-1$. The orbits in $\cP_\lambda$ which are not bad are called {\bf 
   good orbits}. 
 
We define a grading on $\Lambda_\om$ by $|e^A|:= -2\langle 
c_1(TW),A\rangle$. Note that, if $A=i_*(B)$, $B\in H_2(M;\Z)$ then  
$|e^A|=-2\langle c_1(\xi),B\rangle$. For each free homotopy class of 
loops $a$ in $W$, the {\bf contact chain group with 
  coefficients in $\Lambda_\om$} is denoted by 
$C_*^{i^{-1}(a)}(\lambda)$ and is defined as the free 
$\Lambda_\om$-module generated by all good orbits ${\gamma'}\in 
\cP_\lambda^{i^{-1}(a)}$. The grading is given by 
$$ 
|e^A{\gamma'}| := \mu({\gamma'}) -2 \langle c_1(TW),A \rangle + n - 3. 
$$ 
We define the {\bf reduced Conley-Zehnder index} $\overline 
\mu({\gamma'}):=\mu({\gamma'}) + n -3$, so that the grading is 
$|e^A{\gamma'}|=\overline \mu({\gamma'}) - 2 \langle c_1(TW),A \rangle$. 
 
We call the symplectic manifold $(M\times \R,d(e^t\lambda))$ the {\bf 
   symplectization} of $(M,\xi)$. Its symplectomorphism type does not 
depend on $\lambda$, but only on the isotopy class of $\xi$. Let 
$\cJ(\lambda)$ denote the set of {\bf  admissible almost complex 
structures} on $M\times \R$, consisting of 
elements $J_\infty$ satisfying  
\begin{equation}  \label{eq:J} 
\left\{ \begin{array}{rcl} 
J_{\infty \, (p,t)} |_\xi & = & J_0, \\[.2cm] 
J_{\infty \, (p,t)} \frac \p  {\p t} & = & R_\lambda 
\end{array}\right. 
\end{equation} 
on $M\times \R$. 
Here $J_0$ is any compatible complex structure on the symplectic 
bundle $(\xi,d\lambda)$ which is independent of $\theta$ and $t$. 
 
\medskip  
 
>From now on, we choose for each $\gamma'\in\cP_\lambda$ a 
point $P_{\gamma'}$ on the geometric image of $\gamma'$.  
 
\medskip  
 
Let $J_\infty\in \cJ(\lambda)$, ${\og'},{\ug'}, {\gamma'}_1,\ldots,{\gamma'}_k\in 
\cP_\lambda$ and $B\in H_2(M;\Z)$. We define the space  
$$ 
\widehat \cM^B({\og'},{\ug'},\gamma'_1,\ldots,\gamma'_k;J_\infty) 
$$ 
of {\bf punctured $J_\infty$-holomorphic cylinders} as the set of 
tuples  
$$ 
(F,\oL,\uL,L_1,\dots,L_k),  
$$ 
consisting of a solution 
$$ 
F=(f,a):\R \times S^1 \setminus \{z_1,\ldots,z_k\} = \C P^1\setminus 
\{0,\infty,z_1,\ldots,z_k\} \to M\times \R 
$$ 
of the Cauchy-Riemann equation 
\begin{equation} 
   \label{eq:CR} 
   \p_s F + J_\infty \p_\theta F =0, 
\end{equation} 
and of half-lines $L_i\subset T_{z_i}(\R\times S^1)$, $\oL\subset T_0 
\C P^1$, $\uL\subset T_\infty \C P^1$,  
subject to the asymptotic conditions  
\begin{equation} 
   \label{eq:a} 
   \lim_{s\to \pm \infty} a(s,\theta) = \mp\infty, 
\end{equation} 
\begin{equation} 
   \label{eq:f} 
   \lim_{s\to -\infty} f(s,\theta) ={\og'}(-\oT \theta), \qquad  \lim_{s\to 
     +\infty} f(s,\theta) =\ug'(-\uT \theta) 
\end{equation} 
uniformly in $\theta$, and  
\begin{equation} \label{eq:L} 
  \lim_{z\to 0, \ z\in \oL} f(z)=P_{\og'},   \lim_{z\to \infty, \ z\in 
  \uL} f(z)=P_{\ug'}.  
\end{equation} 
Moreover, we require that, given polar coordinates 
$(\rho_i,\theta_i)\in]0,1]\times \R/\Z$ around $z_i$, we have  
\begin{equation} \label{eq:Fpunc} 
\lim_{z\to z_i} a(z)=-\infty, \quad  
\lim_{\rho_i\to 0} f(\rho_i,\theta_i)=\gamma'_i(T_i\theta_i), \quad 
i=1,\ldots,k,  
\end{equation} 
\begin{equation} \label{eq:Lpunc}  
\lim_{z\to z_i, \ z\in L_i} f(z)=P_{\gamma_i'}, \quad i=1,\ldots,k.  
\end{equation}   
In addition, we require that  
\begin{equation} 
   \label{eq:CRA} 
   [(\sigma_{{\ug'}} \cup \sigma_{\gamma'_1} \cup \ldots \cup 
   \sigma_{\gamma'_k}) \# f] =[\sigma_{{\og'}} \# B].  
\end{equation} 
These convergence conditions are to be understood  
for some $R_\lambda$-parametrized representatives of 
$\gamma'_i$, respectively $\og'$, $\ug'$. By the conditions $z\in L_i$  
we mean that $z$ belongs to some curve with endpoint $z_i$ and 
asymptotically tangent to $L_i$ (similarly for $\oL$, $\uL$).  
The half-lines $\oL,\uL,L_1,\dots,L_k$ are called {\bf asymptotic 
  markers}. We refer to 
Figure~\ref{fig:Phi}.(b) on  
page~\pageref{fig:Phi} for a  
representation of these objects (see also the notion of a capped 
punctured $J$-holomorphic cylinder at the end of this section). 
 
\medskip  
 
\noindent {\it Notational convention.} To simplify, we shall use the 
shorthand notation $F\in \widehat 
\cM^B({\og'},{\ug'},\gamma'_1,\ldots,\gamma'_k;J_\infty)$ for a tuple  
$(F,\oL,\uL,L_1,\dots,L_k)$. The same convention 
applies for all subsequent (moduli) spaces.  
 
\begin{remark} {\rm  
Our convention for the asymptotic behaviour is 
different from the usual one in contact homology and is motivated by 
the usual conventions for symplectic homology.  
} 
\end{remark}  
 
\begin{remark} \rm 
   Under the nondegeneracy assumption on 
   ${\og'},{\ug'},\gamma'_1,\ldots,\gamma'_k$,  
   the conditions~\eqref{eq:a}, \eqref{eq:f} and~\eqref{eq:Fpunc} are equivalent to the 
   finiteness of the Hofer energy~\cite[Theorem~1.2]{HWZ} 
$$ 
\cE(F):= \sup _{\phi\in \cC} \int_{\R\times S^1\setminus 
  \{z_1,\ldots,z_k\}} F^*d(\phi\lambda),  
$$ 
where $\cC:= \left\{ \phi \in C^\infty(\R,[0,1]) \, | \, \phi'\ge 0 
\right\}$. We define the {\bf contact action functional} 
$$ 
\cA_\lambda: C^\infty(S^1,M) \to \R \ : \ {\gamma'}\mapsto \int_{S^1} 
{\gamma'}^*\lambda. 
$$ 
For every $F\in \widehat 
\cM^B({\og'},{\ug'},{\gamma'}_1,\ldots,{\gamma'}_k;J_\infty)$ we have  
$\cE(F)=\cA_\lambda({\og'})$. 
\end{remark} 
 
The group of biholomorphisms on the domain of $F$ is $\R\times S^1$ 
and it acts freely on the space $\widehat 
\cM^B({\og'},{\ug'},\gamma'_1,\ldots,\gamma'_k;J_\infty)$ by  
$$ 
h\cdot F := F\circ h^{-1}, \quad h\cdot L:= h_* L 
$$  
with $L\in\{\oL,\uL,L_1,\dots,L_k\}$. The {\bf moduli space of 
  punctured $J_\infty$-holomorphic cylinders} is defined by  
$$ 
\cM^B({\og'},{\ug'},\gamma'_1,\ldots,\gamma'_k;J_\infty):= \widehat 
\cM^B({\og'},{\ug'},\gamma'_1,\ldots,\gamma'_k;J_\infty)/(\R\times S^1).  
$$ 
It follows from equation~\eqref{eq:indDFpunc} below that the virtual 
dimension of this moduli space is  
$$ 
\mu(\og') - \mu(\ug') + 2\langle c_1(\xi),B\rangle - \sum_{i=1}^k 
\mu(\gamma'_i). 
$$ 
 
An almost complex structure $J_\infty\in\cJ(\lambda)$ is called {\bf 
   regular for cylinders} if the linearized operator below is 
surjective for all ${\og'},{\ug'}\in \cP_\lambda$, $B\in H_2(M;\Z)$ and 
$F\in \widehat \cM^B({\og'},{\ug'};J_\infty)$. The linearized operator is 
$$ 
D_F:W^{1,p,d}(\R\times S^1,F^*T(M\times \R))\oplus \R^4 \to L^{p,d}(\R\times 
S^1,F^*T(M\times \R)), 
$$ 
$$ 
D_F \zeta := \nabla_s\zeta +J_\infty \nabla_\theta\zeta +  
(\nabla_\zeta J_\infty)\p_\theta F, 
$$ 
where $\zeta:= \zeta_0 + v_-^1\zeta_-^1 + v_-^2\zeta_-^2 + 
v_+^1\zeta_+^1 + v_+^2\zeta_+^2$, $\zeta_0\in W^{1,p,d}(\R\times 
S^1,F^*T(M\times \R))$, $v_\pm^i\in \R$, $i=1,2$, $p>2$, $d>0$ small 
enough. The sections $\zeta_\pm^i$, $i=1,2$ are asymptotically 
constant with the following asymptotic values: 
$$ 
\zeta_-^1(s,\theta) = R_\lambda, \qquad \zeta_-^2(s,\theta)=\p / {\p 
   t} \qquad \textrm{for } s\le -1, 
$$ 
$$ 
\zeta_-^i(s,\theta) = 0 \qquad \textrm{for } s\ge 1, 
$$ 
$$ 
\zeta_+^i(s,\theta)=\zeta_-^i(-s,\theta), \qquad i=1,2. 
$$ 
The spaces $W^{1,p,d}(\R\times S^1, F^*T(M\times \R))$ and 
$L^{p,d}(\R\times S^1, F^*T(M\times \R))$ are the completions of 
$C^\infty(\R\times S^1,F^*T(M\times \R))$ with respect to the norms 
$$ 
\|\zeta\|_{1,p,d}:= \left( \int_{\R\times S^1} (\|\zeta\|^p + 
   \|\nabla_s\zeta \|^p + \|\nabla_\theta \zeta \|^p) e^{d|s|} \, 
   dsd\theta \right)^{1/p}, 
$$ 
$$ 
\|\zeta\|_{p,d}:= \left( \int_{\R\times S^1} \|\zeta\|^p e^{d|s|} \, 
   dsd\theta \right)^{1/p}. 
$$ 
Because the orbits ${\og'}$, ${\ug'}$ are transversally nondegenerate the 
operator $D_F$ is Fredholm with index~\cite[Proposition~4]{BM} 
\begin{equation} 
   \label{eq:indDF} 
   \ind(D_F) = \mu({\og'}) - \mu({\ug'}) + 2\langle c_1(\xi),B \rangle +2, 
   \ F\in \widehat \cM^B({\og'},{\ug'};J_\infty), 
\end{equation} 
so that the virtual dimension of $\cM^B({\og'},{\ug'};J_\infty)$ is $\mu({\og'}) - 
\mu({\ug'}) + 2\langle c_1(\xi),B \rangle$. 
 
The above discussion can be generalized in a fairly straightforward 
way in order to define {\bf regular almost complex structures for 
  punctured holomorphic cylinders}. The relevant operator now has   
index  
\begin{equation} \label{eq:indDFpunc} 
\ind(D_F)=\mu(\og') - \mu(\ug') + 2\langle c_1(\xi),B\rangle - \sum_{i=1}^k 
\mu(\gamma'_i)+2  
\end{equation}  
for $F\in 
\widehat\cM^B(\og',\ug',\gamma'_1,\ldots,\gamma'_k;J_\infty)$.  
The point is that, unlike in Floer 
homology, the existence of almost complex structures which are regular 
for punctured cylinders is not guaranteed.  
Indeed, since the almost complex structure is not  
domain dependent, a ramified covering of a punctured holomorphic 
cylinder is again holomorphic, and its index may be smaller than the 
dimension of the space of ramified coverings. In that case the 
corresponding linearized operator cannot be surjective.   
 
Let now $J$ be a time-independent almost complex structure on 
$\widehat W$ which is compatible with $\widehat \om$ and whose 
restriction $J_\infty$ to $M\times \R^+$ is translation invariant and 
corresponds to an element of $\cJ(\lambda)$.  
A {\bf $J$-holomorphic plane in $\widehat W$} is a $J$-holomorphic map 
$F:\C=\C P^1\setminus \{\infty\} \to \widehat W$ such that, for $|z|$ 
large enough $F(z)\in M\times \R_+$ and, writing $F=(f,a)$, we have 
$a(z)\to \infty$, $|z|\to \infty$ and there exist 
$\gamma'\in\cP_\lambda$ such that $f(re^{2i\pi\theta})\to 
\gamma'(T\theta)$, $r\to \infty$ uniformly in $\theta$. As for 
cylinders, this convergence condition has to be understood with 
respect to some $R_\lambda$-parametrized representative of 
$\gamma'$. Given $\gamma'\in\cP_\lambda$, $A\in  
H_2(W;\Z)$ we define the {\bf space $\widehat 
  \cM^A({\gamma'},\emptyset;J)$ of $J$-holomorphic planes} as the set 
of pairs $(F,L)$ with $F$ as above and $L\subset T_\infty \C P^1$ a 
half-line, such that $\lim_{z\to\infty,\ z\in L}f(z)=P_{\gamma'}$ and 
$[F]=[\sigma_\gamma \# A]$. The group of biholomorphisms of the domain 
$\C$ consists of affine transformations and has real 
dimension $4$. It acts by  
$$ 
h\cdot F := F\circ h^{-1},\qquad h\cdot L:=h_*L. 
$$ 
The quotient is  
denoted by $\cM^A(\gamma',\emptyset;J)$ and is called the {\bf moduli 
  space of $J$-holomorphic planes}.

The relevant linearized operator now has index 
$$ 
\ind(D_F) = \overline \mu({\gamma'}) + 2\langle c_1(TW),A \rangle +4, \qquad F 
\in \widehat \cM^A ({\gamma'},\emptyset;J), 
$$ 
so that the virtual dimension of $\cM^A ({\gamma'},\emptyset;J)$ is 
$\overline \mu({\gamma'}) + 2\langle c_1(TW),A \rangle$. We say that 
an almost complex structure $J$ on $\widehat W$ is {\bf regular for 
  holomorphic planes} if the linearized operator $D_F$ is surjective 
for every $\gamma'\in\cP_\lambda$, $A\in H_2(W;\Z)$ and $F\in \widehat 
\cM^A ({\gamma'},\emptyset;J)$. Like for punctured holomorphic 
cylinders, the existence of such regular almost complex structures is 
not guaranteed.  
 
\medskip  
 
For the definitions that follow we assume that $J$ is regular for the 
relevant holomorphic curves. A list of examples in which this 
assumption is satisfied is given in Remark~\ref{rmk:transv_exples} 
below.  
 
\medskip  
 
One can associate a sign $\epsilon(F)$ to each element $[F]\in 
\cM^A({\gamma'},\emptyset;J)$ such that  
$\overline \mu({\gamma'}) + 2\langle c_1(TW),A \rangle=0$ (see~\cite{BM}).  
We define a homomorphism  
$$ 
e : C_*^{i^{-1}(a)}(\lambda) \to \Lambda_\om 
$$ 
by  
$$ 
e({\gamma'}) := \sum_{\stackrel {A\in H_2(W;\Z)} {|e^A|=|{\gamma'}|}} \big(\sum 
_{[F]\in\cM^A({\gamma'},\emptyset;J)}  
\epsilon(F) \big) e^A. 
$$ 
 
\begin{remark} \label{rmk:augmentation} {\rm  
 The homomorphism $e$ is obtained from the natural {\bf augmentation} 
 on the differential graded algebra for the full contact homology of 
 $(M,\lambda)$ defined by the symplectic filling $(W,\om)$ 
 (see~\cite{BEE} for details on augmentations).  
} 
\end{remark}  
 
Given $\gamma'\in\cP_\lambda$, we denote by $\kappa_{\gamma'}\in\Z^+$ 
  its {\bf multiplicity}. It is the largest integer such that 
  $\gamma'(\theta+\frac 1 {\kappa_{\gamma'}}) = \gamma'(\theta)$, 
  $\theta\in\R/\Z$.   
 
\medskip   
 
If $k\neq 0$, ${\og'}\neq{\ug'}$ or $0\neq B \in H_2(M;\Z)$, the additive 
group $\R$ acts freely on the moduli space 
$\cM^B({\og'},{\ug'},\gamma'_1,\ldots,\gamma'_k;J_\infty)$ by 
translations in the $t$ direction  
$$ 
t_0 \cdot [(f,a)] := [(f,a+t_0)]. 
$$ 
We denote the quotient by 
$\cM^B({\og'},{\ug'},\gamma'_1,\ldots,\gamma'_k;J_\infty)/\R$. This 
space can be compactified to~$\overline 
\cM^B({\og'},{\ug'},\gamma'_1,\ldots,\gamma'_k;J_\infty)$, a space 
consisting of   
$J_\infty$-holomorphic buildings with a tree 
structure~\cite[Theorem~10.1]{BEHWZ}. If $k=0$, ${\og'}={\ug'}$ and 
$B=0$, the moduli space consists of $\kappa:=\kappa_{\og'}=\kappa_{\ug'}$ 
points. The underlying holomorphic curve is the constant cylinder over 
the orbit ${\og'} = {\ug'}$, and the equivalence class of the pair of 
asymptotic markers $(\oL,\uL)$ is determined by the difference 
$(\Arg(\oL)-\Arg(\uL))/2\pi\in\R/\Z$, which is a multiple of $1/\kappa$. The 
action of $\R$ is in this case trivial.

In the case $\mu ({\og'}) -\mu({\ug'}) +2\langle c_1(\xi), B \rangle -
\sum_{i=1}^k  
\overline \mu({\gamma'}_i) =1$, the moduli space 
$\cM^B({\og'},{\ug'},\gamma'_1,\ldots,\gamma'_k;J_\infty)/\R$ is 
compact and therefore consists  
of a finite number of points. One can associate a sign 
$\epsilon(F)$ to each element $[F]$ of this moduli space~\cite{BM} and 
we define the {\bf linearized contact differential} 
$$  
\p : C_*^{i^{-1}(a)}(\lambda) \to C_{*-1}^{i^{-1}(a)}(\lambda) 
$$ 
by 
\begin{equation} 
   \label{eq:CRdiff} 
\p {\og'} := 
\sum _{\substack{{\ug'},\gamma'_1,\ldots,\gamma'_k,B \\ 
  |{\ug'} e^B| + \sum |\gamma'_i| = |{\og'}| -1 } } 
\hspace{-.3cm} \frac {{\scriptstyle 
    n^B({\og'},{\ug'},\gamma'_1,\ldots,\gamma'_k;J_\infty)}} 
{{\scriptstyle \kappa_{\ug'}\prod_{i=1}^k \kappa_{\gamma_i'}}}  
\ e(\gamma'_1)\ldots e(\gamma'_k) e^B {\ug'}, 
\end{equation} 
where  
$$ 
n^B({\og'},{\ug'},\gamma'_1,\ldots,\gamma'_k;J_\infty) :=  
\sum_{[F]\in 
  \cM^B({\og'},{\ug'},\gamma'_1,\ldots,\gamma'_k;J_\infty)/\R}  
\epsilon(F).  
$$ 
The reader is warned that the classes $B$ in the above sum live in 
$H_2(M;\Z)$, but nevertheless the coefficient in front of ${\ug'}$ is an 
element of $\Lambda_\om$ due to the factors $e(\gamma'_i)$, 
$i=1,\ldots,k$. As a matter of fact, the fraction in~\eqref{eq:CRdiff} 
is an integer because our moduli spaces involve asymptotic markers.  
 
Since $e$ comes from the natural augmentation defined by $(W,\om)$, it 
follows that (see~\cite{BEE}) 
\begin{equation} \label{eq:ed} 
 \p\circ\p = 0 \quad \mbox{and} \quad e\circ \p=0. 
\end{equation}  
We define the {\bf linearized contact homology groups} of the 
pair $(\lambda,J)$ by 
$$ 
HC_*^{i^{-1}(a)}(\lambda,J) := H_*(C_*^{i^{-1}(a)}(\lambda),\p). 
$$ 
 
The linearized contact homology groups 
$HC_*^{i^{-1}(a)}(\lambda,J)$ depend only  
on the symplectic filling $(W,\om)$ of $(M,\xi)$. Since the former is 
part of the data in the context of the present paper, we simply denote 
the resulting homology groups by  
$$ 
HC_*^{i^{-1}(a)}(M,\xi) 
$$  
without reference to $W$.  
 
We shall need in the sequel the following alternative description of 
the differential $\p$ for linearized contact homology. Given $A\in 
H_2(W;\Z)$, $\og',\ug'\in \cP_\lambda$ we define the {\bf 
  moduli space of capped punctured $J$-holomorphic cylinders}  
$$ 
\cM^A_c(\og',\ug';J) 
$$ 
as the set of equivalence classes of pairs $F=(F',F'')$, where 
$F'$ is an element of $\cM^B(\og',\ug',\gamma'_1,\ldots,\gamma'_k; J_\infty)$,  
$\gamma'_1,\ldots,\gamma'_k\in\cP_\lambda$, $B\in H_2(M;\Z)$, $F''$ is a 
collection of $k$ $J$-holomorphic planes in $\widehat W$, of total homology 
class $A-B\in H_2(W;\Z)$, and with asymptotics at their positive punctures 
corresponding to     
$\gamma'_1,\ldots,\gamma'_k$. Recall that $F'$ and $F''$ are endowed 
with asymptotic markers $(\oL',\uL',L_i')$ and $(L_i'')$, 
$i=1,\dots,k$ which determine conformal identifications of the tangent 
spaces at the punctures with asymptotes $\gamma_i'$. Two pairs 
$\big((\oL'_0,\uL'_0,L'_{i,0}),(L''_{i,0})\big)$, 
$\big((\oL'_1,\uL'_1,L'_{i,1}),(L''_{i,1})\big)$ 
corresponding to the same maps $F',F''$ are 
equivalent if $\oL'_0=\oL'_1$, $\uL'_0=\uL'_1$, and   
$(L'_{i,0},L''_{i,0})$  
determine the same conformal identifications as  
$(L'_{i,1},L''_{i,1})$. This last condition is equivalent to the 
equality  
\begin{equation} \label{eq:Arg} 
\Arg(L'_{i,0})-\Arg(L'_{i,1})=\Arg(L''_{i,0})-\Arg(L''_{i,1}), \quad 
i=1,\dots,k 
\end{equation}   
with respect to fixed conformal charts at the corresponding 
punctures. The dimension of this moduli space is  
$$ 
\mu(\og')-\mu(\ug') + 2\langle c_1(TW),A\rangle. 
$$ 
We refer to Figure~\ref{fig:Phi}.(b) for a  
representation of such objects.  
The differential for linearized contact homology can then be rewritten 
as 
\begin{equation} \label{eq:lindiffcapped} 
 \p \og' = \sum _{\substack{\ug',A \\ 
  |\ug' e^A| = |\og'| -1 } } 
\frac 1 {{\scriptstyle \kappa_{\ug'}}} 
\sum_{[F]\in \cM^A_c(\og',\ug';J)/\R} 
  \epsilon(F) e^A \ug'. 
\end{equation}  
The sign $\epsilon(F)$ is defined as the product of the signs of 
the components of $F$.

\begin{remark} \label{rmk:polyfolds} {\rm  
 The proof of the invariance of $HC_*^{i^{-1}(a)}(\lambda,J)$ with 
 respect to 
 $\lambda$ and $J$ makes use of the polyfold formalism currently being 
 developed by Hofer, Wysocki and Zehnder~\cite{H,HWZ2}.  
} 
\end{remark}

\begin{remark} {\rm 
Note that the linearized contact homology is related to the cylindrical  
contact homology defined in~\cite{EGH}, in the following situation.  
Let us assume that   
 $c_1(\xi)=0$, and that there are no Reeb orbits $\gamma'$ with 
 reduced Conley-Zehnder index $\bar \mu(\gamma')=-1,0,1$ and which are 
 contractible in $M$, so that cylindrical contact homology 
 is well-defined. Then, if $c_1(TW)=0$ and if there are no Reeb orbits 
 $\gamma'$ of reduced Conley-Zehnder index $\bar \mu(\gamma')=0$ which 
 are contractible in $W$, the linearized contact homology groups are 
 isomorphic to the cylindrical contact homology groups.   
} 
\end{remark}

\begin{remark}[Note on transversality] \label{rmk:transv_exples} 
 Let $a$ be a free homotopy class of loops in $W$. In order for the 
 linearized contact differential to be well-defined  
 on $i^{-1}(a)$ we need that the almost complex structure $J$ on 
 $\widehat W$ satisfies the following conditions. 
 
\begin{enum} 
\item[$(A)$] $J$ is regular for holomorphic planes belonging to moduli spaces 
$\cM^A(\gamma',\emptyset;J)$ of virtual dimension $\le 0$; 
 
\item[$(B_a)$] $J_\infty$ is regular for punctured holomorphic 
cylinders asymptotic at $\pm\infty$ to closed Reeb orbits in  
$i^{-1}(a)$, belonging to moduli spaces of virtual dimension $\le 2$, 
and which are asymptotic at the punctures to closed  
Reeb orbits $\gamma'$ such that $\cM^A(\gamma',\emptyset;J)\neq 
\emptyset$ and has virtual dimension $0$.  
\end{enum} 
 
We expect these technical 
assumptions to be completely  
removed using the new ongoing approach to transversality by Cieliebak 
and Mohnke (see~\cite{CM} for the symplectic case), or using the  
polyfold theory developed by Hofer, Wysocki and Zehnder~\cite{H,HWZ2}. 
We give below the list of examples known to us in which both these 
conditions are satisfied. In many cases condition~$(A)$ is empty, 
so that linearized contact homology reduces to 
cylindrical contact homology.  
 
\begin{enum} 
\item stabilizations $W:=V\times D^2$ of subcritical Stein manifolds 
 $V$ with $c_1(V)=0$, and in particular the standard balls 
 $B^{2n}$, $n\ge 2$, for the trivial free homotopy 
 class. 
 In this case, if one chooses on $W$ a symplectic  
 form corresponding to sufficiently thin handles, then there are no 
 $J$-holomorphic planes of index $\le 0$ since $\bar \mu(\gamma')\ge n \ge 2$ for 
 all closed Reeb orbits $\gamma'$~\cite[Theorem~3.1~(III), 
 Lemma~4.2]{Y}. Thus condition~$(A)$ is empty and condition~$(B_a)$ 
 has to be checked only for cylinders without punctures. It is proved 
 in~\cite[Lemma~7.5]{Y} that, for each bound $\alpha$ on the contact  
 action, there exists a symplectic structure on $W$ such that 
 condition~$(B_0)$ is satisfied for cylinders whose asymptotes have 
 action at most $\alpha$. This is enough in order to define linearized 
 contact homology in the trivial free homotopy class. 
 
\item negative disc bundles $W$ over symplectically aspherical manifolds 
  $(B,\beta)$, for the trivial free homotopy class. 
  Here $\cL\stackrel 
  \pi\to B$ is a Hermitian line bundle  
  with $c_1(\cL)= -[\beta]$ and $W=\{v\in\cL \ : \ |v|\le 1\}$.  
  The symplectic form is $\widehat \omega=\pi^*\beta+d(r^2\theta)$, 
  where $r$ is the radial coordinate in the fibers and $\theta$ is the 
  angular form, and we choose compatible almost complex structures 
  $J_B$ on $B$, $J$ on $\widehat W$ such that $\pi$ is 
  $(J,J_B)$-holomorphic and $J$ satisfies~\eqref{eq:J}.  
  The Reeb orbits on $M:=\p W$ are the circles in the 
  fibers and, to achieve nondegeneracy, we choose a Morse function 
  $f:B\to\R$, we perturb $\widehat \omega$ to $\widehat 
  \omega_\eps=\widehat \omega +\eps f d(r^2\theta)$, $\eps>0$ and  
  $J$ to $J_\eps$ so that $J_\eps$ satisfies~\eqref{eq:J} for the 
  perturbed Reeb vector field.  
  For each bound $\alpha$ on the contact 
  action, there exists $\eps>0$ such that the closed Reeb orbits with 
  period $\le\alpha$ are the  
  circles in the fibers over $\mathrm{Crit}(f)$. We claim that 
  $J_\eps$ satisfies conditions $(A)$ and $(B_0)$ for curves with 
  asymptotes of period $\le\alpha$, which is enough for our purposes.  
 
  The main point is that all punctured holomorphic curves in $(\widehat 
  W,J)$ and $(M\times \R,J)$ are contained in the fibers of $\cL$, due 
  to the asphericity of $B$. The index of a curve with one positive 
  puncture and $m$ negative punctures is $\dim(B)+2(m-1)$, and a direct 
  computation shows that they are all regular (obvious elements in the 
  kernel of the linearized operator correspond to varying the 
  basepoint and the $m-1$ ramification points). The Morse-Bott 
  analysis in~\cite[Chapter~5]{B} shows that, after a slight 
  perturbation of the evaluation  
  maps on the above moduli spaces, regularity still holds after  
  $\eps$-perturbation. The $J_\eps$-holomorphic curves involved in 
  condition~$(A)$ are the fibers over minima of $f$, and the 
  $J_\eps$-holomorphic curves involved in~$(B_0)$ are trivial 
  cylinders contained in the fibers over $\mathrm{Crit}(f)$, cylinders 
  over gradient trajectories of $f$ of index $1$ or $2$, and once  
  punctured cylinders contained in the fibers over the minima of 
  $f$.  
 

 
\item unit cotangent bundles 
  $W=DT^*L = \{v\in T^*L \ : \ |v|\le 1\}$ of closed Riemannian 
  manifolds $L$ such that~$(A)$ is empty, for a free homotopy class $a$ 
  containing only simple closed geodesics. Condition $(A)$ is empty if 
  either $\dim\,L\ge 4$ or $L$ carries no  
  contractible closed geodesics. In the first case the lift 
  $\gamma'_\alpha$  of a closed geodesic $\alpha$ satisfies $\bar 
  \mu(\gamma'_\alpha)=\ind_{\mathrm{Morse}}(\alpha) + \dim\,L-3\ge 
  1$, whereas the second case includes manifolds with strictly negative 
  sectional curvature and flat tori. Condition~$(B_a)$ is 
  satisfied because any cylinder asymptotic at $\pm\infty$ to a simple 
  orbit is somewhere injective, and the set of almost complex 
  structures which are regular for somewhere injective curves is of 
  the second Baire category in $\cJ(\lambda)$~\cite[Theorem~1.8]{D}.  
\end{enum} 
\end{remark}

\subsection{A non-equivariant construction} \label{sec:nonequiv} 
 
In this section we apply the Morse-Bott formalism of~\cite{BOauto} in 
the context of contact homology. The result is a chain complex closely 
related to the non-equivariant contact homology construction 
of~\cite{CO}.  
 
For $\gamma'\in \cP_\lambda$ we denote by $S'_\gamma$ the circle of 
parametrized closed Reeb orbits representing $\gamma'$. Given 
$\og',\ug',\gamma'_1,\ldots,\gamma'_k\in \cP_\lambda$ and $B\in 
H_2(M;\Z)$ we denote by  
$$ 
\widehat \cM^B(S'_\og,S'_\ug,\gamma'_1,\ldots,\gamma'_k;J_\infty)  
$$ 
the {\bf space of punctured $S^1$-parametrized $J_\infty$-holomorphic\break 
  cylinders}, consisting of tuples $(u,L_1,\dots,L_k)$ such that 
$$ 
u=(f,a):\R\times S^1\setminus \{z_1,\ldots,z_k\} \to M \times \R 
$$ 
satisfies  
$$ 
\p_s u + J_\infty \p_\theta u =0, 
$$ 
the $L_i$'s are asymptotic markers at the punctures $z_i$,  
and we require  
$$ 
\lim_{s\to \pm\infty} a(s,\theta)=\mp\infty, \quad \lim_{s\to -\infty} 
f(s,\cdot)\in S'_\og, \quad \lim_{s\to\infty} u(s,\cdot) \in S'_\ug, 
$$ 
as well as~\eqref{eq:Fpunc},~\eqref{eq:Lpunc}, and~\eqref{eq:CRA} at  
$z_1 \ldots, z_k$. The asymptotic     
conditions on $f$ have to be understood in the sense of~\eqref{eq:f}.  
Note that   
the space of punctured $J_\infty$-holomorphic cylinders $\widehat 
\cM^B(\og',\ug',\gamma'_1,\ldots,\gamma'_k;J_\infty)$ defined in the 
previous section consists of $\kappa_{\og'}\kappa_{\ug'}$ copies of  
the space $\widehat 
\cM^B(S'_\og,S'_\ug,\gamma'_1,\ldots,\gamma'_k;J_\infty)$.  
 
The group 
$\R$ acts on $\widehat 
\cM^B(S'_\og,S'_\ug,\gamma'_1,\ldots,\gamma'_k;J_\infty)$ by 
translations in the $s$-variable and we denote by  
$$ 
\cM^B(S'_\og,S'_\ug,\gamma'_1,\ldots,\gamma'_k;J_\infty) := 
\widehat \cM^B(S'_\og,S'_\ug,\gamma'_1,\ldots,\gamma'_k;J_\infty)/\R 
$$ 
the {\bf moduli space of punctured $S^1$-parametrized 
  $J_\infty$-holomor-\break phic 
  cylinders}.  
 
Since the almost complex structure $J_\infty$ satisfies 
assumption~$(B_a)$, the moduli space   
$\cM^B(S'_\og,S'_\ug,\gamma'_1,\ldots,\gamma'_k;J_\infty)$, 
$\og',\ug'\in i^{-1}(a)$   
is a smooth manifold of dimension  
$$ 
\mu(\og')-\mu(\ug') +2\langle c_1(\xi),B\rangle -  
\sum_{i=1}^k \overline \mu(\gamma'_i) +1.     
$$ 
This differs by $1$ from the dimension of the moduli space  
$\cM^B(\og',\ug',\gamma'_1,\ldots,\gamma'_k;J_\infty)$,  
because the $S^1$-symmetry is now broken.    
 
If $\og'\neq \ug'$ the group $\R$ acts freely on the above moduli 
space by translations in the $t$ direction, and we denote the quotient by 
$\cM^B(S'_\og,S'_\ug,\gamma'_1,\ldots,\gamma'_k;J_\infty)/\R$.  
Each such quotient is equipped with smooth 
evaluation maps $\oev$, $\uev$ with target $S'_\og$, respectively 
$S'_\ug$.   
 
Given $\gamma'\in \cP_\lambda$ we choose a Morse function 
$f'_\gamma:S'_\gamma \to \R$ having two critical points $m'$ and $M'$.  
We denote by $\phi_{f'_\gamma}$ the flow of $\nabla f'_\gamma$ with 
respect to some Riemannian metric on $S'_\gamma$. For 
$p'\in\textrm{Crit}(f'_\gamma)$ we denote the Morse index by  
$$ 
\ind(p'):= \dim W^u(p';f'_\gamma). 
$$ 
We denote by $\gamma'_p$ the Reeb orbit in $S'_\gamma$ which corresponds 
to the critical point $p'\in\textrm{Crit}(f'_\gamma)$. We define the 
grading by  
\begin{equation} \label{eq:gr} 
|\gamma'_p| := |\gamma'| + \ind(p').  
\end{equation}  
 
Let $\og',\ug',{\gamma'}^1_1, \ldots, {\gamma'}^1_{k_1}, \ldots, 
{\gamma'}^m_1, \ldots, {\gamma'}^m_{k_m}  
\in\cP_\lambda$, $p'\in \mathrm{Crit}(f'_{\og})$, $q'\in \mathrm{Crit}(f'_{\ug})$  
and $B\in H_2(M;\Z)$.  
For $m\ge 0$ we denote by  
\begin{equation*}  
\cM_m^B(p',q',{\gamma'}^1_1, \ldots, {\gamma'}^1_{k_1}, \ldots, {\gamma'}^m_1, 
\ldots, {\gamma'}^m_{k_m};\{f'_\gamma\},J_\infty) 
\end{equation*} 
the union over    
 ${\tgamma'}_1,\ldots,{\tgamma'}_{m-1}\in\cP_\lambda$ and 
$B_1+\ldots +B_m=B$ of the fibered products   
\begin{eqnarray} 
&& 
\hspace{-.9cm} 
W^u(p') \times _{\oev} ((\cM^{B_1}(S'_\og,S'_{\tgamma_1},{\gamma'}^1_1, 
\ldots, {\gamma'}^1_{k_1})/\R) \ \times \ \R^+) \nonumber \\ 
&& 
\hspace{-.7cm} 
{_{\varphi_{f_{\tgamma_1}}\!\circ\uev}}\!\times   
_{\oev} 
((\cM^{B_2}(S'_{\tgamma_1},S'_{\tgamma_2},{\gamma'}^2_1, \ldots, 
{\gamma'}^2_{k_2})/\R) \ \times \ \R^+) 
\label{eq:fiberedprod} \\ 
&&  
\hspace{-.7cm} 
\ldots\, 
{_{\varphi_{f_{\tgamma_{m-1}}}\!\!\circ\uev}}\!\!\times 
_{\oev}   
\cM^{B_m}(S'_{\tgamma_{m-1}},\!S'_{\ug},{\gamma'}^m_1, \ldots, 
{\gamma'}^m_{k_m})/\R \  {_{\uev}\times} 
W^s(q') .  \nonumber 
\end{eqnarray} 
In this fibered product, the factors $\R^+$ play the role of flow times 
for the maps $\varphi_{f_{\tgamma_1}}, \ldots, \varphi_{f_{\tgamma_{m-1}}}$.  
By our transversality assumptions, this is a smooth    
manifold of dimension  
\begin{eqnarray*} 
\lefteqn{\dim \, \cM_m^B(p',q',{\gamma'}^1_1, \ldots, {\gamma'}^1_{k_1}, \ldots,  
{\gamma'}^m_1, \ldots, {\gamma'}^m_{k_m};\{f'_\gamma\},J_\infty)} \\ 
& = & \ind(p') - 1 + (\dim \, 
\cM^{B_1}(S'_\og,S'_{\tgamma_1},{\gamma'}^1_1, \ldots, 
{\gamma'}^1_{k_1};J_\infty)/\R +1) - 1 \\  
&&+  \ ... +\dim \,  
\cM^{B_m}(S'_{\tgamma_{m-1}},S'_{\ug},{\gamma'}^m_1, \ldots, 
{\gamma'}^m_{k_m};J_\infty)/\R  \\   
&&  - 1 + (1-\ind(q')) \\ 
& = & \mu(\og') - \mu(\ug') +\ind(p') - \ind(q') -1 \\ 
& & + \ 2\langle 
c_1(\xi),B_1+...+B_m\rangle  
- \sum_{i=1}^m \sum_{j=1}^{k_i} \bar\mu({\gamma'}^i_j) \\ 
& = & |\og'_p| - |\ug'_q| - \sum_{i=1}^m \sum_{j=1}^{k_i} \bar\mu({\gamma'}^i_j)  
+ 2 \langle c_1(\xi),B\rangle -1. 
\end{eqnarray*}  
We denote  
\begin{eqnarray*} 
\lefteqn{\cM^B(p',q',{\gamma'}_1, \ldots, {\gamma'}_k;\{f'_\gamma\},J_\infty)} \\ 
&:=& \bigcup_{m\ge 0} 
\cM_m^B(p',q',{\gamma'}^1_1, \ldots, {\gamma'}^1_{k_1}, \ldots,  
{\gamma'}^m_1, \ldots, {\gamma'}^m_{k_m};\{f'_\gamma\},J_\infty) 
\end{eqnarray*}  
with $\{ {\gamma'}^1_1, \ldots, {\gamma'}^1_{k_1}, \ldots,  
{\gamma'}^m_1, \ldots, {\gamma'}^m_{k_m} \} = \{ {\gamma'}_1, \ldots, {\gamma'}_k \}$, 
and we call this {\bf the moduli space of punctured $S^1$-parametrized 
  broken $J_\infty$-holomorphic cylinders}, whereas 
$$ 
\cM_m^B(p',q',{\gamma'}^1_1, \ldots, {\gamma'}^1_{k_1}, \ldots,   
{\gamma'}^m_1, \ldots, {\gamma'}^m_{k_m};\{f'_\gamma\},J_\infty) 
$$ 
will be called {\bf 
the moduli space of punctured $S^1$-parametrized broken 
$J_\infty$-holomorphic cylinders with $m$ 
sublevels}. We refer to Figure~\ref{fig:Phi}.(c) on 
page~\pageref{fig:Phi} for a representation of the elements of these 
moduli spaces (see also the capped version at the end of this section).  
 
The compactified moduli spaces $\overline 
\cM^B(p',q',{\gamma'}_1, \ldots,  
{\gamma'}_k;\{f'_\gamma\},J_\infty)$ admit a topology which is 
described similarly to 
that of $\overline 
\cM^A(p,q;H,\{f_\gamma\},J)$  
in Section~\ref{sec:MBsymp}.  
The boundary (and corners) of  
$\overline \cM^B(p',q',{\gamma'}_1, \ldots,{\gamma'}_k;\{f'_\gamma\},J_\infty)$ 
contain configurations which involve broken gradient trajectories, like 
for $\overline \cM^A(p,q;H,\{f_\gamma\},J)$, but 
also configurations which involve  
punctured $S^1$-parametrized broken $J_\infty$-holomorphic cylinders 
and genus zero $J_\infty$-holomorphic buildings of arbitrary height 
(in the sense of~\cite[Section~7.2]{BEHWZ}), with exactly one positive 
puncture.   
 
\begin{remark} {\rm  
Since we are considering parametrized Reeb orbits $\og'$, $\ug'$, all 
moduli spaces under consideration are orientable provided  
that the asympotes $\gamma'_1,\ldots,\gamma'_k$ are good. In the 
sequel we shall restrict ourselves to such asymptotes at the 
punctures.  
Then one can associate 
a sign $\epsilon(F)$ to each $[F]\in 
\cM^B(p',q',\gamma'_1,\dots,\gamma'_k;\{f'_\gamma\},J_\infty)$ by 
using coherent orientations~\cite{BM} and the fiber-sum 
rule~\cite[\S4]{BOauto}.   
} 
\end{remark}

Let $a$ be a free homotopy class of loops in $W$. We 
define a differential complex, which we call the {\bf 
  $S^1$-parametrized contact complex}, by setting   
$$ 
BC^{i^{-1}(a)}_*(\lambda) := \bigoplus 
_{\gamma'\in\cP^{i^{-1}(a)}_\lambda}  
\Lambda_\om \langle  
\gamma'_m,\gamma'_M\rangle.  
$$ 
We define the {\bf $S^1$-parametrized contact differential}  
$$ 
\delta:BC^{i^{-1}(a)}_*(\lambda) \to 
BC^{i^{-1}(a)}_{*-1}(\lambda)  
$$ 
by 
\begin{equation} 
   \label{eq:delta} 
\delta \og'_p := \sum 
  _{\substack{{\ug'_q},{\gamma'}_1,\ldots,{\gamma'}_k,B \\  
  |{\ug'_q} e^B| + \sum |{\gamma'}_i| = |{\og'_p}| -1 } } 
\hspace{-.5cm} \frac {{\scriptstyle 
    n^B(p',q',\gamma'_1,\ldots,\gamma'_k;\{f'_\gamma\},J_\infty)}} 
{{\scriptstyle \prod_{i=1}^k \kappa_{\gamma_i'}}}  
\ e({\gamma'}_1)\ldots e({\gamma'}_k) e^B {\ug'_q}, 
\end{equation} 
where 
$$ 
n^B(p',q',\gamma'_1,\ldots,\gamma'_k;\{f'_\gamma\},J_\infty) 
:= \sum_{[F]\in 
  \cM^B(p',q',{\gamma'}_1,\ldots,{\gamma'}_k;\{f'_\gamma\},J_\infty)}  
\epsilon(F). 
$$ 
The classes $B$ in the above sum live in    
$H_2(M;\Z)$, but nevertheless the coefficient in front of ${\ug'}_q$ is an 
element of $\Lambda_\om$ due to the factors $e({\gamma'}_i)$, 
$i=1,\ldots,k$. Moreover, the fraction in~\eqref{eq:delta} is an 
integer.  
 
The map $\delta$ satisfies 
$\delta\circ \delta=0$. As explained above, the boundary of  
the $1$-dimensional moduli spaces $\cM^B(p',q',{\gamma'}_1, \ldots, 
{\gamma'}_k;\{f'_\gamma\},J_\infty)$ consists of configurations 
with a single broken gradient trajectory, and of configurations 
involving a rigid  
punctured $J_\infty$-holomorphic cylinder. The 
count of configurations  
of the first type using the augmentation $e$ is equal to $\delta\circ 
\delta$, and the count of configurations of the second type using the 
augmentation $e$ involves $e\circ \p=0$ and hence vanishes. As a 
consequence, the composition $\delta\circ \delta$ vanishes as well.  
The homology groups 
$H_*(BC^{i^{-1}(a)}_*(\lambda),\delta)$ are actually isomorphic to 
the non-equivariant contact homology groups of~\cite{CO}.   
 
\begin{remark} {\rm  
Let $\alpha \in \R^+$ be such that $\alpha\notin 
\mathrm{Spec}(M,\lambda)$. One can define subcomplexes  
$$ 
BC^{i^{-1}(a),\le\alpha}_*(\lambda) := \bigoplus 
_{\gamma'\in\cP^{i^{-1}(a),\le\alpha}_\lambda} \Lambda_\om \langle   
\gamma'_m,\gamma'_M\rangle \ \subseteq \ BC_*^{i^{-1}(a)}(\lambda) 
$$ 
and, for $\alpha\to\infty$, we obviously have  
$$ 
\lim_{\stackrel \longrightarrow \alpha} 
H_*(BC_*^{i^{-1}(a),\le\alpha}(\lambda),\delta) = 
H_*(BC_*^{i^{-1}(a)}(\lambda),\delta).  
$$ 
} 
\end{remark}  
 
We can give an alternative description of the $S^1$-parametrized 
contact differential as follows. Given $A\in 
H_2(W;\Z)$, $\og',\ug'\in \cP_\lambda$, $p'\in\mathrm{Crit}(f'_\og)$, 
$q'\in\mathrm{Crit}(f'_\ug)$ we define the {\bf  
  moduli space of capped punctured $S^1$-parametrized broken 
  $J$-holomorphic cylinders}    
$$ 
\cM^A_c(p',q';\{f'_\gamma\},J) 
$$ 
as the set of equivalence classes of pairs $u=(u',u'')$, with 
$u'$ an element of the moduli space 
$\cM^B(p',q',\gamma'_1,\ldots,\gamma'_k; \{f'_\gamma\},  
J_\infty)$,  
$\gamma'_1,\ldots,\gamma'_k\!\in\!\cP_\lambda$, $B\!\in\!H_2(M;\Z)$, and $u''$ a 
collection of $J$-holomorphic planes in $\widehat W$, of total homology 
class $A-B\in H_2(W;\Z)$, and with asymptotics at their positive punctures  
corresponding to     
$\gamma'_1,\ldots,\gamma'_k$. The elements $u'$, $u''$ are endowed 
with asymptotic markers $L'_i$, $L''_i$, $i=1,\dots,k$. Two collections 
of asymptotic markers $(L'_{i,0},L''_{i,0})$, $(L'_{i,1},L''_{i,1})$ 
are equivalent if they satisfy~\eqref{eq:Arg}.  
The dimension of this moduli space is   
$$ 
|\og'_p|-|\ug'_q| + 2\langle c_1(TW),A\rangle-1. 
$$ 
We refer to Figure~\ref{fig:Phi}.(c) on page~\pageref{fig:Phi} for a 
pictorial representation of these objects. The $S^1$-parametrized 
contact differential can then be rewritten as 
\begin{equation} \label{eq:S1paramdiffcapped} 
 \delta \og'_p = \sum _{\substack{\ug'_q,A \\ 
  |\ug'_q e^A| = |\og'_p| -1 } } 
\sum_{u\in \cM^A_c(p',q';\{f'_\gamma\},J)} 
  \epsilon(u) e^A \ug'_q. 
\end{equation}  
The sign $\epsilon(u)$ is defined as the product of the signs of the 
components of $u$.  
 
\begin{remark} \label{rmk:capped}  
{\rm  
One can define in an obvious manner moduli spaces  
$$ 
\cM^A_c(S'_\og,S'_\ug;J) 
$$  
of capped punctured $S^1$-parametrized 
$J_\infty$-holomorphic cylinders, having dimension $\mu(\og') - \mu(\ug') 
+ 2\langle c_1(TW),A \rangle + 1$, so that the moduli spaces 
$\cM^A_c(p',q';\{f'_\gamma\},J)$ are obtained from 
$\cM^A_c(S'_\og,S'_\ug;J)$ via a fibered product construction  
analogous to~\eqref{eq:fiberedprod}.     
} 
\end{remark}

 
\section{Filtrations} \label{sec:MBcomplex} 
 
This section is organized as follows. We first exhibit natural 
filtrations on the Morse-Bott complex for symplectic homology, as well 
as on the $S^1$-parametrized contact complex. We then describe the 
differential on the $0$-th page of the associated spectral 
sequences. The crucial result in that direction is 
Proposition~\ref{prop:reg}. We then define an isomorphism 
of filtered complexes in Proposition~\ref{prop:filtered}, and explain 
the main steps for the proof of Theorem~\ref{thm:intro}.  
 
Let $a$ be a free homotopy class of loops in $W$. For $a\neq 0$ or 
$a=+$ we define  
\begin{equation}  \label{eq:Bk} 
B_k C_*^a(H) := 
\bigoplus _{\substack{ 
  \gamma\in \cP^a(H) \\ 
  A\in H_2(W;\Z)\\ 
  \mu(\gamma)-2\langle c_1(TW),A\rangle =k }} 
\langle e^A\gamma_m,e^A\gamma_M\rangle. 
\end{equation} 
We claim that 
$$ 
\oplus _{k\le \ell} B_kC_*^a(H), \qquad \ell\in \Z 
$$ 
  forms an increasing filtration on $BC_*^a(H)$. 
Indeed, the sum~(\ref{eq:MBdiff_gamma}) involves 
elements $\og_p,\ug_q$ such that $\mu(\og)+\ind(p)-\mu(\ug)-\ind(q) 
+2\langle c_1(TW),A\rangle =1$. Since 
$$ 
\ind(p)-\ind(q)\in\{-1,0,1\} 
$$ 
it follows that 
$$ 
\mu(\og)-\mu(\ug)+2\langle c_1(TW),A\rangle \in \{0,1,2\}. 
$$ 
As a consequence, the differential 
$d$ on $BC_*^a(H)$ can be written as 
\begin{equation} 
   \label{eq:MBsplit} 
   d = d^0 + d^1 + d^2, 
\end{equation} 
with 
$$ 
d^i:B_k C_*^a(H) \to B_{k-i} C_*^a(H), \qquad i=0,1,2. 
$$

Similarly, we define a filtration on $BC_*^{i^{-1}(a),\le\alpha}(\lambda)$  
\begin{equation} \label{eq:Bk'} 
B_k C_*^{i^{-1}(a),\le\alpha}(\lambda) :=  
\bigoplus _{\substack{ 
  \gamma'\in \cP_\lambda^{i^{-1}(a),\le\alpha} \\ 
  A\in H_2(W;\Z)\\ 
  \mu(\gamma')-2\langle c_1(TW),A\rangle =k }} 
\langle e^A\gamma'_m,e^A\gamma'_M\rangle. 
\end{equation}  
The same argument as above shows that  
$$ 
\oplus _{k\le \ell} B_kC_*^a(\lambda), \qquad \ell\in \Z 
$$ 
forms an increasing filtration on 
$BC_*^{i^{-1}(a),\le\alpha}(\lambda)$, and the differential $\delta$ 
can be written as  
$$ 
\delta = \delta^0 + \delta^1 + \delta^2 
$$ 
with 
$$ 
\delta^i:B_k C_*^{i^{-1}(a),\le\alpha}(\lambda) \to B_{k-i} 
C_*^{i^{-1}(a),\le\alpha}(\lambda), \qquad i=0,1,2.  
$$ 
 
 
\begin{remark}[on the pictorial representation of the moduli spaces] {\rm  
We have drawn in Figure~\ref{fig:Phi} on page~\pageref{fig:Phi} 
elements of the moduli spaces involved in $\delta^i$, $i=0,1,2$. The 
justification for $\delta^0$ and $\delta^1$ is given below, whereas 
the justification for $\delta^2$ is given in Section~\ref{sec:D}. It 
is a consequence of Proposition~\ref{prop:reg} and Remark~\ref{rmk:d0} 
that we can choose the Hamiltonian $H$ so that the same type of moduli 
spaces are involved also in the differentials $d^i$, $i=0,1,2$.   
} 
\end{remark}  
 
We now further elaborate on the differential $\delta^0$, corresponding to 
orbits $\og',\ug'\in\cP_\lambda^{i^{-1}(a),\le \alpha}$ such that  
$\mu(\og')-\mu(\ug')+2\langle c_1(TW),A\rangle =0$.  
The decomposition~\eqref{eq:Bk'} of 
$B_kC_*^{i^{-1}(a),\le\alpha}(\lambda)$  
as a direct sum 
induces a splitting of the differential $\delta^0$ as  
$$ 
\delta^0=\sum_{\gamma'\in\cP_\lambda^{i^{-1}(a),\le \alpha}} \delta^0_{\gamma'}, 
$$ 
with $\im(\delta^0_{\gamma'}) \subset \Lambda_\om \langle 
\gamma'_m,\gamma'_M \rangle$.  
 
We claim that 
$\delta^0_{\ug'}(\og')=0$ if $\ug'\neq \og'$. Let us recall for that 
purpose the moduli spaces of capped punctured $S^1$-parametrized 
$J_\infty$-holomorphic cylinders $\cM^A_c(\og'_p,\ug'_q;J)$ defined in 
Remark~\ref{rmk:capped}. Each such space necessarily has dimension  
$1$. However, it carries a nontrivial $S^1$-action, as well as a  
free $\R$-action on the target if $\og'\neq\ug'$, in which case it 
must be empty. As a consequence, we must have 
$\delta^0(\gamma')=\delta^0_{\gamma'}(\gamma')$. Since the difference 
of actions of the asymptotes at $\pm\infty$ for elements of 
$\cM^A_c(S'_\og,S'_\ug;J_\infty)$ is zero, the latter spaces are nonempty 
only if $A=0$ and there are no punctures. Hence $\delta^0$ counts 
gradient trajectories and, in particular, we have 
$\delta^0(\gamma'_M)=0$.  
The next result is a straightforward adaptation of Lemma~4.25 
in~\cite{BOauto}.  
 
\begin{proposition} \label{prop:delta0} 
 If $\gamma'\in\cP_\lambda^{i^{-1}(a),\le\alpha}$ is a good orbit, 
 then  
$$ 
\delta^0(\gamma'_m)=0. 
$$ 
If $\gamma'\in\cP_\lambda^{i^{-1}(a),\le\alpha}$ is a bad orbit, then  
$$ 
\delta^0(\gamma'_m)=\pm 2\gamma'_M. 
$$ 
\end{proposition}  
 
Let us now discuss the differential $d^0$, which corresponds to 
orbits $\og,\ug\in\cP^a(H)$ such that $\mu(\og)-\mu(\ug)+2\langle 
c_1(TW),A \rangle = 0$. As above, it also splits as  
$$ 
d^0 = \sum_{\gamma\in\cP^a(H)} d^0_\gamma 
$$ 
with $\textrm{im}(d^0_\gamma) \subset \Lambda_\om \langle 
\gamma_m,\gamma_M \rangle$.  
 
One important situation is $\og=\ug$. In this case the moduli spaces 
of Floer trajectories $\cM^A(\og,\og;H,J)$, $A\neq 0$ are empty, 
whereas the space $\cM^0(S_\og,S_\og;H,J)$ consists of constant cylinders 
and is naturally parametrized by $S_{\og}$. The Morse-Bott 
differential $d^0_\og(\og_p)$, $p \in {\rm Crit}(f_\og)$ is given  
by a count of gradient trajectories and therefore  
$$ 
 d^0_\og(\og_M)=0, 
$$ 
while $d^0_\og(\og_m)$ is a multiple of $\og_M$. The next statement 
is a reformulation of Lemma~4.25 in~\cite{BOauto}.  
 
\begin{proposition} \label{prop:d0}  
If $\og\in\cP_\lambda^{\le \alpha}$ is a good orbit then 
$$ 
  d^0_\og \og_m =0. 
$$ 
If $\og\in \cP_\lambda^{\le\alpha}$ is a bad orbit then 
$$ 
  d^0_\og \og_m = \pm 2 \og_M. 
$$ 
\end{proposition} 
 
We shall describe in Section~\ref{sec:slow} a 
procedure for slowing down the Hamiltonian $H$ which produces a family 
$H^R$, $R>0$.  
 
\begin{proposition} \label{prop:reg} 
 Assume the almost complex structure $\widehat J$ on $\widehat W$ satisfies 
 conditions~$(A)$ and~$(B_a)$ in Section~\ref{sec:contact}. If $R>0$  
 is large enough, then we have  
$\cM^A(S_\og, S_\ug;H^R,\widehat J)=\emptyset$ if 
 $\og\neq \ug\in\cP_\lambda^{\le\alpha}$, $\og,\ug\in i^{-1}(a)$   
and  $\mu(\og)-\mu(\ug)+2\langle c_1(TW),A\rangle =0$. 
\end{proposition}  
 
The proof of Proposition~\ref{prop:reg} is given at the end of  
Section~\ref{sec:slow}.  
 
\begin{remark} \label{rmk:d0} {\rm  
By a limiting argument, one sees that the conclusion of 
Proposition~\ref{prop:reg} still holds if $J\in\Jreg(H^R)$ is a small 
time-dependent perturbation of $\widehat J$. As a consequence, one can 
find regular almost complex structures for which  
$d^0(\og_p)=d^0_\og(\og_p)$, $\og\in\cP_\lambda^{\le 
\alpha}$, $p \in {\rm Crit}(f_\og)$ and Proposition~\ref{prop:d0} 
holds with $d^0_\og$ replaced by $d^0$.  
} 
\end{remark}  
 
\begin{remark} \label{rmk:sdep} {\rm  
We assume in this remark that the almost complex structure  
 $\widehat J$ on $\widehat W$ satisfies  
 conditions~$(A)$ and~$(B_a)$ in Section~\ref{sec:contact}.  
Let $H_-^R\le H_+^R$ be slow Hamiltonians as constructed in 
Section~\ref{sec:slow}, and let $H_s^R$, $s\in\R$ be an increasing 
homotopy. If $R>0$ is large enough 
 and if $\p_sH_s^R$ is small enough in $C^0$-norm, we have 
 $\cM^A(S_\og,S_\ug;H_s^R,\widehat J)=\emptyset$ for  
 $\og\in\cP^{\le\alpha}(H_-^R)$, 
 $\ug\in\cP^{\le\alpha}(H_+^R)$ such that $\og,\ug\in i^{-1}(a)$ and  
$\mu(\og)-\mu(\ug)+2\langle c_1(TW),A\rangle=-1$.  
This is essentially a consequence of the fact that the homotopy 
 $H_s^R$ is independent of $\theta$, and is proved using the arguments 
 for the proof of Proposition~\ref{prop:reg} in Section~\ref{sec:slow} 
 (the dimensions of the corresponding 
 moduli spaces have to be shifted by $1$ due to the $s$-dependence of 
 the equation). Moreover, the conclusion of Remark~\ref{rmk:d0} above  
 continues to hold in the $s$-dependent situation.  
 
This implies that the continuation morphism  
 $$ 
 \sigma_{H_+^R,H_-^R}:BC_*^a(H_-^R)\to BC_*^a(H_+^R) 
 $$ 
defined in Section~\ref{sec:MBsymp} preserves the filtrations.  
Indeed, this map 
preserves the degree and is obtained by a count of rigid configurations in  
 $\cM^A(p,q;H_s^R,\{f_\gamma^\pm\},J_s)$, which can 
 be of the following types:  
\begin{enumerate} 
 \item \label{type:1} $\cM^A(m,M;H_s^R,\{f_\gamma^\pm\},J_s)$, with 
$\mu(\og)-\mu(\ug)+2\langle c_1(TW),A\rangle = -1$; 
 \item $\cM^A(m,m;H_s^R,\{f_\gamma^\pm\},J_s)$ or 
$\cM^A(M,M;H_s^R,\{f_\gamma^\pm\},J_s)$, such that  
$\mu(\og)-\mu(\ug)+2\langle c_1(TW),A\rangle = 0$; 
 \item $\cM^A(M,m;H_s^R,\{f_\gamma^\pm\},J_s)$, with 
$\mu(\og)-\mu(\ug)+2\langle c_1(TW),A\rangle = 1$.  
\end{enumerate} 
Our discussion shows that rigid configurations of 
Type~\ref{type:1} do not exist, hence $\sigma_{H_+^R,H_-^R}$ preserves the 
filtrations.   
} 
\end{remark}  
 
We shall describe in Section~\ref{sec:stretch} a procedure for 
stretching the neck in the neighbourhood of $\p W$ which  
produces a deformed almost complex structure $J^\tau$ and a 
deformed Hamiltonian $H^{\tau,R}$, $\tau>0$ on $\widehat W$, so that  
$J^\tau \in\Jreg(H^{\tau,R})$, the conclusion of 
Proposition~\ref{prop:reg} still holds and the Floer trajectories of 
$H^{\tau,R}$ are ``close''  
to punctured Floer trajectories in the symplectization, capped with 
rigid holomorphic planes in $\widehat W$. 
>From now on we denote by $H$, $J$ the Hamiltonian $H^{\tau,R}$ and 
the almost complex structure $J^\tau$ with the above properties.

Let $(E^r_d,\bar d^r)$ and $(E^r_\delta,\bar \delta^r)$, $r\ge 0$ be the spectral sequences 
corresponding to the above filtrations on the complexes $(BC_*^a(H),d)$ and 
$(BC_*^{i^{-1}(a),\le\alpha}(\lambda),\delta)$ respectively. 
As a consequence of the previous discussion on $\delta^0$ and $d^0$ 
we infer that 
\begin{equation} \label{eq:E^1delta} 
E^1_\delta:=H_*(BC_*^{i^{-1}(a),\le\alpha}(\lambda),\delta^0) = 
C_*^{i^{-1}(a),\le \alpha}(\lambda)  
\otimes H_*(S^1) 
\end{equation} 
and  
\begin{equation} \label{eq:E^1} 
E^1_d:=H_*(BC_*^a(H),d^0) \simeq C_*^{i^{-1}(a),\le \alpha}(\lambda) 
\otimes H_*(S^1) 
\end{equation} 
as $\Lambda_\om$-modules. As a consequence we have $E^1_\delta\simeq 
E^1_d$ as $\Lambda_\om$-modules. The next statement implies in 
particular that we have an isomorphism of differential complexes 
$(E^1_\delta,\bar \delta^1)\simeq (E^1_d,\bar d^1)$.

\begin{proposition} \label{prop:filtered}  
 There is an isomorphism of filtered complexes  
$$ 
\Phi:BC_*^{i^{-1}(a),\le\alpha}(\lambda) \stackrel \sim \to BC_*^a(H) 
$$ 
which decreases the degree by $n-3$.  
\end{proposition}  
 
The proof of Proposition~\ref{prop:filtered} is given in 
Section~\ref{sec:filtered}.  
 
\begin{corollary} \label{cor:spectral}  
 The map $\Phi$ induces an isomorphism of spectral sequences 
$$ 
(E^r_\delta,\bar \delta^r) \stackrel \sim \to (E^r_d,\bar d^r), \qquad r\ge 0. 
$$ 
\end{corollary}  
 
The differential $\bar \delta^1$ is closely connected to the 
differential $\p$ for linearized contact homology. More precisely,  
$\bar \delta^1$ is given by a count of elements in moduli spaces   
$\cM^A_c(p',q';\{f'_\gamma\},J)$, $p'\in\mathrm{Crit}(f'_\og)$, 
$q'\in\mathrm{Crit}(f'_\ug)$ with $\ind(p')=\ind(q')$, so that 
either $W^u(p')$ or $W^s(q')$ is reduced to one point.  
If $p'=M'$ is a maximum, then $W^u(p')$ is reduced to one point and the  
moduli space $\cM^A_c(p',q';\{f'_\gamma\},J)$ is diffeomorphic to  
$\kappa_{\og'}$ copies of the quotient 
$\cM^A_c(S'_\og,S'_\ug;J)/S^1$. If $q'=m'$ is a minimum, then 
$W^s(q')$ is reduced to one point and the moduli space 
$\cM^A_c(p',q';\{f'_\gamma\},J)$ is diffeomorphic to   
$\kappa_{\ug'}$ copies of the quotient 
$\cM^A_c(S'_\og,S'_\ug;J)/S^1$. Recalling that $\cM^A_c(\og',\ug';J)$ 
is diffeomorphic to $\kappa_{\og'}\kappa_{\ug'}$ copies of 
$\cM^A_c(S'_\og,S'_\ug;J)/S^1$, and denoting by 
$\#\cM^A_c(S'_\og,S'_\ug;J)/S^1$ the algebraic count of elements of 
the moduli space  
$\cM^A_c(S'_\og,S'_\ug;J)/S^1$, we obtain   
$$ 
\p(\og')=\sum_{A,\ug'} \kappa_{\og'} \# 
\cM^A_c(S_{\og'},S_{\ug'};J)/S^1 e^A \ug', 
$$ 
$$ 
\bar \delta^1(\og'_M) = \sum_{A,\ug'} \kappa_{\og'} \# 
\cM^A_c(S_{\og'},S_{\ug'};J)/S^1 e^A \ug'_M, 
$$ 
$$ 
\bar \delta^1(\og'_m) = \sum_{A,\ug'} \kappa_{\ug'} \# 
\cM^A_c(S_{\og'},S_{\ug'};J)/S^1 e^A \ug'_m. 
$$ 
We define an automorphism $\Theta$ of 
$C_*^{i^{-1}(a),\le\alpha}(\lambda)\otimes H_*(S^1)$ by  
\begin{equation} \label{eq:Theta} 
\Theta(\gamma'\otimes M):=\gamma'\otimes M, \qquad \Theta(\gamma'\otimes 
m):=\frac 1 {\kappa_{\gamma'}} \gamma'\otimes m. 
\end{equation}  
Here $M$, $m$ are the generators of $H_0(S^1)$, respectively 
$H_1(S^1)$. Then  
$$ 
\bar \delta ^1 = \Theta^{-1}\circ \p \circ \Theta.  
$$ 
In particular, $\Theta$ induces an isomorphism  
\begin{equation} \label{eq:E2delta} 
 \bar \Theta: E^2_\delta \stackrel \sim \longrightarrow 
   HC_*^{i^{-1}(a),\le\alpha}(\lambda,J) \otimes  
   H_*(S^1). 
\end{equation}  
 
Theorem~\ref{thm:intro} follows from the previous 
considerations and the commutative diagram below, with vertical 
arrows being isomorphisms (see Section~\ref{sec:proofleq} for the full 
details).   
$$ 
\xymatrix 
@C=16pt 
{ 
& & E^2_{\delta;*,0} \ar[r]^-{\bar \delta^2} \ar[d]_{\Phi} & 
E^2_{\delta;*-2,1} \ar[d]^{\Phi} &  
& \\ 
\ldots \ar[r] & SH_*^a(H,J) \ar[r] & E^2_{d;*,0} \ar[r]^-{\bar d^2} & E^2_{d;*-2,1} 
\ar[r] & SH_{*-1}^a(H,J) \ar[r] & \ldots 
} 
$$

 
\section{Floer- and holomorphic 
  cylinders} \label{sec:Floerhol}  
 
The goal of this section is to prove Proposition~\ref{prop:reg} 
and to prepare the proof of Proposition~\ref{prop:filtered}. For 
Proposition~\ref{prop:reg}  
we modify the Hamiltonian to another one with the same 
asymptotic slope but which varies very slowly in $M\times \R$, and 
analyze the limit of Floer trajectories as the rate of variation 
goes to zero. For Proposition~\ref{prop:filtered} we must, loosely 
speaking, confine Floer trajectories in $\widehat W$ near the 
boundary, so that we can view them in the symplectization $M\times 
\R$. In Section~\ref{sec:stretch} we stretch 
the neck near $M$ and show, by a compactness argument, 
that the rigid Floer trajectories connecting nonconstant orbits of 
$X_H$ are in bijective correspondence with punctured solutions of 
Floer's equation in the symplectization $M\times \R$, capped at the 
punctures with rigid holomorphic planes in $\widehat W$.  
 
We recall that we work under the 
standing assumptions $(A)$ and $(B_a)$ in  
Section~\ref{sec:contact}, Remark~\ref{rmk:transv_exples}.

\subsection{Slowing down the Hamiltonian} \label{sec:slow}  
 
We start with the following two perturbative lemmas.    
 
\begin{lemma} \label{lem:shift} 
  Let $c:\R\to \R_+$ be a smooth increasing function with the property 
  that $c'(s)\neq 0$ if $c(s)\in\mathrm{Spec}(M,\lambda)$. Let 
  $u=(f,a):\R\times S^1\setminus  
  \{z_1,\ldots,z_k\}\to M\times \R$ be a solution of the equation  
  \begin{equation}  \label{eq:c}  
   \p_s u + \widehat J_\infty \p_\theta u - c(a(s,\theta)) \frac \p {\p t} =0 
  \end{equation}  
 which, in polar coordinates $(\rho_i,\theta_i)$ around $z_i$, satisfies 
  \begin{equation} \label{eq:asyzi}  
  \lim_{z\to z_i} a(z)=-\infty, \qquad \lim_{\rho_i\to 0} 
  f(\rho_i,\theta_i)=\gamma'_i(T_i \theta_i), 
  \end{equation}  
for some $\gamma'_1,\ldots,\gamma'_k\in\cP_\lambda$, as well as one of 
the following asymptotic conditions at $\pm\infty$.  
\begin{enumerate}  
 \item \begin{equation} \label{eq:asyinfinite} 
        \lim_{s\to \pm\infty} a(s,\cdot)=\mp\infty, \quad 
  \lim_{s\to\pm\infty} f(s,\theta)=\gamma'_\pm(T_\pm\theta),  
       \end{equation}  
 \item \begin{equation} \label{eq:asyinfinitefinite} 
        \lim_{s\to -\infty} a(s,\cdot)=+\infty, \ 
        \lim_{s\to+\infty} a(s,\cdot) = a_0, \  
  \lim_{s\to\pm\infty} f(s,\theta)=\gamma'_\pm(T_\pm\theta),  
       \end{equation}  
 \item \begin{equation} \label{eq:asyfiniteinfinite} 
        \lim_{s\to -\infty} a(s,\cdot)=a_0, \ 
        \lim_{s\to+\infty} a(s,\cdot) = -\infty, \ 
  \lim_{s\to\pm\infty} f(s,\theta)=\gamma'_\pm(T_\pm\theta), 
       \end{equation}  
\end{enumerate}  
for $\gamma'_\pm\in\cP_\lambda^{\le\alpha}$. Given $c_0>0$ we denote 
$u_{c_0}:=(f,a-c_0s)$, so that $u_{c_0}$ satisfies the equation  
\begin{equation} \label{eq:c0} 
\p_s u + \widehat J_\infty \p_\theta u - (c(a(s,\theta))-c_0) \frac \p {\p t} =0. 
\end{equation}  
 
Assume $c_0$ is such that $c'(s)\neq 0$ if $c(s)-c_0\in 
\mathrm{Spec}(M,\lambda)$.  
Then the linearized operator    
$D_{u,c}$ corresponding to~\eqref{eq:c} is surjective if and only if 
the linearized operator $D_{u_{c_0},c-c_0}$ corresponding 
to~\eqref{eq:c0} is surjective.  
\end{lemma}  
 
\proof We first treat case 1. The asymptotic conditions at 
  $z_1,\ldots,z_k$ are obviously preserved, whereas   
   the asymptotic conditions at $\pm\infty$ are still satisfied because 
  $c_0>0$. The two linearized operators have the same expressions, 
  and their domain and target   
  are canonically identified since the components of $u$ and 
  $u_{c_0}$ along $M$ are the same. The conclusion follows.  
 
Cases 2. and 3. are similar so we treat only Case 2.     
There is one asymptotic condition which changes, namely $\lim_{s\to \infty} 
a(s,\cdot)=-\infty$, and this causes a change in the domain of the 
linearized operators. The domain of $D_{u_{c_0},c-c_0}$ contains an 
additional $1$-dimensional summand $V$ corresponding to a degeneracy of 
the asymptote $\gamma'_+$ in the direction $\frac \p {\p t}$. However, 
since~\eqref{eq:c0} is invariant under translation in the $s$-variable, 
we have $\p_s u_{c_0}\in \ker\, D_{u_{c_0},c-c_0}$ and $\p_s u_{c_0}$ 
has a nontrivial component along this additional $1$-dimensional 
summand $V$. Since the quotient $\mathrm{dom}(D_{u_{c_0},c-c_0})/V$ is 
canonically identified with $\mathrm{dom }(D_{u,c})$, the conclusion 
follows.  
\hfill{$\square$}  
 
\begin{lemma} \label{lem:pert}  
  For any $C>0$ there exists $\epsilon=\epsilon(C) >0$ such that, for 
  any function $c:\R\to [0,C]$ with $\mathrm{Supp}(c')\subset 
  [-1,1]$, such that $c'(s)\neq 0$ if $c(s)\in 
  \mathrm{Spec}(M,\lambda)$, and such that  
  $$ 
  0\le c'<\epsilon, 
  $$ 
  the almost complex structure $\widehat J_\infty$ is  
  regular for all solutions $u=(f,a):\R\times S^1\setminus 
  \{z_1,\ldots,z_k\}\to M\times \R$ of equation~\eqref{eq:c}  
  satisfying~\eqref{eq:asyzi} and one of the additional asymptotic 
  conditions~\eqref{eq:asyinfinite}, \eqref{eq:asyinfinitefinite} 
  or~\eqref{eq:asyfiniteinfinite}.  
\end{lemma}  
 
\proof  
  We prove the Lemma by contradiction.  
  Without loss of generality 
  we can assume that there exist sequences $\epsilon_n \to 0$, $c_n$ and  
  $u_n$ as in the statement of the Lemma, where all the $u_n$'s have 
  the same asymptotes, such that the linearized operators 
  $D_{u_n,c_n}$ are not surjective. The sequence $c_n$ converges to 
  some constant $c\in[0,C]$ and we choose $c_0>0$ such that 
  $c-c_0\notin \mathrm{Spec}(M,\lambda)$. In particular, the constant 
  $c_0$ satisfies the conditions of Lemma~\ref{lem:shift} for each 
  $c_n$ if $n$ is large enough. Renaming the sequences $u_{n,c_0}$ and 
  $c_n-c_0$ as $u_n$, $c_n$, the operators $D_{u_n,c_n}$ are not 
  surjective by Lemma~\ref{lem:shift} and the maps $u_n$ satisfy the 
  asymptotic conditions~\eqref{eq:asyzi} 
  and~\eqref{eq:asyinfinite}. We also rename $c-c_0$ as $c$. 
 
  After passing to a subsequence~\cite[Section 10.2]{BEHWZ} (see also the proof of 
  Step 1 in Proposition~\ref{prop:stretch}), $u_n$ converges to a   
  broken curve $u$ whose components either solve $\p_s w+\widehat J_\infty 
  \p_\theta w - c\frac \p {\p t}=0$ or are $\widehat J_\infty$-holomorphic.  
  The linearized operator at each such component is surjective by 
  Lemma~\ref{lem:shift} and the regularity 
  assumption~$(B_a)$ for $\widehat J_\infty$. For 
  $n$ large enough, $u_n$ is approximated by a gluing construction as 
  follows. Let us denote $\dbar_{c_n}:=\p_s + \widehat J_\infty \p_\theta -  
  c_n\frac \p {\p t}$. We replace each holomorphic component $v=(f,a)$ by 
  $v_{c_n}:=(f,a+c_n(-\infty)s)$ and each component $w=(f,a)$ solving 
  $\dbar_c w=0$ by $w_{c_n}:=(f,a+(c_n(a)-c)s)$.  
  The first remark is that $\|\dbar_{c_n}(w_{c_n})\|$ is arbitrarily 
  small for $n$ large enough. Here we use an $L^p$-norm with 
  exponential weights for a metric which has cylindrical ends near the 
punctures $z_i$. These weights do not play any role in estimating 
  $\|\dbar_{c_n}(w_{c_n})\|$ because $\mathrm{supp}(c'_n) \subset 
[-1,1]$ and $\bar\partial_{c_n}(w_{c_n})$ vanishes in a fixed 
neighborhood of the punctures $z_i$ in view of the asymptotic condition 
\eqref{eq:asyzi}. In fact, the function $\bar\partial_{c_n}(w_{c_n})$  
has compact support because it satisfies~\eqref{eq:asyinfinite}.  
Decomposing $T(M\times \R)=\langle \frac \p {\p t} \rangle \oplus 
\langle R_\lambda \rangle \oplus \xi$ we obtain 
\begin{equation*} 
\dbar_{c_n} w_{c_n} \!=\! \left( \begin{array}{c} \p_s(a+(c_n(a)-c)s) - 
    \lambda(\p_\theta)f - c_n(a) \\ 
  \lambda(\p_s f) + \p_\theta (a+(c_n(a)-c)s) \\ 
  0 \end{array} \right)  
\!=\! \left( \begin{array}{c} c_n'(a) s \p_s a \\ 
  c_n'(a) s \p_\theta a \\ 
  0 \end{array} \right), 
\end{equation*}  
and this quantity goes to zero when $n\to\infty$ since it is bounded 
by a constant multiple of $\epsilon_n$. The second remark is  
  that each linearized operator $D_{w_{c_n},c_n}$ is close to the 
  linearized operator $D_{w,c}$ because they again differ by a term 
  involving $c_n'(a)$. Hence each $D_{w_{c_n},c_n}$ is surjective, 
  whereas each linearized operator $D_{v_{c_n},c_n(-\infty)}$ is 
  surjective by the argument in Lemma~\ref{lem:shift}. We now apply 
the standard gluing construction  
to the broken curve with components $v_{c_n}$, $w_{c_n}$ and get a 
solution $\widetilde u_n$ solving $\dbar_{c_n}(\widetilde u_n)=0$ and having 
a surjective linearized operator $D_{\widetilde u_n,c_n}$ with 
uniformly bounded right inverse with respect to $n$.  
Moreover, $\widetilde u_n$ is arbitrarily close to $u_n$ as $n\to\infty$. As a 
conclusion, the linearized operator $D_{u_n,c_n}$ is surjective for 
$n$ large enough, a contradiction.  
\hfill{$\square$}  
 
\medskip  
 
We now describe a family of admissible Hamiltonians $H^R:\widehat W\to 
\R$ such that, on $M\times [0,\infty[$, we have 
$H^R(p,t)=h^R(t)=\rho(t) e^t + \tC$ for some $\tC \in \R$,   
with $\rho(t)=\alpha_0$, $\alpha_0<\min \mathrm{Spec}(M,\lambda)$ if $t$ 
is close to $0$ and $\rho(t)=\alpha$ for $t\ge t_R$ large 
enough.  
 
Let $\epsilon = \epsilon(\alpha)$ be such that Lemma~\ref{lem:pert} holds for 
$C = \alpha$. After possibly choosing a smaller value for 
$\epsilon$, we can assume without loss of generality that, for all 
$T_1\neq T_2\in \mathrm{Spec}(M,\lambda)$ such that $T_1,T_2\le\alpha$ 
we have $|T_1 -T_2|\ge 2\epsilon$. We require the function $h^R$ to 
have the following property:  
 
\begin{enumerate}  
 \item[{}] The interval $[0,t_R]$ is the concatenation of $N = \big[\frac 
 {2\alpha} \epsilon\big]$ intervals $[-\frac R 2, \frac R 2]$ on which the 
 function $c=\rho +\rho'$ satisfies the hypotheses of Lemma~\ref{lem:pert}. 
\end{enumerate}  
 
\proof[Proof of Proposition~\ref{prop:reg}] 
 We proceed by contradiction and we assume there exist sequences 
 $R_n\to \infty$ and $u_n\in \cM^A(S_\og,S_\ug;H^{R_n},\widehat J)$, 
 $\og,\ug\in\cP(H)$  
 with $\og\neq \ug$ and $\mu(\og)-\mu(\ug) + 2\langle c_1(TW),A\rangle 
 =0$. Note that the index of each of the operators $D_{u_n}$ is equal 
 to $1$. After passing to a subsequence \cite[Section 10.2]{BEHWZ}  
 (see also the proof of Step 1 in Proposition~\ref{prop:stretch}),  
 $u_n$ converges to a broken curve whose 
 components in $M\times \R$ are either $\widehat J_\infty$ holomorphic 
 curves, or satisfy~\eqref{eq:c} and condition~\eqref{eq:asyzi}, as 
 well as one of the asymptotic conditions~\eqref{eq:asyinfinite}, 
 \eqref{eq:asyinfinitefinite} 
 or~\eqref{eq:asyfiniteinfinite}. Moreover, the components in 
 $\widehat W$ are $\widehat J$-holomorphic planes. By the regularity 
 assumptions~$(A)$ and~$(B_a)$ on $\widehat J$ and Lemma~\ref{lem:pert} we 
 know that the almost complex structures $\widehat J_\infty$ and 
 $\widehat J$ are regular for all these components. As a consequence, 
 the Fredholm index of the 
 linearized operators at the components  
 satisfying~\eqref{eq:c} is at most $1$, and the same holds for 
 the dimension of their kernel. On the other hand, since $\og\neq 
 \ug$ there is at least one such component which is not a 
 (reparametrized) vertical cylinder. Since equation~\eqref{eq:c} is 
 invariant under reparametrizations in both variables $s$ and 
 $\theta$, we deduce that the kernel of the corresponding linearized 
 operator is at least $2$-dimensional, a contradiction.  
\hfill{$\square$}

\subsection{Stretching the neck near the boundary} \label{sec:stretch} 
 
In this section we assume without loss of generality that $H\in\cH'$ 
is of the form  
$H=\alpha_0e^t+\beta_0$ on a neighbourhood $\p W\times 
[-\epsilon,\epsilon]$ of $\p W$ in $\widehat W$, where 
$0<\alpha_0< \min\,\textrm{Spec}(M,\lambda)$. Moreover, we 
assume that  
\begin{equation} \label{eq:slow}  
|\frac d {dt} (e^{-t}h'(t))| <1. 
\end{equation}  
 
We define a deformation $(\widehat W^\tau,\widehat\om^\tau)$, $\tau  
\ge 0$ by  
\begin{equation*} 
 \widehat W^\tau := W \,\bigcup\, (M \times [-\tau,0]) \,\bigcup\, (M 
 \times \R^+) 
\end{equation*}  
and  
\begin{equation*} 
  \widehat \om^\tau :=\left\{\begin{array}{rl} e^{-\tau}\om, & \mbox{ 
  on } W, \\ 
  d(e^t\lambda), & \mbox{ on } M\times [-\tau,\infty[ 
  .\end{array}\right.   
\end{equation*} 
 
Let $\beta:\R\to[0,1]$ be a smooth function such  
that $\beta(t)=1$ for $t\le -C$, $\beta(t)=0$ for $t\ge 0$,  
where $C>0$ is some large constant.  
We define a family of Hamiltonians $\{H^\tau\}$, $H^\tau:\widehat W^\tau\to \R$ by  
\begin{equation} \label{eq:stretch} 
H^\tau:=\left\{\begin{array}{rl}Êe^{-\tau}H, & \mbox{ on } W, \\ 
(e^{-\tau}\beta(t)+ 1-\beta(t))(\alpha_0e^t+\beta_0), & \mbox{ on } M\times [-\tau,0], \\ 
h(t), & \mbox{ on } M\times \R^+. 
\end{array}\right. 
\end{equation} 
We define $H_\infty:M\times \R\to\R$ by  
$H_\infty:=\lim_{\tau\to\infty} H^\tau|_{M\times [-\tau,\infty[}$. We then have  
$H_\infty(t)=0$ for $t\le -C$ and, if $C$ is large enough, we can choose $\beta$ so that  
the functions $e^{-t}\frac d {dt} H^\tau(t)$, $\tau\ge 0$ are increasing on $[-\tau,\infty[$.  
 
\begin{figure}[hpt] 
\centering 
\input{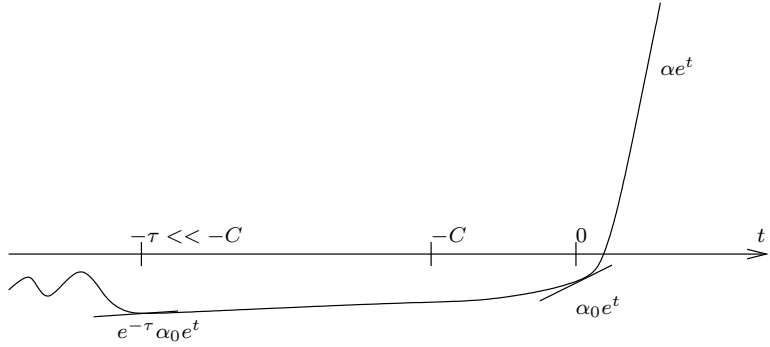} 
\caption{The modified Hamiltonian $H^\tau$ \label{fig:graph}} 
\end{figure}

Let $\widehat J$ be an almost complex structure satisfying the 
regularity assumptions~$(A)$ and~$(B_a)$.  
Similarly to Proposition~3.5 in~\cite{BOauto}, we can choose a small   
time-dependent perturbation $J_\infty$ of $\widehat J_\infty$ which is 
localized   
in an arbitrarily small neighbourhood of the $1$-periodic orbits of $X_H$ and  
away from the orbits themselves, and  
which is regular for punctured Floer trajectories  
\begin{equation*}  
u=(f,a):\R\times S^1\setminus \{z_1,\ldots,z_k\}\to M\times \R 
\end{equation*} 
satisfying  
\begin{equation} \label{eq:F} 
\p_s u + J_\infty(\p_\theta u - X_H(u))=0, 
\end{equation} 
\begin{equation} \label{eq:Fasy} 
\lim_{s\to\pm\infty}u(s,\cdot)=\gamma_\pm(\cdot), \quad \gamma_\pm\in\cP(H), 
\end{equation} 
as well as~\eqref{eq:Fpunc} for some $\gamma'_i\in\cP_\lambda^{\le\alpha}$, 
$i=1,\ldots,k$. We denote by $J$ the resulting almost complex 
structure on $\widehat W$, which coincides with $\widehat J$ on $W$ 
and with $J_\infty$ on $M\times \R^+$.   
 
We define a family $\{\widehat J^\tau\}$, $\tau\ge 0$ of $\widehat
\om^\tau$-compatible   
almost complex structures on $\widehat W^\tau$ by  
\begin{equation*} 
  \widehat J^\tau :=\left\{\begin{array}{rl} \widehat J, & \mbox{ 
  on } W, \\ 
  \widehat J_\infty, & \mbox{ on } M\times [-\tau,\infty[. 
  \end{array}\right.   
\end{equation*} 
As above, we can choose a small time-dependent perturbation $\{J^\tau\}$ of  
the family $\{\widehat J^\tau\}$ which is supported in an arbitrarily
small neighbourhood   
of the $1$-periodic orbits of $X_H$, which satisfies  
$$ 
J^\tau|_{M\times \R^+} \to J_\infty, \quad \tau\to\infty,
$$ 
and such that each $J^\tau$ is regular for $H^\tau$.

The sets $\cP(H)$ and $\cP(H^\tau)$, $\tau\ge 0$ are naturally
identified. These sets also contain constant elements (critical points
of $H$ and $H^\tau$), but we shall tacitly consider in the rest of this
section only nonconstant elements. The same convention applies to
$\cP(H_\infty)$. 

For any $\tau\ge 0$, the set
$\cF_{\mathrm{reg}}(H^\tau,J^\tau)$ of regular collections of perfect
Morse functions $\{f_\gamma^\tau:S_\gamma\to\R, \, \gamma\in\cP(H)\}$
defined in Section~\ref{sec:MBsymp} is of the second Baire
category in the space of collections of perfect Morse
functions~\cite[Lemma~3.6]{BOauto}. Hence, given a sequence 
$\tau_\nu\to\infty$, $\nu\to\infty$, there exists a collection
$\{f_\gamma\}$ which belongs to $\cF_{\mathrm{reg}}(H^{\tau_\nu},J^{\tau_\nu})$ for any
$\nu$. The moduli spaces of
Morse-Bott broken trajectories
$$
\cM^A(p,q;H^{\tau_\nu},\{f_\gamma\},J^{\tau_\nu})
$$
are then
well-defined, for any $\og,\ug\in\cP(H)$, $p\in\mathrm{Crit}(f_\og)$,
$q\in\mathrm{Crit}(f_\ug)$, and $A\in H_2(W;\Z)$.

Our goal in this section is to establish a bijective correspondence
with the moduli spaces $\cM^A_c(p,q;H_\infty,\{f_\gamma\}, J_\infty)$ 
of capped punctured Morse-Bott broken trajectories which we now
define.  

 Given elements $\og,\ug\in\cP(H)$,   
$\gamma'_1,\ldots,\gamma'_k\in\cP_\lambda^{\le \alpha}$ and $B\in 
H_2(M;\Z)$, we denote by   
$$ 
\widehat\cM^B(S_\og,S_\ug,\gamma'_1,\ldots,\gamma'_k;H_\infty,J_\infty) 
$$ 
the {\bf space of punctured Floer trajectories}, consisting of tuples 
$(u,L_1,\dots,L_k)$ such that  
$$ 
u=(f,a):\R\times S^1\setminus \{z_1,\ldots,z_k\}\to M\times \R 
$$ 
satisfies~(\ref{eq:F}--\ref{eq:Fasy}) for  
some $\gamma_-\in S_\og$, $\gamma_+\in S_\ug$, as well 
as~(\ref{eq:Fpunc}--\ref{eq:Lpunc}) and the relation  
\begin{equation} \label{eq:B}  
[(\sigma_\ug \cup \sigma_{\gamma_1} \cup \ldots \cup 
\sigma_{\gamma_k}) \# f ] = [\sigma_\og \# B]. 
\end{equation}  
The additive group $\R$  
acts freely by shifts on the domain and we define the {\bf moduli space of  
punctured Floer trajectories} by  
\begin{eqnarray*}  
\lefteqn{\cM^B(S_\og,S_\ug,\gamma'_1,\ldots,\gamma'_k;H_\infty,J_\infty)} 
\\  
& := & 
\widehat\cM^B(S_\og,S_\ug,\gamma'_1,\ldots,\gamma'_k;H_\infty,J_\infty)/\R. 
\end{eqnarray*} 
For a choice of regular $J_\infty$ as above, this is a smooth manifold of  
dimension (see also~\cite[\S3.3]{S}) 
$$ 
\mu(\og)-\mu(\ug) +2\langle c_1(\xi),B\rangle - \sum_{i=1}^k \bar\mu(\gamma'_i).
$$ 
Moreover, there are smooth evaluation maps 
\begin{equation} \label{eq:evB}
\oev,\uev:\cM^B(S_\og,S_\ug,\gamma'_1,\ldots,\gamma'_k;H_\infty,J_\infty)\to
S_\og,S_\ug.
\end{equation}
 
Given $A\in H_2(W;\Z)$ we define the {\bf moduli space of capped 
  punctured Floer trajectories}   
$$ 
\cM^A_c(S_\og,S_\ug;H_\infty,J) 
$$ 
as the set of equivalence classes of pairs $(u,F)$, where  
$u$ is an element of the moduli space $\cM^B(S_\og,S_\ug,\gamma'_1,\ldots,\gamma'_k;H_\infty,J_\infty)$,  
$\gamma'_1,\ldots,\gamma'_k\in\cP_\lambda^{\le \alpha}$, $B\in 
H_2(M;\Z)$, $F$ is a collection of $J$-holomorphic planes in $\widehat 
W$ of total homology class $A-B\in H_2(W;\Z)$, and with top asymptotes 
$\gamma'_1,\ldots,\gamma'_k$. Two sets of asymptotic markers at the 
punctures asymptotic to $\gamma'_1,\dots,\gamma'_k$ are equivalent if 
they satisfy~\eqref{eq:Arg}. This moduli space is a smooth
manifold of dimension  
$$ 
\mu(\og)-\mu(\ug) +2\langle c_1(TW),A\rangle. 
$$ 
Again, there are smooth evaluation maps 
\begin{equation} \label{eq:evAc} 
\oev,\uev:\cM^A_c(S_\og,S_\ug;H_\infty,J)\to
S_\og,S_\ug.
\end{equation}

The nonconstant elements of $\cP(H_\infty)$ and $\cP(H)$ are the
same. We denote by $\cF_{\mathrm{reg}}(H_\infty,J_\infty)$ the set whose
elements are collections of perfect Morse functions
$\{f_\gamma:S_\gamma\to\R,\, \gamma\in\cP(H)\}$ which satisfy the
transversality properties in~\eqref{eq:transvf} with respect to the
evaluation maps in~\eqref{eq:evB} (or, equivalently,
in~\eqref{eq:evAc}). The proof of Lemma~3.6 in~\cite{BOauto} carries
over verbatim and shows that $\cF_{\mathrm{reg}}(H_\infty,J_\infty)$
is of the second Baire category in the space of collections of perfect
Morse functions.  

Given $\{f_\gamma\}\in \cF_{\mathrm{reg}}(H_\infty,J_\infty)$,
$\og,\ug\in\cP(H)$, $p\in\mathrm{Crit}(f_\og)$,
$q\in\mathrm{Crit}(f_\ug)$, $A\in H_2(W;\Z)$, and $m\ge 0$, we
define the {\bf moduli space $\cM^A_{c;m}(p,q;H_\infty,\{f_\gamma\}, J)$ 
of capped punctured Morse-Bott broken trajectories with $m$ sublevels}
as the union for $\tgamma_1,\dots,\tgamma_{m-1}\in\cP(H)$ and
$A_1+\dots+A_m=A$ of the fibered products 
\begin{eqnarray*} 
&& 
\hspace{-.7cm}W^u(p)  
\times_{\oev} 
(\cM^{A_1}_c(S_{\og}\,,S_{\tgamma_1})\!\times\!\R^+) 
{_{\varphi_{f_{\tgamma_1}}\!\circ\uev}}\!\times   
_{\oev} 
(\cM^{A_2}_c(S_{\tgamma_1},S_{\tgamma_2})\!\times\!\R^+) \\ 
&&  
{_{\varphi_{f_{\tgamma_2}}\!\circ\uev}\times_{\oev}} \ldots\, 
{_{\varphi_{f_{\tgamma_{m-1}}}\!\!\circ\uev}}\!\!\times 
_{\oev}   
\cM^{A_m}_c(S_{\tgamma_{m-1}},\!S_{\ug})  
{_{\uev}\times} W^s(q). 
\end{eqnarray*} 
In the above formula one has to read
$$
\cM^{A_1}_c(S_{\og}\,,S_{\tgamma_1})=\cM^{A_1}_c(S_{\og}\,,S_{\tgamma_1};H_\infty,J)
\ \mbox{etc.}
$$
The {\bf moduli space of capped punctured Morse-Bott broken 
trajectories} is
$$ 
\cM^A_c(p,q;H_\infty,\{f_\gamma\},J):=\bigcup_{m\ge 0} 
\cM_m^A(p,q;H_\infty,\{f_\gamma\},J). 
$$ 
This is a smooth manifold of dimension 
\begin{eqnarray*} 
\lefteqn{\dim \, \cM^A_c(p,q;H_\infty,\{f_\gamma\},J)} \\ 
& = & \ \mu(\og) + \ind(p) - \mu(\ug) - \ind(q) + 2\langle 
c_1(TW),A\rangle - 1. 
\end{eqnarray*} 

In the statement of the next result we fix a sequence
$\tau^\nu\to\infty$, $\nu\to\infty$, and assume that  
\begin{equation}  \label{eq:Freg}
\{f_\gamma\}\in \bigcap_\nu
\cF_{\mathrm{reg}}(H^{\tau_\nu},J^{\tau_\nu}) \cap
\cF_{\mathrm{reg}}(H_\infty,J_\infty).
\end{equation}

\begin{proposition} \label{prop:stretch} 
There exists $\nu_0$ such that the following holds for any $\nu\ge
\nu_0$. For any 
$\og,\ug\in\cP(H)$, $p\in\mathrm{Crit}(f_\og)$, 
$q\in\mathrm{Crit}(f_\ug)$, and $A\in H_2(W;\Z)$ such that 
\begin{equation} \label{eq:dimension0}
\mu(\og) + \ind(p) - \mu(\ug) - \ind(q) + 2\langle 
c_1(TW),A\rangle - 1=0,
\end{equation}
the (discrete and finite) moduli spaces 
$$
\cM^A(p,q;H^{\tau_\nu},\{f_\gamma\},J^{\tau_\nu}) \quad and \quad 
\cM^A_c(p,q;H_\infty,\{f_\gamma\},J)
$$
are in natural bijective correspondence which preserves the signs of
their elements. 
\end{proposition}

Our proof of Proposition~\ref{prop:stretch} involves an extension of
the Correspondence Theorem~3.7 in~\cite{BOauto} which states that, in
dimension $0$, 
the moduli spaces of Morse-Bott broken trajectories are in
sign-preserving bijective correspondence with the moduli spaces of
Floer trajectories for a suitable time-dependent perturbation of the
Hamiltonian. This time-dependent perturbation is 
constructed as follows~\cite[Proposition~2.2]{CFHW}. For each
$\gamma\in\cP(H)$ we denote
by $\ell_\gamma\in\Z^+$ the maximal positive integer 
such that $\gamma(\theta+ 1 / \ell_\gamma)=\gamma(\theta)$,
$\theta\in S^1$. For each $S_\gamma$, $\gamma\in\cP(H)$ we choose a chart
$S^1\times \R^{2n-1}\ni(\widetilde \theta,p)$ in a
neighbourhood of the geometric image of $\gamma$ such that
$\widetilde \theta\circ\gamma(\theta)=\ell_\gamma\theta$ and $p\circ
\gamma(\theta)=0$. We also choose a smooth cut-off function
$\rho_\gamma=\rho_{S_\gamma}:S^1\times \R^{2n-1}\to \R$. For
$\delta>0$ small enough the time-dependent Hamiltonian 
\begin{equation} \label{eq:Hdelta}
   H_\delta(\theta,\widetilde \theta,p):= H +
\delta\sum_{S_\gamma}\rho_\gamma(\widetilde \theta,p)
   f_\gamma(\widetilde \theta-\ell_\gamma\theta)
\end{equation}
has only nondegenerate orbits. More precisely, out of each circle
$S_\gamma$ of elements of $\cP(H)$ survive precisely the orbits
$\gamma_m$ and $\gamma_M$ starting at the minimum,
respectively maximum of $f_\gamma$, with Conley-Zehnder
index $\mu(\gamma_p)=\mu(\gamma)+\mathrm{ind}(p)$.

\proof[Proof of Proposition~\ref{prop:stretch}] The Correspondence
Theorem~3.7 in~\cite{BOauto} implies the 
following. For each $\nu$, there exists $\delta_\nu>0$
such that, for any $\delta\in]0,\delta_\nu]$, the moduli space
$\cM^A(p,q;H^{\tau_\nu},\{f_\gamma\},J^{\tau_\nu})$ is in 
sign-preserving bijective correspondence with the moduli space of
Floer trajectories
$\cM^A(\og_p,\ug_q;H^{\tau_\nu}_\delta,J^{\tau_\nu})$ for the
time-dependent perturbation described by~\eqref{eq:Hdelta}. Here the standing
assumptions are~\eqref{eq:Freg} and~\eqref{eq:dimension0}, and we can
assume without loss of generality that $\delta_\nu\to 0$,
$\nu\to\infty$. 

Let $H^\nu:=H^{\tau_\nu}_{\delta_\nu}$,
$J^\nu:=J^{\tau_\nu}$. We prove Proposition~\ref{prop:stretch} in three steps. 

\smallskip 
\noindent {\bf Step~1.} {\it After passing to a subsequence, 
any sequence $v_\nu\in\cM^A(\og_p,\ug_q;H^\nu,J^\nu)$ converges
to an element of the moduli space $\cM^A_c(p,q;H_\infty,\{f_\gamma\},J)$.} 

There are two types of degenerations for the Floer
trajectories $v_\nu$ as $\nu\to\infty$. The first type, due to the fact
that the perturbation $\delta_\nu\to 0$, is that $v_\nu$ can spend an
amount of time $T_\nu\simeq T/\delta_\nu\to\infty$ in the neighbourhood of a
nontrivial periodic orbit of $H_\infty$. The second type, due to
the fact that we stretch the neck, is that 
$\|dv_\nu\|_{L^\infty}$ can become unbounded on suitable sequences
$(s_\nu,\theta_\nu)\in\R\times S^1$. 

The first type of degeneration is
dealt with in~\cite[Proposition~4.7]{BOauto}. The second type of
degeneration is the object of the compactness theorems for
SFT~\cite[Section~10.2]{BEHWZ}, applied in the following setting. 
By the mapping cylinder construction, the elements of the moduli space
$\cM^A(\og_p,\ug_q;H^\nu,J^\nu)$ can be interpreted as 
holomorphic sections of a symplectic fibration over $\R\times S^1$ 
with fiber $\widehat W^{\tau_\nu}$ for some almost complex 
structure $\widetilde J^\nu$. This is an almost complex manifold 
with symmetric cylindrical ends in the sense of~\cite[Sections~2.2 
and~3.2]{BEHWZ}. Although the case of a non-compact fiber is not 
explicitly considered in~\cite{BEHWZ}, this causes no problem here
because $\widehat W^{\tau_\nu}$ is convex at infinity. 

The important observation now is
that these two types of degenerations happen for {\it regions of $\R\times
S^1$ which have disjoint images via $v_\nu$}. More precisely, for the
second type of degeneration, the image of a 
neighbourhood of $(s_\nu,\theta_\nu)$ is necessarily contained in the
region $W\cup (M\times [-\tau_\nu,-C])\subset \widehat W^{\tau_\nu}$
for $\nu$ large enough (otherwise it would give rise to a nonconstant
$J_\infty$-holomorphic sphere in $M\times \R$). On the other hand, for
the first type of degeneration, the image of a cylinder of size
$T_\nu$ staying close to a nonconstant periodic orbit of
$H_\infty$ is contained in $M\times \R_+$. As a consequence, the
arguments used to deal with the first type of degeneration do not
interfere with those used to deal with the second type of
degeneration. We obtain, after passing to a subsequence, that 
$v_\nu$ converges to a tuple 
$$
\big(c_m,[\tu_m,F_m],c_{m-1},[\tu_{m-1},F_{m-1}],\ldots,[\tu_1,F_1],c_0\big),
\quad m\ge 0
$$
with the following properties (see also~\cite[Definition~4.1]{BOauto}):

\begin{itemize} 
 \item $\tu_i=(\tf_i,\ta_i): 
\Sigma_i\setminus \{z^i_1,\ldots,z^i_{k_i}\}\to M\times \R$,  
$i=1,\ldots,m$, $k_i\ge 0$ with  
 $$ 
 \Sigma_i=\R\times S^1 \quad \mbox{or} \quad  
 \Sigma_i=\R\times S^1\amalg\R\times S^1=:\oSigma_i\amalg\uSigma_i, 
 $$ 
 satisfying $\p_s\tu_i+J_\infty(\p_\theta
\tu_i-X_{H_\infty}(\tu_i))=0$ and the asymptotic conditions  
 \begin{itemize} 
  \item if $\Sigma_i=\R\times S^1$, then $\lim_{s\to-\infty}
\tu_i(s,\cdot)\in S_{\gamma_i}$,  
  $\lim_{s\to\infty} \tu_i(s,\cdot)\in S_{\gamma_{i-1}}$,
 \item if $\Sigma_i=\oSigma_i\amalg\uSigma_i$, then\\  
 $\lim_{s\to-\infty} \tu_i|_{\oSigma_i}(s,\cdot)\in S_{\gamma_i}$,  
  $\lim_{s\to\infty} \tu_i|_{\uSigma_i}(s,\cdot)\in S_{\gamma_{i-1}}$,\\ 
  $\lim_{s\to\infty} \tf_i|_{\oSigma_i}(s,\theta)=\og'_i(-\oT_i\theta)$,  
  $\lim_{s\to-\infty} 
\tf_i|_{\uSigma_i}(s,\theta)=\ug'_i(\uT_i\theta)$,\\ 
  $\lim_{s\to\infty} \ta_i|_{\oSigma_i}(s,\cdot)=-\infty$,  
  $\lim_{s\to-\infty}\ta_i|_{\uSigma_i}(s,\cdot)=-\infty$, 
 \end{itemize} 
 as well as~\eqref{eq:Fpunc} at $z^i_1,\ldots,z^i_{k_i}$ with asymptotes  
 ${{\tgamma}'{}}^i_1,\ldots,{{\tgamma}'{}}^i_{k_i}$.  
  
 Here $\gamma_i\in\cP(H_\infty)$, $i=0,\ldots,m$ with $\gamma_m=\og$,
$\gamma_0=\ug$, the orbits $\og'_i,\ug'_i\in\cP_\lambda^{\le \alpha}$,  
 $i=1,\ldots,m$ have periods $\oT_i,\uT_i$, and we have
 ${{\tgamma}'{}}^i_1,\ldots,{{\tgamma}'{}}^i_{k_i}\in\cP_\lambda^{\le\alpha}$. 
  
\item $c_0:[-1,+\infty[\to S_{\gamma_0}$, $c_i:[-T_i/2,T_i/2]\to
S_{\gamma_i}$, $i=1,\dots,m-1$ for some $T_i>0$, and 
$c_m:]-\infty,1]\to S_{\gamma_m}$, satisfy $\dot c_i=\nabla
f_{\gamma_i}\circ c_i$, $i=0,\dots,m$. 

\item $\oev(\tu_i)=\uev(c_i)$, $\uev(\tu_i)=\oev(c_{i-1})$,
$i=1,\dots,m$, and moreover we have $c_0(+\infty)=q$, $c_m(-\infty)=p$.  

 \item $F_i$, $i=1,\ldots,m$ is a  
genus zero $J$-holomorphic building in $\widehat W$  
of height $0|1|k_+$, $k_+\ge 0$ , whose top  
asymptotes are ${{\tgamma}'{}}^i_1,\ldots,{{\tgamma}'{}}^i_{k_i}$ and, 
possibly, $\og'_i,\ug'_i$   
if $\Sigma_i=\oSigma_i\amalg\uSigma_i$, and whose underlying 
nodal Riemann surface has exactly $k_i$, respectively $k_i+1$  
connected components.  
 
\item if $\Sigma_i=\oSigma_i\amalg\uSigma_i$ the following additional  
matching condition has to be satisfied. There exists a finite sequence  
$(F^i_1,\ldots,F^i_{n_i})$, $n_i\ge 1$ of connected components of the building  
$F_i$, on which we have marked punctures $\oz^i_j,\uz^i_j$, $j=1,\ldots,n_i$ such  
that  
\begin{itemize} 
\item $F^i_1$ is asymptotic to $\og'_i$ at $\oz^i_1$, $F^i_{n_i}$ is 
  asymptotic to $\ug'_i$ at $\uz^i_{n_i}$,  
\item the component $F^i_j$, $j=1,\ldots,n_i-1$  can be glued to $F^i_{j+1}$ 
at $\uz^i_j$, $\oz^i_{j+1}$,  
\item upon gluing in this way all the $F^i_j$, $j=1,\ldots,n_i$ 
and capping all the remaining punctures except $\oz^i_1$, $\uz^i_{n_i}$, the  
asymptotic markers at $\oz^i_1$, $\uz^i_{n_i}$ are opposite to each other with  
respect to the conformal structure on the cylinder.  
\end{itemize} 
Note that this last condition is present because the almost complex structure  
in Floer's equation~\eqref{eq:F} depends on $\theta$ (also, it ensures that  
pregluing produces an approximate solution).     
\end{itemize} 

We must show that any such tuple belongs to 
$\cM^A_c(p,q;H_\infty,\{f_\gamma\},J)$.  
For that purpose we notice that it is enough to show that the domain
of each $\tu_i$   
is connected. If this is the case, the dimension
condition~\eqref{eq:dimension0} implies  
that all the holomorphic buildings $F_i$ have height $k_+=0$, and we
recover the definition of an element of $\cM^A_c(p,q;H_\infty,\{f_\gamma\},J)$.  

Let us assume by contradiction that the domain of some  
$\tu_i$ has two components. Let $\tu_i=(f_i,a_i)$ and let $t_{i-1}$ be the level  
on which the asymptote $\gamma_{i-1}$ is located. Then the period $T_{i-1}$ of  
$\gamma_{i-1}$ is given by $T_{i-1}=e^{-t_{i-1}}h'(t_{i-1})$. 
Let $s_0\in\R$ be such that $\tu_i|_{\uSigma_i}([s_0,\infty[\times
S^1)$ is contained   
in a neighbourhood of $\gamma_{i-1}$ where $J_\infty=\widehat
J_\infty$ and therefore preserves $\xi$.   
We claim that there exists a point 
$(s,\theta)\in[s_0,\infty[\times S^1\subset \uSigma_i$ such that 
$a_i(s,\theta)>t_{i-1}$. Let  
$\bar a_i(s):=\int_{S^1}a_i(s,\theta) d\theta$. For $s\ge s_0$,    
Floer's equation $\p_s \tu_i + J_\infty\p_\theta\tu_i -
e^{-a_i}h'(a_i)\frac \p {\p t} =0$   
can be rewritten as 
\begin{eqnarray} 
 \p_s a_i - \lambda(\p_\theta f_i) - e^{-a_i}h'(a_i) & = & 0, \label{eq:Fa} \\ 
 \lambda(\p_s f_i) + \p_\theta a_i & = & 0, \nonumber \\ 
 \pi_\xi\circ df_i\circ j - J_\infty\pi_\xi\circ df_i & = & 0, \label{eq:Fxi} 
\end{eqnarray} 
where $\pi_\xi:TM\to\xi$ is the projection along $R_\lambda$.  
By~\eqref{eq:Fa} and using the fact that the function $t\mapsto e^{-t}h'(t)$  
is increasing we get that  
\begin{eqnarray*} 
\p_s \bar a_i (s)& = & \int_{S^1}\lambda(\p_\theta f_i(s,\theta)) d\theta +  
\int _{S^1} e^{-a_i(s,\theta)}h'(a_i(s,\theta)) d\theta \\ 
& \le & \int_{S^1}\lambda(\p_\theta f_i(s,\theta)) d\theta + T_{i-1}.  
\end{eqnarray*} 
On the other hand, for each $s_0\le s\le s'$ we have by Stokes' theorem  
\begin{equation} \label{eq:mon} 
\int_{S^1}\lambda(\p_\theta f_i(s',\theta)) d\theta -  
\int_{S^1}\lambda(\p_\theta f_i(s,\theta)) d\theta 
= \int_{[s,s']\times S^1} f^*d\lambda \ge 0, 
\end{equation} 
where the last inequality follows from~\eqref{eq:Fxi} and from the 
compatibility of $J_\infty$   
with $d\lambda$. Because $\lim_{s\to\infty}f_i(s,\cdot)=\gamma_{i-1}$, it 
follows that  
$\lim_{s\to\infty}\int _{S^1} \lambda(\p_\theta f_i(s,\theta))
d\theta=-T_{i-1}$ and,   
by the monotonicity relation~\eqref{eq:mon}, we get  
$$ 
\p_s \bar a_i (s) \le 0, \quad s\ge s_0.  
$$ 
Since $\lim_{s\to\infty} a_i(s,\theta) = t_{i-1}$, we have $\bar
a_i(s) \ge t_{i-1}$ for $s \ge s_0$.   
But if $a_i(s,\theta)$ is constant equal to $t_{i-1}$, then $\int f^*d\lambda=0$,  
which means that $\tu_i$ is constant, a contradiction. Hence we must have  
$a_i(s,\theta) > t_{i-1}$ for some $(s,\theta) \in [s_0, \infty) 
\times S^1$ as claimed.    
 
On the other hand, $\lim_{s\to-\infty} a_i(s,\theta) =-\infty$, so
that $a_i$ has a local maximum   
on $\R\times S^1$, which is impossible by the maximum principle. This
final contradiction  
shows that each $\tu_i$ has a domain which is connected.

\smallskip 
\noindent {\bf Step~2.} {\it We show that, for $\nu>0$ large enough, the elements 
of $\cM^A(\og_p,\ug_q;H^\nu,J^\nu)$ can be obtained  by a gluing construction
from the elements of $\cM^A_c(p,q;H_\infty,\{f_\gamma\},J)$. }

 For $\tau>0$ large enough let $H_\infty^\tau$ be the trivial extension 
 to $\widehat W^\tau$ of the Hamiltonian $H_\infty|_{M\times [-\tau,\infty[}$.  
For $\delta>0$ small enough, let $H_{\infty, \delta}^\tau$ be the perturbation
of $H_\infty^\tau$ described by~\eqref{eq:Hdelta}.

Combining the gluing arguments from Proposition 4.22 in~\cite{BOauto} with
aguments for gluing holomorphic curves along nondegenerate Reeb orbits in SFT
(\cite{B} and~\cite[Proposition~5]{BM}), we produce a Floer trajectory in 
$\widehat W^\tau$ for $H^\tau_{\infty,\delta}$, for $\tau>0$ large enough and 
$\delta>0$ small enough. As in Step~1, the important observation is 
that the gluing constructions take place in disjoint regions of $\widehat W^\tau$. 
More precisely, the gluing construction of punctured Floer trajectories
with gradient trajectories along nonconstant periodic 
orbits of $H_\infty$ takes place in $M \times \R_+$, while the gluing construction of
punctured Floer trajectories which are holomorphic near the punctures with
holomorphic planes takes place in $M \times [-\tau, -C]$.
Therefore, the gluing estimates do not interfere and we obtain the
desired Floer trajectory for $H^\tau_{\infty, \delta}$.
Since the elements of $\cM^A_c(p,q;H_\infty,\{f_\gamma\},J)$ are
rigid, the same arguments as in the uniqueness part
of~\cite[Proposition~4.22]{BOauto} show that we obtain in 
this way all elements in $\cM^A(\og_p,\ug_q;H^\tau_{\infty,\delta},J^\tau)$.

On the other hand, for $\tau>0$ large enough the Floer 
trajectories for $H^\tau_{\infty, \delta}$ and for $H_\delta^\tau$ are in 
bijective correspondence. The claim in Step~2 follows.  

\smallskip 
\noindent {\bf Step~3.} Combining  Step~1 and Step~2, 
we obtain a sign-preserving bijective correspondence between the $0$-dimensional 
moduli spaces $\cM^A(\og_p,\ug_q;H^\nu,J^\nu)$ and
$\cM^A_c(p,q;H_\infty,\{f_\gamma\},J)$, 
when $\nu>0$ is large enough.
As explained in the beginning of the proof, the former moduli space is
in sign-preserving bijective correspondence 
with the moduli space
$\cM^A(p,q;H^{\tau_\nu},\{f_\gamma\},J^{\tau_\nu})$, and the
conclusion follows. 
\hfill{$\square$}

\begin{remark} \label{rmk:sdep-capped}  
{\rm  
Given an increasing homotopy $H_s$ between admissible Hamiltonians 
$H_-\le H_+$, the continuation morphism    
$$ 
\sigma_{H_+,H_-}:BC_*^a(H_-)\to BC_*^a(H_+) 
$$  
can also be described in 
terms of capped punctured Floer trajectories. More precisely, 
equation~\eqref{eq:stretch} applied to the homotopy $H_s$ gives rise 
to increasing homotopies $H_s^\tau$ from $H_-^\tau$ to $H_+^\tau$ and, 
in the limit $\tau\to\infty$, to an increasing homotopy $H_{s,\infty}$ 
from $H_{-,\infty}$ to $H_{+,\infty}$ on $M\times \R$. The latter 
homotopy is supported in 
the region $t>-C$ for some large enough constant  
$C>0$. Given a generic homotopy $J_s$ of almost complex structures on 
$\widehat W$, we denote by $J_s^\tau$, respectively $J_{s,\infty}$ the 
corresponding homotopies of almost complex structures on $\widehat 
W^\tau$, respectively $M\times \R$. The moduli space  
$$ 
\cM^A_c(p,q;H_{s,\infty},\{f_\gamma^\pm\},J_{s,\infty}) 
$$  
of {\bf capped punctured $s$-dependent Morse-Bott broken 
tra\-jec\-to\-ries} is defined similarly to the moduli space  
$\cM^A(p,q;H_s^\tau,\!\{f_\gamma^\pm\},\!J_s^\tau)$ of $s$-dependent Morse-Bott broken 
trajectories introduced in Section~\ref{sec:MBsymp}. The modification of the 
definition is straightforward, with Floer trajectories for $H_-^\tau$, 
$H_+^\tau$, $H_s^\tau$ being replaced by capped punctured Floer trajectories for 
$H_{-,\infty}$, $H_{+,\infty}$, and respectively 
$H_{s,\infty}$. The obvious $s$-dependent version of 
Proposition~\ref{prop:stretch} still holds true with the same proof, 
so that the moduli spaces 
$\cM^A_c(p,q;H_{s,\infty},\{f_\gamma^\pm\},J_{s,\infty})$ and 
$\cM^A(p,q;H_s^\tau,\{f_\gamma^\pm\},J_s^\tau)$ of dimension $0$
are in sign-preserving bijective correspondence.   
} 
\end{remark}


\section{A filtered isomorphism} \label{sec:filtered} 
 
In this section we construct the chain complex 
isomorphism in Proposition~\ref{prop:filtered}  
by counting ``mixed moduli spaces'' consisting of rigid 
punctured cylinders which are asymptotic to a Reeb orbit at one end 
and to a Hamiltonian orbit at the other end. We prove in 
Proposition~\ref{prop:commdiag} that this isomorphism is compatible 
with the continuation maps in Floer homology.  
 
For this construction, we omit the free homotopy classes of loops  
from the notation.    
Let $\rho:\R\to [0,1]$ be a smooth increasing function such that 
$\rho(s)=0$ if $s<<0$ and  
$\rho(s)=1$ if $s>>0$. We define  
$$ 
H^\rho_\infty:\R\times (M\times\R) \to \R,\qquad (s,p,t)\mapsto \rho(s)H_\infty(t). 
$$ 
 
Given $\og',\gamma'_1,\ldots,\gamma'_k\in\cP_\lambda^{\le \alpha}$,  
$\ug\in\cP(H)$, $B\in H_2(M;\Z)$ we define the  
{\bf space of punctured interpolating  trajectories} 
$$ 
\cM^B(S'_\og,S_\ug,\gamma'_1,\ldots,\gamma'_k;H^\rho_\infty,J_\infty) 
$$ 
to consist of tuples $(u,L_1,\dots,L_k)$ such that  
$u=(f,a):\R\times 
S^1\setminus\{z_1,\ldots,z_k\}\to M\times \R$ satisfies  
\begin{equation} \label{eq:Fsnew} 
\p_s u + J_\infty(\p_\theta u - X_{H^\rho_\infty}(u))=0, 
\end{equation} 
\begin{equation} \label{eq:Fsasy-new} 
\lim_{s\to-\infty} a(s,\theta)=+\infty, \quad  
\lim_{s\to-\infty} f(s,\cdot)\in S'_\og,   
\end{equation} 
\begin{equation} \label{eq:Fsasy+new} 
\lim_{s\to\infty}u(s,\cdot)=\ug(\cdot+\utheta), \mbox{ for some } \utheta\in S^1, 
\end{equation} 
and equations~(\ref{eq:Fpunc}--\ref{eq:Lpunc}),~\eqref{eq:B}. Note that we have  
natural evaluation maps 
$$ 
\oev:
\cM^B(S'_\og,S_\ug,\gamma'_1,\ldots,\gamma'_k;H^\rho_\infty,J_\infty)\to
S'_\og  
$$   
and 
$$ 
\uev:\cM^B(S'_\og,S_\ug,\gamma'_1,\ldots,\gamma'_k;H^\rho_\infty,J_\infty)\to S_\ug. 
$$ 
As in Proposition~3.5 in~\cite{BOauto} we can choose $J_\infty$ regular for all the  
elements of  
$\cM^B(S'_\og,S_\ug,\gamma'_1,\ldots,\gamma'_k;H^\rho_\infty,J_\infty)$ 
and so that it still satisfies   
all the previous regularity assumptions. Note 
that, due to the last term in the left hand side of~\eqref{eq:Fsnew},  
the additive group $\R$ \emph{does not} act on the space $\cM^B$, and therefore we 
call it from now on the {\bf moduli space of punctured interpolating 
  trajectories}. This is a smooth manifold of dimension~\cite[\S3.3]{S} 
\begin{eqnarray*} 
\lefteqn{\dim 
  \cM^B(S'_\og,S_\ug,\gamma'_1,\ldots,\gamma'_k;H^\rho_\infty,J_\infty)} \\   
& = & \  
\mu(\og') - \mu(\ug) +2\langle c_1(\xi),B\rangle  
- \sum_{i=1}^k \bar\mu(\gamma'_i)+1.  
\end{eqnarray*} 
 
\begin{figure}[hpt] 
\centering 
\input{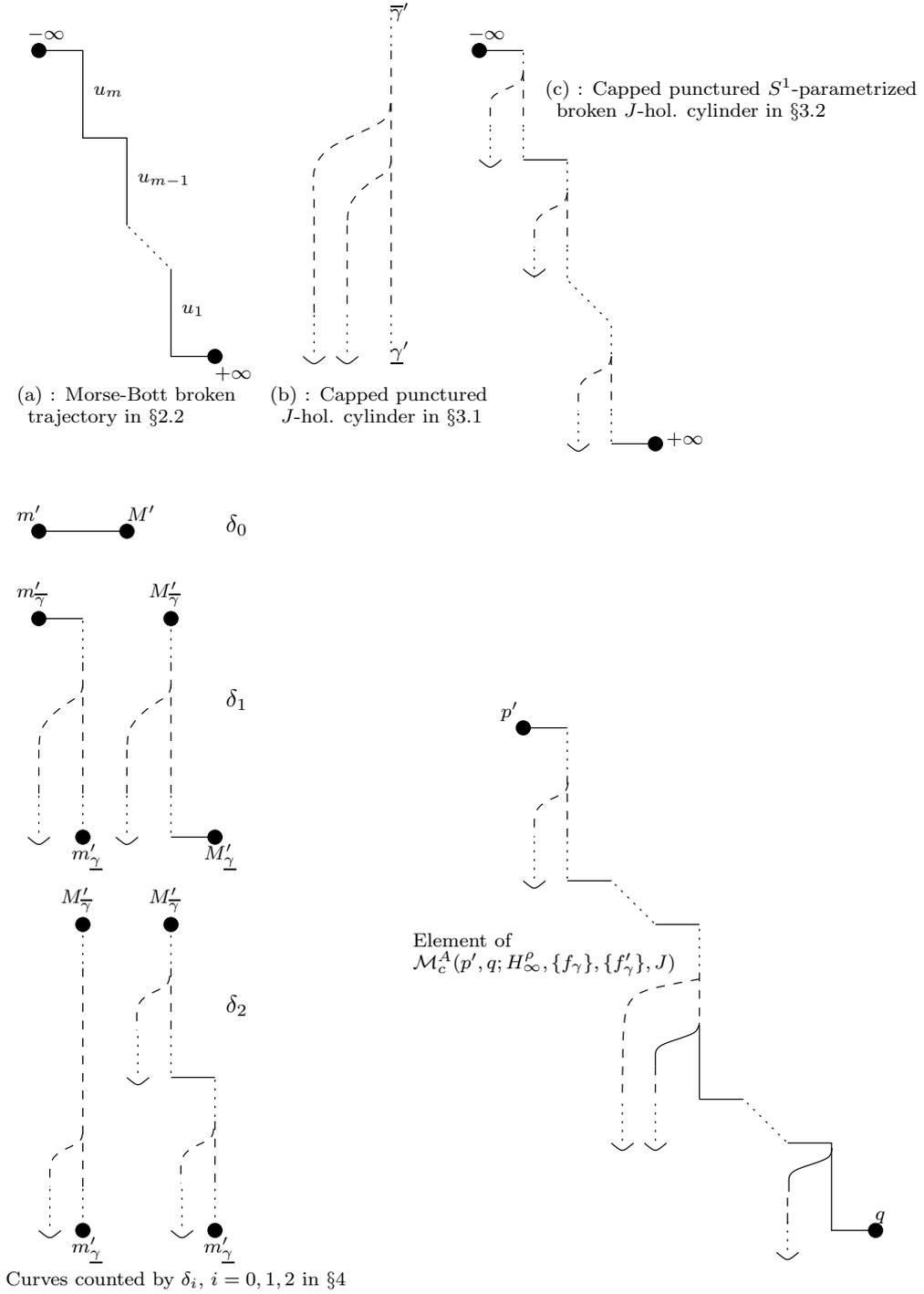} 
\caption{Moduli spaces \label{fig:Phi}} 
\end{figure}

Given $m,\ell\ge 0$, let $\og',{\gamma'}^1_1, \ldots, 
{\gamma'}^1_{k_1}, \ldots,  
{\gamma'}^{m+\ell+1}_1, \ldots, {\gamma'}^{m+\ell+1}_{k_{m+\ell+1}} 
\in\cP_\lambda^{\le\alpha}$, $\ug \in \cP(H)$, 
$p'\in\textrm{Crit}(f'_\og)$, $q\in 
\mathrm{Crit}(f_{\ug})$ and $B\in H_2(M;\Z)$.  
We denote by  
\begin{equation*}  
\cM_{m,\ell}^B(p',q,{\gamma'}^1_1, \ldots, {\gamma'}^1_{k_1}, \ldots, 
{\gamma'}^{m+\ell+1}_1,  
\ldots, 
{\gamma'}^{m+\ell+1}_{k_{m+\ell+1}};H^\rho_\infty,\{f_\gamma,f'_\gamma\},J_\infty)   
\end{equation*} 
the union over $\tgamma'_1,\ldots,\tgamma'_m\in\cP_\lambda$, 
$\tgamma_{m+1},\ldots,\tgamma_{m+\ell}\in \cP(H)$ and 
$B_1+\ldots +B_{m+\ell+1}=B$ of the fibered products   
\begin{eqnarray*} 
&& 
W^u(p')\times _{\oev} ((\cM^{B_1}(S'_\og,S'_{\tgamma_1},{\gamma'}^1_1, 
\ldots, {\gamma'}^1_{k_1};J_\infty)/\R)  \times \R^+)  
{_{\varphi_{f'_{\tgamma_1}}\circ\uev}}\times _{\oev} \\ 
&&  
\ldots ((\cM^{B_m}(S'_{\tgamma_{m-1}},S'_{\tgamma_m},{\gamma'}^m_1, \ldots, 
{\gamma'}^m_{k_m};J_\infty)/\R)  \times  \R^+) 
{_{\varphi_{f'_{\tgamma_m}}\circ\uev}}\times _{\oev} \\ 
&&  
(\cM^{B_{m+1}}(S'_{\tgamma_m},S_{\tgamma_{m+1}},{\gamma'}^{m+1}_1
\!\!\!\!\!,\ldots,    
{\gamma'}^{m+1}_{k_{m+1}};H_\infty^\rho,J_\infty) \times \R^+) 
{_{\varphi_{f_{\tgamma_{m+1}}}\circ\uev}}\times _{\oev} \\ 
&&  
(\cM^{B_{m+2}}(S_{\tgamma_{m+1}},S_{\tgamma_{m+2}},{\gamma'}^{m+2}_1 
\!\!\!\!\!\!\!,\ldots, 
{\gamma'}^{m+2}_{k_{m+2}};H_\infty,J_\infty)\!\times\!\R^+\!) 
{_{\varphi_{f_{\tgamma_{m+2}}}\circ\uev}}\times _{\oev}  \\ 
&&  
\ldots 
\cM^{B_{m+\ell+1}}(S_{\tgamma_{m+\ell}},S_\ug,{\gamma'}^{m+\ell+1}_1, 
\ldots,  
{\gamma'}^{m+\ell+1}_{k_{m+\ell+1}};H_\infty,J_\infty)  
{_{\uev}} \times W^s(q). 
\end{eqnarray*} 
By our transversality assumptions, this is a smooth manifold of dimension 
\begin{eqnarray*} 
\lefteqn{ 
\hspace{-1cm}\dim \, \cM_{m,\ell}^B(p',q,{\gamma'}^1_1,..., 
  {\gamma'}^1_{k_1},..., {\gamma'}^{m+\ell+1}_1 
\hspace{-.5cm},..., 
  {\gamma'}^{m+\ell+1}_{k_{m+\ell+1}};
H^\rho_\infty,\{f_\gamma,f'_\gamma\},J_\infty)}   
  \\  
& = & \ind(p') - 1 + (\dim \, 
\cM^{B_1}(S'_\og,S'_{\tgamma_1})/\R  + 1) - 1 \\ 
&&+ \ ... + 
(\dim \,\cM^{B_m}(S'_{\tgamma_{m-1}},S'_{\tgamma_m})/\R +1) - 1  \\  
&& + \ (\dim \,\cM^{B_{m+1}}(S'_{\tgamma_m},S_{\tgamma_{m+1}}) +1) - 1 
  \\  
&& + \ (\dim \,\cM^{B_{m+2}}(S_{\tgamma_{m+1}},S_{\tgamma_{m+2}}) +1) - 
  1  \\ 
&& + \ \ldots + \ \dim \,\cM^{B_{m+\ell+1}}(S_{\tgamma_{m+\ell}},S_\ug) - 
  1 + (1-\ind(q)) \\  
& = & \mu(\og') - \mu(\ug)  +\ind(p') - \ind(q) \\ 
&& + \ 2\langle 
c_1(\xi),B_1+...+B_{m+\ell+1}\rangle  
- \sum_{i=1}^{m+\ell+1} \sum_{j=1}^{k_i} \bar\mu({\gamma'}^i_j)  \\ 
& = & \mu(\og'_p) - \mu(\ug_q) - \sum_{i=1}^{m+\ell+1} 
  \sum_{j=1}^{k_i} \bar\mu({\gamma'}^i_j)   
+ 2 \langle c_1(\xi),B\rangle. 
\end{eqnarray*}  
We denote  
\begin{eqnarray*} 
\lefteqn{\cM^B(p',q,\gamma'_1, \ldots, 
  \gamma'_k;H^\rho_\infty,\{f_\gamma,f'_\gamma\},J_\infty):=} \\  
&& \hspace{-.2cm}\bigcup_{m,\ell\, \ge 0} 
\hspace{-.2cm}\cM_{m,\ell}^B(p',q,{\gamma'}^1_1,..., 
  {\gamma'}^1_{k_1},..., {\gamma'}^{m+\ell+1}_1  
\hspace{-.5cm},..., 
  {\gamma'}^{m+\ell+1}_{k_{m+\ell+1}};H^\rho_\infty,\{f_\gamma,f'_\gamma\},J_\infty) 
\end{eqnarray*}  
with $\{  {\gamma'}^1_1, \ldots, 
  {\gamma'}^1_{k_1}, \ldots, {\gamma'}^{m+\ell+1}_1,  
\ldots, 
  {\gamma'}^{m+\ell+1}_{k_{m+\ell+1}}\} = \{ \gamma'_1, \ldots, \gamma'_k \}$, 
and we call this {\bf the moduli space of punctured interpolating 
  Morse-Bott broken  
trajectories}, whereas the spaces $\cM_{m,\ell}^B$ are called {\bf 
the moduli spaces of punctured interpolating Morse-Bott broken 
  trajectories with $m$ holomorphic sublevels and $\ell$ Floer sublevels}.   
We refer to Figure~\ref{fig:Phi} for a representation of the elements 
  of these moduli spaces (see also their capped version below).  
The description of    
the topology of the compactified moduli spaces $\overline
  \cM^B(p',q,\gamma'_1, \ldots,  
  \gamma'_k;H^\rho_\infty,\{f_\gamma,f'_\gamma\},J_\infty)$ is 
  entirely similar to that of the compactified moduli spaces  
$\overline \cM^A(p,q;H,\{f_\gamma\},J)$ in Section~\ref{sec:MBsymp} and  
$\overline \cM^B(p',q',{\gamma'}_1, \ldots, 
{\gamma'}_k;\{f'_\gamma\},J_\infty)$ in Section~\ref{sec:nonequiv}.

The system of coherent orientations in~\cite[Section~4.4]{BOauto} orients  
zero-dimensional moduli spaces as above (even if the intermediate  
orbits $\tgamma'_i$ are bad) and hence there is  a sign 
$\oeps(u)$ associated to each element $u$ in a zero-dimensional 
moduli space.    
 
We define 
\begin{equation*}  
\Phi : BC_*^{\le\alpha}(\lambda) \to BC_*(H) 
\end{equation*}  
by  
\begin{equation} \label{eq:bigPhi}  
\Phi(\og'_p) := \sum_{u} \frac 1 {{\scriptstyle \prod_{i=1}^k \kappa_{\gamma'_i}}}  
\ \oeps(u) e(\gamma'_1) \ldots e(\gamma'_k) e^B \ug_q, 
\end{equation}  
where we sum over $u \in \cM^B(p',q,\gamma'_1, \ldots, 
  \gamma'_k;H^\rho_\infty,\{f_\gamma,f'_\gamma\},J_\infty)$.  
 
Let us derive now an alternative expression for 
the map $\Phi$ in terms of capped moduli spaces.   
Given $A\in H_2(W;\Z)$, $\og'\in \cP_\lambda^{\le\alpha}$, 
$\ug\in\cP(H)$, $p'\in\textrm{Crit}(f'_\og)$, $q\in 
\textrm{Crit}(f_\ug)$ we define the {\bf moduli space of capped 
  punctured interpolating Morse-Bott broken trajectories}  
\begin{equation} \label{eq:cpiMB} 
\cM^A_c(p',q;H_\infty^\rho,\{f_\gamma,f'_\gamma\},J) 
\end{equation} 
as the set of equivalence classes of pairs $F=(F',F'')$, where 
$F'$ is an element of the moduli space $\cM^B(p',q,\gamma'_1,\ldots,\gamma'_k; 
H_\infty^\rho,\{f_\gamma,f'_\gamma\},J_\infty)$,   
$\gamma'_1,\ldots,\gamma'_k\in\cP_\lambda$, $B\in H_2(M;\Z)$ 
and $F''$ a 
collection of $J$-holomorphic planes in $\widehat W$, of total homology 
class $A-B\in H_2(W;\Z)$, and whose top asymptotes are 
$\gamma'_1,\ldots,\gamma'_k$. Two sets of asymptotic markers at the 
punctures asymptotic to $\gamma'_1,\dots,\gamma'_k$ are equivalent if 
they satisfy~\eqref{eq:Arg}. The dimension of this moduli space is  
\begin{equation} \label{eq:dimcappedpuncintMBbr}  
\mu(\og'_p) - \mu(\ug_q) + 2\langle c_1(TW),A\rangle. 
\end{equation}  
In the case of dimension zero one can associate a sign $\oeps(F)$ to 
each element  
$F\in\cM^A_c(p',q;H_\infty^\rho,\{f_\gamma,f'_\gamma\},J)$.  
The map $\Phi$ can then be rewritten as 
\begin{equation} \label{eq:Phialt}  
\Phi(\og'_p) := \sum_{F \in 
  \cM^A_c(p',q;H^\rho_\infty,\{f_\gamma,f'_\gamma\},J)}   
\oeps(F) e^A \ug_q. 
\end{equation}  
 
We shall also need in the proof of Proposition~\ref{prop:filtered} the 
following variant of the moduli space of capped punctured 
interpolating Morse-Bott broken trajectories. Given  
$A\in H_2(W;\Z)$, $\og',\gamma'_1\in \cP_\lambda^{\le\alpha}$, 
$\ug\in\cP(H)$, $p'\in\mathrm{Crit}(f'_\og)$, $q\in 
\textrm{Crit}(f_\ug)$ we denote by  
$$ 
\cM^A_{c,1}(p',q,\gamma'_1;H_\infty^\rho,\{f_\gamma,f'_\gamma\},J) 
$$ 
the set of equivalence classes of pairs $(F',F'')$, where 
$F'$ is an element of the moduli space 
$\cM^B(p',q,\gamma'_1,\gamma'_2,\ldots,\gamma'_k;H_\infty^\rho, 
\{f_\gamma,f'_\gamma\},J)$ with         
$\gamma'_2,\ldots,\gamma'_k\in \cP_\lambda^{\le\alpha}$, $B\in 
H_2(M;\Z)$, $F''$ is a collection of $J$-holomorphic planes in $\widehat 
W$ of total homology class $A-B\in H_2(W;\Z)$, and with top asymptotes 
$\gamma'_2,\ldots,\gamma'_k$. Two sets of asymptotic markers at the 
punctures asymptotic to $\gamma'_2,\dots,\gamma'_k$ are equivalent if 
they satisfy~\eqref{eq:Arg}. One should think of this 
as being the moduli space of punctured interpolating Morse-Bott broken 
trajectories which are capped at all but one of the punctures.

\proof[Proof of Proposition~\ref{prop:filtered}] 
 The map $\Phi$ preserves filtrations since all the moduli spaces in the 
 definition of     
 $\cM^A_c(p',q;H^\rho_\infty,\{f_\gamma,f'_\gamma\},J)$ are 
 regular and hence $\mu(\og)-\mu(e^A\ug)\ge 0$. As a matter of fact, 
 since~\eqref{eq:Phialt} involves only moduli spaces 
 such that $\mu(\og'_p)-\mu(\ug_q) +2\langle c_1(TW),A \rangle=0$, and 
 since $\ind(p)-\ind(q)\in\{1,0,-1\}$, we have 
 $\mu(\og)-\mu(e^A \ug)\in\{0,1\}$. Moreover, the 
 relevant moduli spaces are empty if $p'$ is a minimum and $q$ is a 
 maximum.  
 
 We now prove that $\Phi$ is a morphism of differential complexes, 
 i.e. that it satisfies the relation  
\begin{equation} \label{eq:Phi-diff} 
 \Phi \circ \delta = d \circ \Phi. 
\end{equation}  
Let us consider the $1$-dimensional components of the moduli space 
of capped punctured interpolating Morse-Bott broken trajectories
$\cM^A_c(\og'_p,\ug_q;H^\rho_\infty,\{f_\gamma,f'_\gamma\},J)$. We   
 claim that its boundary has the form  
\begin{eqnarray} \label{eq:boundary}  
\lefteqn{\bigcup _{ \stackrel{\gamma'\in \cP_\lambda^{\le\alpha} } 
{ r'\in\mathrm{Crit}(f'_\gamma) }} 
\hspace{-.3cm} \cM^{A_1}_c(p',r';\{f'_\gamma\}, J) \times 
 \cM^{A-A_1}_c(r',q;H^\rho_\infty,\{f_\gamma,f'_\gamma\},J)} \\ 
& \cup &  
\bigcup_{\stackrel{\gamma\in \cP(H)} 
{r\in\mathrm{Crit}(f_\gamma)}} 
 \cM^{A_1}_c(p',r;H^\rho_\infty,\{f_\gamma,f'_\gamma\},J) \times 
\cM^{A-A_1}_c(r,q;H_\infty,\{f_\gamma\},J) \nonumber \\ 
& \cup &  
\bigcup _{\gamma'_1,\gamma'_2\in \cP_\lambda^{\le\alpha}} 
\big[ \cM^{A_1}_{c,1}(p',q,\gamma'_1;H_\infty^\rho,\{f_\gamma,f'_\gamma\},J) 
 \nonumber \\ 
& & \hspace{2cm} 
 \times \ 
 \cM^{A_2}_c(\gamma'_1,\gamma'_2;J)/\R \times 
 \cM^{A-A_1-A_2}(\gamma'_2,\emptyset;J)\big], \nonumber 
\end{eqnarray}  
with $[\dots]$ denoting triples modulo the equivalence relation given 
by~\eqref{eq:Arg} for the 
asymptotic markers at the punctures with common asymptotics 
$\gamma'_1$,$\gamma'_2$.  
Indeed, the compactness theorem for SFT~\cite[Theorem~10.2]{BEHWZ}    
ensures that the boundary has the form above with 
all moduli spaces replaced by their compactifications. In our 
situation the boundary is $0$-dimensional and the above moduli spaces 
are already compact.   
 
Let us discuss the contribution of the three types of boundary components 
above. The count of the elements of the first type corresponds to 
$\Phi\circ \delta$, whereas the count of the 
elements of the second type corresponds to $d \circ \Phi$.  
 
We claim 
that the count of the elements of the third type is equal to 
$0$. Indeed, the count of elements of  
$$ 
\big[\cM^{A_2}_c(\gamma'_1,\gamma'_2;J)/\R \times 
 \cM^{A-A_1-A_2}(\gamma'_2,\emptyset;J)\big] 
$$ 
is equal to $e \circ \p (\gamma'_1) =0$ by~\eqref{eq:ed}. Here 
$[\dots]$ denotes pairs modulo the equivalence relation given 
by~\eqref{eq:Arg} for the 
asymptotic markers at the punctures with common asymptote 
$\gamma'_2$.  
This shows that $\Phi \circ \delta \pm d \circ 
\Phi =0$. Relation~\eqref{eq:Phi-diff} follows from the fact  
that  
orientations are coherent with the gluing operation, and we refer  
to~\cite[p.69]{FH1} for details (the situation is analogous to the 
one in Floer homology involving a homotopy of Hamiltonians, i.e. an 
equation which is not translation invariant).  
 
We now prove that $\Phi$ is an isomorphism. The first important 
observation is that $\Phi$ increases the contact action $\gamma\mapsto \int 
\gamma^*\lambda$.  
This is due to the fact that $\int u^* d\lambda \ge 0$ if $u$ 
satisfies~\eqref{eq:Fsnew}. Indeed, 
$i(X_{H_\infty^\rho})d\lambda =0$ so 
that $d\lambda(\p_s u,\p_\theta u)= d\lambda(\p_s u,J_\infty \p_s 
u)\ge 0$.

We arrange the generators of the complex 
$BC_*^{i^{-1}(a),\le\alpha}(\lambda)$ in increasing order 
according to their action. Then the matrix of $\Phi$ in this basis is lower 
triangular and we claim that the entries on the diagonal are all equal 
to $\pm1$, which implies that $\Phi$ is an isomorphism.  
 
Given $\gamma\in\cP_\lambda^{\le\alpha} \equiv \cP(H)$, 
$p'\in\textrm{Crit}(f'_\gamma)$, $p\in\textrm{Crit}(f_\gamma)$ with 
$\ind(p)=\ind(p')$, we must determine the elements
$u\in\cM^B(p',p,\gamma'_1,\ldots,\gamma'_k; 
H_\infty^\rho,\{f_\gamma,f'_\gamma\},J_\infty)$.   
Since  
$$ 
\int u^*d\lambda = \cA(\gamma) - \cA(\gamma) - \sum_{i=1}^k 
\cA(\gamma'_i) \ge 0 
$$  
we infer that $k=0$. Moreover $\int u^* d\lambda 
=0$, which implies that $u(s,\theta)$ is a vertical cylinder of the 
form  
$$ 
u(s,\theta)=(a(s,\theta),\gamma \circ  b(s,\theta)).  
$$ 
 
Equation~\eqref{eq:Fsnew} is then equivalent to the system 
\begin{equation} \label{eq:abnew}  
\left\{ \begin{array}{rcl} \p_s a - \p_\theta b - \rho(s) H_\infty'(a) e^{-a} 
    & = & 0, \\ 
  \p_\theta a + \p_s b & = & 0. 
\end{array} \right. ,  
\end{equation}  
and we claim that any solution satisfying 
$$ 
\lim_{s\to-\infty} a(s,\cdot)=+\infty, \quad \lim_{s\to +\infty} 
a(s,\cdot) = a_0,  
$$ 
\begin{equation*}  
\lim_{s\to\pm\infty} b(s,\theta)= -T\theta +\theta_\pm, \quad 
T=e^{-a_0}h'(a_0)  
\end{equation*}  
has the following property:  
\begin{enum}  
\item if $p'$, $p$ are both minima and 
$\theta_+$ corresponds to the orbit $\gamma_p$ (and $\theta_-$ is 
arbitrary), then  
$$ 
a(s,\theta)=a(s), \quad b(s,\theta)=-T\theta + \theta_+; 
$$ 
\item if $p'$, $p$ are 
both maxima and $\theta_-$ corresponds to the orbit $\gamma'_p$ (and 
$\theta_+$ is arbitrary), then  
$$ 
a(s,\theta)=a(s), \quad b(s,\theta)=-T\theta + \theta_-. 
$$ 
\end{enum}  
 
We treat only the first case since the other case is entirely 
analogous. We linearize~\eqref{eq:abnew} and obtain  
\begin{equation} \label{eq:eqlinnew} 
\p_s\zeta + J_0\p_\theta\zeta +S \zeta =0, \qquad \zeta:\R\times S^1 \to \R^2, 
\end{equation}  
with $S(s,\theta)=- \rho(s)\left(\begin{array}{cc} 
    \p_a(e^{-a}H_\infty'(a)) & 0 \\ 0 & 0 \end{array}\right)$. Our 
hypothesis~\eqref{eq:slow} on $h$ ensures that  
$$ 
\| S \| <1.  
$$ 
By~\cite[Proposition~4.2]{SZ} applied with $\tau=1$ we infer that 
every solution $\zeta$ of~\eqref{eq:eqlinnew} is independent of 
$\theta$. Hence solutions of~\eqref{eq:eqlinnew} satisfy  
$$ 
\p_s \zeta_1 = \rho(s)\p_a(e^{-a}H_\infty'(a)) \zeta_1, \quad \p_s \zeta_2=0. 
$$ 
In particular we see that, if $\lim_{s\to+\infty}\zeta_1(s,\cdot)=0$, 
we necessarily have $\zeta_1\equiv 0$ because 
$\rho(s)\p_a(e^{-a}H_\infty'(a))$ is strictly positive for $s$ large.  
 
Since equation~\eqref{eq:abnew} is independent of $\theta$, we 
know that, for any solution $(a,b)$ the vector $\zeta=(\p_\theta 
a,\p_\theta b)$ solves~\eqref{eq:eqlinnew}. In particular $\p_\theta 
b=-T$ (by the asymptotic conditions on $b$) and $\p_s b = -\p_\theta 
a=0$, so that $b(s,\theta)=-T \theta + \theta_+$ as claimed.  
 
We have already proved that $\p_\theta a=0$, hence $a(s,\theta)=a(s)$ 
solves the differential equation 
$a'+e^{-a_0}H_\infty'(a_0)-\rho(s)(e^{-a}H_\infty'(a))=0$. Since 
$e^{-a} H_\infty'(a)$ is increasing with $a$, we must have $a(s)=a_0$ 
when $\rho(s)=1$. Hence the solution $a(s)$ is unique.  
 
As a conclusion, we have showed that  
$\cM^A_c(p',p;H_\infty^\rho,\{f_\gamma,f'_\gamma\},J)$ contains a 
unique element, which is the ``trivial'' cylinder over $\gamma$ described 
above.  
\hfill{$\square$}

\begin{remark}[An inverse for $\Phi$] \label{rmk:inversePhi} 
{\rm We chose to construct the 
    morphism $\Phi$ by counting mixed moduli spaces of curves with 
    asymptotes $\og'_p$ at $-\infty$ and $\ug_q$ at $+\infty$. In 
    heuristic terms, we considered moduli spaces of curves ``going 
    down'' from Reeb orbits to Hamiltonian orbits, and for that 
    purpose we needed to consider the cut-off Hamiltonian 
    $H_\infty^\rho$.  
 
   We might have as well proceeded the other way around,   
   by constructing a filtered isomorphism  
  $$ 
 \Psi:BC_*^a(H) \stackrel \sim \to BC_*^{i^{-1}(a),\le\alpha}  
  $$ 
  obtained by counting mixed moduli spaces of curves with asymptotes 
  $\og_p$ at $-\infty$ and $\ug'_q$ at $+\infty$, i.e. ``going down'' 
  from Hamiltonian orbits to Reeb orbits. For that purpose 
  equation~\eqref{eq:Fsnew} would have had to be replaced by  
$$ 
\p_su +J_\infty(\p_\theta u - X_{H_\infty}(u)) =0.  
$$ 
The one difference with respect to~\eqref{eq:Fsnew}   
is that this equation does not depend on the variable 
$s$ and therefore $\R$ acts by translations on the spaces of 
solutions. Nevertheless, the same dimension formula holds 
for the corresponding moduli space after the quotient by this  
$\R$-action.

One can prove that $\Psi$ is an inverse for $\Phi$ up to 
(filtered) homotopy. 
} 
\end{remark}  
 
We introduce now an algebraic concept which    
is useful when working 
with spectral sequences. Let $(C_*,\delta),(D_*,d)$ be differential 
complexes endowed with increasing filtrations $F_\ell C_*$ and 
$F_\ell D_*$ with $\ell\in\Z$. Let $f,g:C_*\to D_*$ be filtration 
preserving chain maps and $K:C_*\to D_{*+1}$ be a chain homotopy 
between $f$ and $g$, so that $f-g=K\delta+dK$. We say that $K$ is a 
{\bf chain homotopy of order $k\ge 0$} if $K(F_\ell C_*)\subset 
F_{\ell+k}D_{*+1}$ for all $\ell \in \Z$.  
 
Let $(E^r_\delta,\bar \delta^r)$ and $(E^r_d,\bar d^r)$, $r\ge 0$ be 
the spectral sequences associated to the given filtrations on 
$(C_*,\delta)$ and respectively $(D_*,d)$. Then $f$ and $g$ 
induce chain maps $f_r,g_r:E^r_\delta \to E^r_d$ for all $r\ge 
0$. The next result motivates the concept of a chain homotopy of order $k$. 
 
\begin{proposition}[{\cite[Exercise~3.8, p.~87]{McC}}]  
 \label{prop:homotopy} 
If $K:C_*\to 
  D_{*+1}$ is a chain homotopy of order  
  $k$ between the filtered chain maps $f,g:C_*\to D_*$, then  
$f_k,g_k:E^k_\delta\to E^k_d$ are chain homotopic, and 
the maps $f_r,g_r:E^r_\delta\to E^r_d$ coincide for any $k<r\le 
\infty$.  
\end{proposition}  
 
The next Proposition shows that the filtered isomorphism $\Phi$ of 
Proposition~\ref{prop:filtered}  
is compatible with the continuation morphisms in Floer homology. Let 
$\tau>0$ be large  
enough and, with the notations of Section~\ref{sec:stretch}, let  
$H^\alpha=H^{\alpha,\tau}\le H^{\alpha'}=H^{\alpha',\tau}$  
be admissible Hamiltonians on $\widehat W^\tau$ with asymptotic slopes 
$\alpha,\alpha'\notin\mathrm{Spec}(M,\lambda)$. 
Let $H_s=H_s^\tau$, $s\in\R$ be an increasing 
homotopy from $H^\alpha$ to $H^{\alpha'}$, such that $\p_sH_s$ is small  
enough in $C^0$-norm.  Let $J_s=J_s^\tau$ be a generic 
regular homotopy of 
admissible almost complex structures on $\widehat W^\tau$ which is a small 
perturbation of the time independent almost complex structure 
postulated in our standing assumptions~$(A)$ and~$(B_a)$. Thus, 
Proposition~\ref{prop:stretch} and Remark~\ref{rmk:sdep-capped} apply
to the Hamiltonians $H^\alpha$, $H^{\alpha'}$, and respectively 
$H_s$.

\begin{proposition} \label{prop:commdiag} 
The diagram below commutes up to a chain homotopy of order $1$ 
\begin{equation}\label{eq:last_diag*} 
\xymatrix 
@C=25pt 
@R=15pt 
{BC_*^{i^{-1}(a),\le \alpha}(\lambda) \ar[r]^-{\Phi^\alpha}  
\ar[d]_-{\kappa_{\alpha',\alpha}} & 
BC_*^a(H^\alpha) \ar[d]^-{\sigma_{\alpha',\alpha}} \\ 
BC_*^{i^{-1}(a),\le \alpha'}(\lambda)  \ar[r]_-{\Phi^{\alpha'}} &  
BC_*^a(H^{\alpha'}) 
} 
\end{equation}  
Here  
$\Phi^\alpha$, $\Phi^{\alpha'}$ are the filtered isomorphisms of 
Proposition~\ref{prop:filtered} for the Hamiltonians $H^\alpha$, 
respectively $H^{\alpha'}$, the map $\sigma_{\alpha',\alpha}$ is the 
continuation morphism defined in Section~\ref{sec:MBsymp}, and 
$\kappa_{\alpha',\alpha}$ is the inclusion.  
\end{proposition}  
 
\proof 
All morphisms in the diagram~\eqref{eq:last_diag*} preserve the 
filtrations on the corresponding chain complexes. This was proved 
in Proposition~\ref{prop:filtered} for $\Phi^\alpha$ and 
$\Phi^{\alpha'}$, it was proved in Remark~\ref{rmk:sdep} for the 
continuation  
morphism $\sigma_{\alpha',\alpha}$, and it follows directly from the 
definition for the inclusion $\kappa_{\alpha',\alpha}$.  
 
We denote by $H_{s,\infty}:M\times \R\to\R$ the homotopy from 
$H^\alpha_\infty$ to $H^{\alpha'}_\infty$ defined via 
formula~\eqref{eq:stretch} applied to $H_s$ (see also 
Remark~\ref{rmk:sdep-capped}).  
As above, let $\rho:\R\to [0,1]$ be a smooth increasing function  
such that $\rho(s)=0$ if $s<< 0$ and $\rho(s)=1$ if $s>> 0$. Given 
$r\in \R$ we define the Hamiltonian family $H_{s,\infty}^{\rho,r}$, 
$s\in \R$ by 
$$ 
H_{s,\infty}^{\rho,r} : \R \times (M\times \R) \to \R, \qquad  
(s,p,t)\mapsto \rho(s-r)H_{s,\infty}(t). 
$$ 
Given $\og'\in\cP_\lambda^{\le\alpha}$, $\ug\in\cP(H^{\alpha'})$, 
$p'\in\mathrm{Crit}(f'_\og)$, $q\in \mathrm{Crit}(f_\ug)$, and an 
$s$-dependent family of admissible almost complex structures 
$J^r=J^r_s$, the 
{\bf moduli space of capped punctured $s$-dependent interpolating 
  Morse-Bott broken trajectories}  
$$ 
\cM^A_c(p',q;H_{s,\infty}^{\rho,r},\{f_\gamma,f'_\gamma\},J^r),  
$$ 
is defined similarly to the moduli space of capped punctured 
interpolating Morse-Bott broken trajectories~\eqref{eq:cpiMB}, 
using $H_{s,\infty}^{\rho,r}$ instead of $H_\infty^\rho$. For a 
generic choice of the collection of perfect Morse functions 
$\{f_\gamma, f'_\gamma\}$ and of the $s$-dependent almost 
complex structure $J^r$, this is a smooth manifold of 
dimension   
$$ 
\mu(\og'_p)-\mu(\ug_q) + 2\langle c_1(TW),A\rangle.  
$$ 
 
Let us now consider the moduli space  
$$ 
\cM^A_c := 
\bigcup_{r\in \R}  
\cM^A_c(p',q;H_{s,\infty}^{\rho,r},\{f_\gamma,f'_\gamma\},J^r) 
$$ 
for  
$\mu(\og'_p)-\mu(\ug_q) + 2\langle c_1(TW),A\rangle=0$ and 
$\og'\in\cP_\lambda^{\le \alpha}$, $\ug\in\cP(H^{\alpha'})$.  
For a generic choice of a smooth $1$-parameter family of almost complex 
structures $J^r$, the space $\cM^A_c$ is a smooth 
$1$-dimensional manifold, and its boundary splits into a disjoint 
union  
$$ 
\partial \cM^A_c = \partial ^+\cM^A_c \cup \partial ^- 
\cM^A_c \cup \partial ^0 \cM^A_c, 
$$ 
where $\partial^ \pm\cM^A_c$ correspond to $r\to 
\pm\infty$ and $\partial^0 \cM^A_c$ corresponds to finite values of 
$r$. Since $H_{s,\infty}^{\rho,r}=H_\infty^{\alpha',\rho(\cdot - r)}$ 
for $r>>0$, the set $\partial ^+\cM^A_c$ is in bijective 
correspondence with    
$$ 
\cM^A_c(p',q;H_\infty^{\alpha',\rho(\cdot - 
r)},\{f_\gamma,f'_\gamma\},J) 
$$  
and the count of its elements 
gives rise to a morphism which is chain homotopic to 
$\Phi^{\alpha'}\circ \kappa_{\alpha',\alpha}$ 
in~\eqref{eq:last_diag*} (the chain homotopy comes from the need to 
further deform $\rho(\cdot- r)$ to $\rho$). On the other hand,  
for $r<<0$ the count of elements in $\partial^-\cM^A_c$ gives rise to 
a morphism which is chain homotopic to the morphism obtained by 
counting the elements of  
$$ 
\cM^{A_1}_c(p',q_1;H_\infty^{\alpha,\rho(\cdot-r)},\{f_\gamma,f'_\gamma\},J) 
\times \cM^{A-A_1}_c(q_1,q;H_{s,\infty},\{f_\gamma\},J). 
$$ 
Here $q_1\in \mathrm{Crit}(f_{\gamma_1})$ satisfies the equation 
$\mu(\og')-\mu(\gamma_1)  
+ \ind(p') - \ind(q_1)+2\langle c_1(TW),A_1\rangle=0$.  
This last morphism is in turn  
chain homotopic to $\sigma_{\alpha',\alpha}\circ \Phi^\alpha$.  
Finally, the set $\partial ^0\cM^A_c$ is in bijective correspondence with 
the union of  
$$ 
\cM^{A_1}_c(p',q_1;H_{s,\infty}^{\rho,r},\{f_\gamma,f'_\gamma\},J^r) 
\times \cM^{A-A_1}_c(q_1,q;H_\infty^{\alpha'},\{f_\gamma\},J) 
$$ 
and  
$$ 
\cM^{A_1}_c(p',q'_1;\{f'_\gamma\},J_\infty) \times 
\cM^{A-A_1}_c(q'_1,q;H_{s,\infty}^{\rho,r},\{f_\gamma,f'_\gamma\},J^r),  
$$ 
for $r\in \R$, $q_1\in\mathrm{Crit}(f_{\gamma_1})$, $q'_1\in 
\mathrm{Crit}(f'_{\gamma'_1})$ such that 
$\mu(\og')-\mu(\gamma_1)+\ind(p')-\ind(q_1)+2\langle 
c_1(TW),A_1\rangle =-1$, respectively 
$\mu(\gamma'_1)-\mu(\ug)+\ind(q'_1)-\ind(q)+2\langle 
c_1(TW),A-A_1\rangle =-1$.  
 
The count of the elements of the above  
moduli spaces of index $-1$ gives rise to a chain homotopy 
between the morphisms corresponding to the count of elements of 
$\partial^+\cM^A_c$ and $\partial ^-\cM^A_c$ respectively. More 
precisely, let  
$$ 
K:BC_*^{i^{-1}(a),\le\alpha}(\lambda)\to BC_{*+1}^a(H^{\alpha'}) 
$$  
be defined by  
$$ 
K(\og'_p)=\sum_{r\in\R} \  \sum_{F \in 
  \cM^A_c(p',q;H^{\rho,r}_{s,\infty},\{f_\gamma,f'_\gamma\},J^r)}   
\oeps(F) e^A \ug_q, 
$$ 
where the second sum runs over elements such that 
$\mu(\og')-\mu(\ug)+\ind(p')-\ind(q)+2\langle c_1(TW), A\rangle 
=-1$. It follows from our transversality assumptions and the 
discussion above that there are only a finite number of values of the 
parameter $r\in\R$ for which the second sum is nonempty.  
 
We claim that the chain homotopies constructed above have order 
$1$. We argue only for $K$, the other cases being similar.  
We need to show that $\mu(\og)-\mu(e^A\ug)\ge -1$. The moduli spaces 
involved in the  
definition of 
$\cM^A_c(p',q;H^{\rho,r}_{s,\infty},\{f_\gamma,f'_\gamma\},J^r)$ are 
of one of the following three types: 
\begin{itemize}
\item $\cM^{A_i}_c(S'_{\gamma_{i-1}},S'_{\gamma_i};J_\infty)$ with
$i<0$, 
\item
$\cM^{A_0}_c(S'_{\gamma_{-1}},S_{\gamma_0};H^{\rho,r}_{s,\infty},J^r)$,
\item $\cM^{A_j}_c(S_{\gamma_{j-1}},S_{\gamma_j};H^{\alpha'}_\infty,J)$ 
with $j>0$. 
\end{itemize}
The moduli spaces of the first and third type are regular, so 
that $\mu(\gamma'_{i-1})-\mu(e^{A_i}\gamma'_i)\ge 0$, and  
$\mu(\gamma_{j-1})-\mu(e^{A_j}\gamma_j)\ge 0$. On the other hand, due 
to the presence of the parameter $r\in\R$, the 
index of the Fredholm problem for the moduli space of the second type 
is $1$ bigger than the index of the Fredholm problem obtained by replacing 
$H^{\rho,r}_{s,\infty}$ with $H^\rho_\infty$ (see also the 
definition of $\Phi$). Since our moduli spaces are regular and have 
dimension at least $1$ (see Remark~\ref{rmk:d0}), we infer that 
$\mu(\gamma'_{-1}) - \mu(e^{A_0}\gamma_0)\ge -1$. This proves the 
claim, and the Proposition.   
\hfill{$\square$}

 
\section{Proof of the Main Theorem}  
 
\subsection{The long exact sequence} \label{sec:proofleq}  
 
\proof[Proof of Theorem~\ref{thm:intro}] 
 Let $\alpha>0$ be such that $\alpha \notin 
 \mathrm{Spec}(M,\lambda)$. Let $H\in \cH'$ be an admissible 
 Hamiltonian of maximal slope $\alpha$. By  
 Sections~\ref{sec:slow} and~\ref{sec:stretch} and with notation as there, 
 we modify $H$ to $H^{R,\tau}$ and the almost complex structure $J$ 
 to $J^\tau$ which is regular for Floer's equation    
 so that Proposition~\ref{prop:reg} holds, and so that the Floer 
 trajectories are close to punctured Floer trajectories in the 
symplectization, capped with rigid holomorphic planes in $\widehat W$. 
We have performed the ``slowing down'' and ``stretch of 
 the neck'' operations separately in order to emphasize the key ideas 
 for each of them, but it is clear that they can be performed 
 simultaneously in order to obtain such a $H^{R,\tau}$. We denote in 
 the sequel $H=H^{R,\tau}$ and $J=J^\tau$, with both parameters $R$, 
 $\tau$ being large enough.  
 
The spectral sequence $(E^r_d,\bar d^r)$ for symplectic homology is 
supported in two lines. Its $E^2$ page has the form  
$$ 
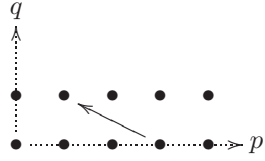
\xymatrix 
@C=7pt  
@R=7pt@W=1pt@H=1pt 
{ q & & & & & \\ 
 &  &  &  &  & \\ 
\bullet & \bullet & \bullet & \bullet & \bullet & \\ 
\bullet \ar@{.>}[rrrrr] \ar@{.>}[uuu] & \bullet & \bullet & \bullet \ar[ull] 
& \bullet & p  
} 
$$  
and the only possibly nontrivial differentials are $\bar d^2:E^2_{k,0}\to 
E^2_{k-2,1}$. By definition of $E^3$ we have exact sequences  
$$ 
0\to E^3_{k,0} \to E^2_{k,0} \stackrel {\bar d^2} \to E^2_{k-2,1} \to 
E^3_{k-2,1} \to 0. 
$$ 
The spectral sequence converges to $SH_*^a(H,J)$ and 
$E^3_d=E^\infty_d$ for dimensional reasons, so that 
we have exact sequences   
$$ 
0\to E^3_{k-1,1} \to SH_k^a(H,J) \to E^3_{k,0} \to 0 
$$ 
by definition of convergence. This information can be put together 
into a long exact sequence by discarding the $E^3$ terms  
$${ 
\def\objectstyle{\scriptstyle} 
\def\labelstyle{\scriptstyle} 
\xymatrix 
@C=4.5pt 
@R=10pt@W=1pt@H=1pt 
{ 
E^2_{k-1,1} \ar@{.>}[rr] \ar[dr] & & SH_k^a(H,J) \ar@{.>}[rr] 
\ar[dr] & & E^2_{k,0} \ar[rr]^{\bar d^2} & & E^2_{k-2,1} \ar@{.>}[rr] 
\ar[dr]  
& & SH_{k-1}^a(H,J) \\ 
& E^3_{k-1,1} \ar[ur] \ar[dr]-<3mm,0mm> \ar@{<-}[dl]+<0mm,0mm> 
& & E^3_{k,0}  
\ar@{<-}[dl]+<4mm,0mm> \ar[ur] \ar[dr]+<1mm,0mm>  
& & & & E^3_{k-2,1} \ar[ur] \ar[dr]-<3mm,0mm>  
\ar@{<-}[dl]-<1.6mm,0mm> & \\ 
0 \quad & & \ 0 \ \ 0 & & \qquad 0 & & 0 \qquad & & \! \! 0  
} 
} 
$$ 
 
We have already seen in Corollary~\ref{cor:spectral} that the map 
 $\Phi$ induces an isomorphism of spectral sequences 
 $\Phi:(E^r_\delta,\bar \delta^r)\stackrel \sim \to (E^r_d,\bar d^r)$, $r\ge 0$. 
In particular we have the following commutative diagram, with vertical 
 arrows being isomorphisms  
$$ 
  \xymatrix{ 
E^2_\delta \ar[r]^{\bar \delta^2} \ar[d]^\simeq_{\Phi} & E^2_\delta 
  \ar[d]_\simeq^{\Phi} \\  
E^2_d \ar[r]^{\bar d^2} & E^2_d  
} 
$$ 
Combining the isomorphism~\eqref{eq:E2delta}   
$$ 
\bar \Theta:E^2_\delta \stackrel \sim \longrightarrow 
HC_*^{i^{-1}(a),\le\alpha}(\lambda,J) \otimes H_*(S^1) 
$$ 
with the two previous diagrams we get a long exact sequence 
\begin{eqnarray} \label{eq:seqalpha} 
\hspace{-8mm}\dots \!\to \!SH_k^a(H,J) \!\to\!  
 HC_{k+(n-3)}^{i^{-1}(a),\le\alpha}(\lambda,J) 
 \!\stackrel D \to\! HC_{k-2+(n-3)}^{i^{-1}(a),\le\alpha}(\lambda,J)  
 \!\to \\ 
 \hspace{5cm}\quad \to  
 SH_{k-1}^a(H,J) \to \dots & \nonumber  
\end{eqnarray} 
with $D=\bar\Theta \circ \bar \delta^2 \circ \bar\Theta^{-1}$. The shift in degree 
is due to the fact that $\Phi$  
decreases degrees by $n-3$.  
Since the limiting slope of $H$ equals $\alpha$ we have  
$$ 
SH_*^a(H,J)\simeq SH_*^{a,\le\alpha}(W,\om), 
$$ 
 where the 
 latter notation stands for a direct limit on $\cH'$ of Floer homology 
 groups truncated by the values of the action functional in the range 
 $]-\infty,\alpha]$.  
 
We claim that the exact sequences~\eqref{eq:seqalpha}  
form a natural directed system, i.e. for 
$\alpha<\alpha'\notin\mathrm{Spec}(M,\lambda)$ the continuation  
morphisms $\sigma_{\alpha',\alpha}$ in symplectic homology induced by 
an increasing homotopy of Hamiltonians as in Section~\ref{sec:MBsymp}, 
and the extension morphisms $\kappa_{\alpha',\alpha}$ in linearized  
contact homology induced by inclusion, fit into a commutative  
diagram   
\begin{equation} \label{eq:seqdirect} 
\xymatrix 
@C=11pt 
{ 
\hspace{-6mm} 
\cdots  SH_k^a(H^\alpha,J^\alpha) \ar[r] \ar[d]^{\sigma_{\alpha',\alpha}} &   
HC_{k+(n-3)}^{i^{-1}(a),\le\alpha}(\lambda,J) \ar[r]^-{D}  
\ar[d]^{\kappa_{\alpha',\alpha}} & 
HC_{k-2+(n-3)}^{i^{-1}(a),\le\alpha}(\lambda,J)  
\ar[d]^{\kappa_{\alpha',\alpha}}  \, \cdots \\ 
\hspace{-6mm} 
\cdots  SH_k^a(H^{\alpha'},J^{\alpha'}) \ar[r] & 
HC_{k+(n-3)}^{i^{-1}(a),\le\alpha'}(\lambda,J) \ar[r]^-{D} & 
HC_{k-2+(n-3)}^{i^{-1}(a),\le\alpha'}(\lambda,J)  \, \cdots 
} 
\end{equation} 
Here $H^\alpha=H^{R(\alpha),\tau(\alpha)}$, 
$H^{\alpha'}=H^{R(\alpha'),\tau(\alpha')}$, 
$J^\alpha=J^{\tau(\alpha)}$ and $J^{\alpha'}=J^{\tau(\alpha')}$ are 
as explained in the beginning of the proof. The asymptotic slope of 
$H^\alpha$ is equal to $\alpha$, the asymptotic slope of 
$H^{\alpha'}$ is equal to $\alpha'$ and, by 
Proposition~\ref{prop:stretch}, we can assume 
without loss of generality   
that $\tau(\alpha)=\tau(\alpha')=\tau$ so that the 
Hamiltonians and almost complex structures are defined on $\widehat 
W^\tau$ as in Section~\ref{sec:stretch}. Moreover, we can also assume 
without loss of generality that $H^\alpha\le H^{\alpha'}$.  
 
By Remark~\ref{rmk:sdep} in Section~\ref{sec:MBcomplex}, the 
continuation morphism  
 $$ 
 \sigma_{\alpha',\alpha}:BC_*^a(H^\alpha)\to BC_*^a(H^{\alpha'}) 
 $$ 
preserves the filtrations, and hence induces a morphism of 
spectral sequences   
$$ 
\sigma_{\alpha',\alpha}: (E^r_d(H^\alpha),\bar d^r) \to 
(E^r_d(H^{\alpha'}),\bar d^r), \qquad r\ge 0.  
$$ 
The inclusion 
$\kappa_{\alpha',\alpha}:BC_*^{i^{-1}(a),\le\alpha}(\lambda)\to 
BC_*^{i^{-1}(a),\le\alpha'}(\lambda)$ 
preserves the filtrations by definition, and therefore also induces a 
morphism between the associated spectral sequences  
$$ 
\kappa_{\alpha',\alpha}:(E^{r,\le\alpha}_\delta,\bar \delta^r)\to 
(E^{r,\le\alpha'}_\delta,\bar \delta^r), \qquad r\ge 0. 
$$ 
In view of Proposition~\ref{prop:homotopy}, commutativity 
of~\eqref{eq:seqdirect} follows from the commutativity up to a chain 
homotopy of order $1$ of the following diagram of 
morphisms of filtered complexes 
\begin{equation}\label{eq:last_diag} 
\xymatrix 
@C=25pt 
@R=15pt 
{BC_*^{i^{-1}(a),\le \alpha}(\lambda) \ar[r]^-{\Phi^\alpha}  
\ar[d]_-{\kappa_{\alpha',\alpha}} & 
BC_*^a(H^\alpha) \ar[d]^-{\sigma_{\alpha',\alpha}} \\ 
BC_*^{i^{-1}(a),\le \alpha'}(\lambda)  \ar[r]_-{\Phi^{\alpha'}} &  
BC_*^a(H^{\alpha'}) 
} 
\end{equation}  
Here $\Phi^\alpha$ and $\Phi^{\alpha'}$ are the filtered isomorphisms 
constructed in Proposition~\ref{prop:filtered}, for the Hamiltonians 
$H^\alpha$, respectively $H^{\alpha'}$.  
 
The commutativity of~\eqref{eq:last_diag} up to a chain 
homotopy of order $1$ is precisely the content of  
Proposition~\ref{prop:commdiag} above. Therefore~\eqref{eq:seqdirect} is  
commutative and, passing to the direct limit on 
 $\alpha$, we obtain an exact sequence   
\begin{eqnarray} \label{eq:mainthmwithlambdaJ} 
\hspace{-7mm}\ldots \!\to\! SH_k^a(W,\om) \!\to\! 
HC_{k+(n-3)}^{i^{-1}(a)}(\lambda,J)  
\!\stackrel D \to\! HC_{k-2+(n-3)}^{i^{-1}(a)}(\lambda,J) \!\to \\ 
\hspace{4.5cm} \to  
SH_{k-1}^a(W,\om) \to \ldots & \nonumber 
\end{eqnarray}   
Here we used the fact that the direct limit functor is exact.  
By the invariance of linearized contact homology with respect to 
the pair $(\lambda,J)$ (see Remark~\ref{rmk:polyfolds}) we obtain the 
exact sequence in the statement of Theorem~\ref{thm:intro}.  
 
Finally, the description of the differential $D$ claimed in the 
statement of Theorem~\ref{thm:intro} is the content of 
Proposition~\ref{prop:D} below. This finishes the proof. 
\hfill{$\square$}  
 
\begin{remark} \label{rmk:strongerstatement}  
  {\rm  
What we have actually proved is that $\Phi$ induces an isomorphism of 
degree $3-n$ between the exact sequences  
{\scriptsize 
$$ 
\xymatrix 
@C=7pt 
{... \ar@{.>}[r]  & H_*(BC_*^{i^{-1}(a)}(\lambda),\delta) 
  \ar@{.>}[r] \ar@{.>}[d]^\Phi & HC_*^{i^{-1}(a)}(\lambda) \ar[r]^D 
  \ar[d]^{\Phi\circ \bar\Theta^{-1}} 
  & HC_{*-2}^{i^{-1}(a)}(\lambda) \ar@{.>}[r] \ar[d]^{\Phi\circ\bar\Theta^{-1}} & 
  H_{*-1}(BC_*^{i^{-1}(a)}(\lambda),\delta) \ar@{.>}[r] \ar@{.>}[d]^\Phi & 
  ... \\ 
... \ar[r] & SH_*(W,\om) \ar[r] & E^2_{d;*,0} \ar@{.>}[r]^{\bar d^2} & 
E^2_{d;*-2,1} \ar[r] & SH_{*-1}(W,\om) \ar[r] & ... 
} 
$$ 
} 
 
In order to establish Theorem~\ref{thm:intro} we have only used the 
undotted arrows in the above diagram. We shall not explain in this 
paper the significance of the discarded terms $E^2_d$ and 
$H_*(BC_*^{i^{-1}(a)}(\lambda),\delta)$. This will be the topic of the 
forthcoming papers~\cite{BO} and~\cite{CO} respectively.  
} 
\end{remark}

\subsection{The differential $D$} \label{sec:D}  
 
The purpose of this section is to give a description of $D$ which does 
not make use of the auxiliary Morse functions $f'_\gamma$, and  
thus  
complete the proof of Theorem~\ref{thm:intro}.   
 
Given half-lines $\oL\subset T_0\C P^1$, $\uL\subset T_\infty \C P^1$, 
we define half-lines $\oL_\infty\subset T_\infty \C P^1$, 
$\uL_0\subset T_0\C P^1$ by choosing a global polar coordinate on $\C 
P^1\setminus \{0,\infty\}$ and requiring  
$$ 
\Arg(\oL_\infty)=\Arg(\oL), \qquad \Arg(\uL_0)=\Arg(\uL). 
$$ 
Given half-lines $L_0\subset T_0\C P^1$, $L_\infty\subset T_\infty \C 
P^1$, and a map $F=(f,a):\R\times S^1=\C P^1\setminus\{0,\infty\}\to 
M\times \R$ satisfying~\eqref{eq:a}, \eqref{eq:f}, we define 
$\ev(L_0)=\lim_{z\to 0,\ z\in L_0} f(z)$, $\ev(L_\infty)=\lim_{z\to 
  \infty,\ z\in L_\infty} f(z)$, so that $\ev(L_0)$ belongs to the 
geometric image of $\og'$ and $\ev(L_\infty)$ belongs to the geometric 
image of $\ug'$. We also recall that we have chosen a point 
$P_{\gamma'}$ on the geometric image of each $\gamma'\in\cP_\lambda$.  
 
Given 
$\og',\ug'\in\cP_\lambda$, $A\in H_2(W;\Z)$ we denote by   
$$ 
\cM^A_{1,c}(P_{\og'},P_{\ug'}; J) \subset \cM^A_c(\og',\ug';J) 
$$ 
the subset of equivalence classes of pairs 
$[u',F']\in\cM^A_c(\og',\ug';J)$ such that the asymptotic markers 
$\oL$ at $0$ and $\uL$ at $\infty$ satisfy  
$$ 
\oL_\infty=\uL, \mbox{ or equivalently } \uL_0=\oL. 
$$ 
The decoration ``$1$'' for the moduli space 
is motivated by the fact that it consists of curves with one sublevel.  
 
Given $\og',\ug'\in\cP_\lambda$, $A\in H_2(W;\Z)$ we denote by   
$$ 
\cM^{A,+}_{2,c}(P_{\og'},P_{\ug'};J) \subset 
\hspace{-.5cm}  
\bigcup_{\gamma'\in\cP_\lambda, A_1\in H_2(W;\Z)} 
\hspace{-.5cm}  
\big[\cM^{A_1}_c(\og',\gamma';J)\big] \times 
\big[\cM^{A-A_1}_c(\gamma',\ug';J) \big] 
$$ 
the subset of pairs of equivalence classes $([u'], [u''])$ 
for the equivalence relation given by ignoring the asymptotic markers 
$\uL'$, $\oL''$ corresponding to the common asymptote $\gamma'$,  
such that the cyclic order of the points  
$(P_{\gamma'},\ev(\oL'_\infty),  
\ev(\uL''_0))$ is the same as the natural  
orientation of the geometric image of $\gamma'$. The decorations 
``$2$'' and~``$+$'' for the moduli space are motivated by the fact that it 
consists of curves with two sublevels and satisfying an additional 
``positive'' cyclic order condition.  
In the situation  
$\mu(\og')-\mu(\ug')+2\langle c_1(TW),A\rangle =2$ and for a generic 
choice of the points $P_{\gamma'}$ the moduli spaces 
$\cM^A_{1,c}(P_{\og'},P_{\ug'};J)$ and $\cM^{A,+}_{2,c}(P_{\og'},P_{\ug'};J)$ 
are rigid and one can associate a sign $\epsilon(u)$ to each of their 
elements via coherent orientations and fibered products.  
 
For each free homotopy class $a$ in $W$ we define a 
map   
\begin{equation*}  
 \Delta: C_*^{i^{-1}(a)}(\lambda) \to C_{*-2}^{i^{-1}(a)}(\lambda), 
\end{equation*}  
\begin{equation} 
  \label{eq:Delta} 
  \Delta(\og') = \sum_{\substack{\ug',A \\ 
  |\ug' e^A| = |\og'| -2 } } 
\frac 1 {{\scriptstyle \kappa_{\ug'}}} 
\sum_{u\in \cM^A_{1,c}(P_{\og'},P_{\ug'};J) \cup 
  \cM^{A,+}_{2,c}(P_{\og'},P_{\ug'};J)}  
  \epsilon(u) e^A \ug'. 
\end{equation}

\begin{proposition} \label{prop:D}  
 The map $\Delta$ defined by~\eqref{eq:Delta} is a chain map, and 
 induces in homology the map $D$ in the long exact sequence of 
 Theorem~\ref{thm:intro}.     
\end{proposition}  
 
\proof Let us first reinterpret the previous 
moduli spaces in terms of moduli spaces of capped punctured 
$S^1$-parametrized holomorphic cylinders. Given 
$\og',\ug'\in\cP_\lambda$, $A\in H_2(W;\Z)$ we denote by   
$$ 
\widetilde \cM^A_{1,c}(P_{\og'},P_{\ug'}; J) \subset 
\cM^A_c(S'_\og,S'_\ug;J)  
$$ 
the subset consisting of pairs $u=(u',F')\in\cM^A_c(S'_\og,S'_\ug;J)$ such 
that  
$$ 
\oev(u')=P_{\og'}, \qquad \uev(u')=P_{\ug'}. 
$$ 
It follows from the definition that there is a bijective 
correspondence  
$$ 
\cM^A_{1,c}(P_{\og'},P_{\ug'}; J) \sim \widetilde 
\cM^A_{1,c}(P_{\og'},P_{\ug'}; J).  
$$ 
Given $\og',\ug'\in\cP_\lambda$, $A\in H_2(W;\Z)$ we denote by   
$$ 
\widetilde \cM^A_c(P_{\og'},S'_\ug ; J), \widetilde \cM^A_c(S'_\og,P_{\ug'} ; J) 
\subset \cM^A_c(S'_\og,S'_\ug;J)  
$$ 
 the subsets consisting of pairs $u=(u',F')\in\cM^A_c(S'_\og,S'_\ug;J)$ 
 such that  
$$ 
\oev(u')=P_{\og'}, \mbox{ respectively } \uev(u')=P_{\ug'}. 
$$ 
If $\mu(\og')-\mu(\ug')+2\langle c_1(TW),A\rangle =1$ these moduli spaces 
are rigid. Let  
$$ 
\widetilde \cM^{A,+}_{2,c}(P_{\og'},P_{\ug'};J) \subset 
\hspace{-.5cm} 
\bigcup_{\gamma'\in\cP_\lambda, A_1\in H_2(W;\Z)}  
\hspace{-.5cm} 
\widetilde \cM^{A_1}_c(P_{\og'},S'_\gamma;J) \times 
\widetilde \cM^{A-A_1}_c(S'_\gamma,P_{\ug'};J)  
$$ 
be the subset consisting of pairs $(\overline u, \underline u)$ such 
that the cyclic order of the points $(P_{\gamma'},\uev(\overline u), 
\oev(\underline u))$ is the same as the one induced by the chosen 
orientation of $S'_\gamma$. It follows from the definition that there 
is a bijective correspondence  
$$ 
\cM^{A,+}_{2,c}(P_{\og'},P_{\ug'};J)\sim \widetilde 
\cM^{A,+}_{2,c}(P_{\og'},P_{\ug'};J).  
$$ 
Hence $\widetilde \Delta:=\Theta^{-1} \circ \Delta \circ \Theta: 
C_*^{i^{-1}(a)}(\lambda)\otimes H_0(S^1)\to 
C_{*-2}^{i^{-1}(a)}(\lambda)\otimes H_1(S^1)$, where $\Theta$ is 
defined in~\eqref{eq:Theta}, acts by  
$$ 
\widetilde \Delta(\og'_M) = \sum_{\substack{\ug',A \\ 
  |\ug' e^A| = |\og'| -2 } } 
\sum_{u\in \widetilde \cM^A_{1,c}(P_{\og'},P_{\ug'};J) \cup 
  \widetilde \cM^{A,+}_{2,c}(P_{\og'},P_{\ug'};J)}  
  \epsilon(u) e^A \ug'_m. 
$$ 
 
  For a generic choice of Morse functions $f'_\gamma$, 
  $\gamma'\in\cP_\lambda$ the map $\bar\delta^2$ is induced in homology by the 
  map $\delta^2$ in the decomposition of the $S^1$-parametrized 
  differential $\delta$, and does not depend on 
  the choice of the collection $\{f'_\gamma\}$. It is therefore 
  enough to show that, for a suitable choice of this collection, the 
  map induced by $\delta^2$ on  
  $(E^1_\delta,\bar \delta^1)$ is $\widetilde \Delta$ itself.  
 
  Let us fix $\alpha>0$ such that 
  $\alpha\notin\mathrm{Spec}(M,\lambda)$. Let $\widetilde \Delta^\alpha$ 
  and $\delta^{2,\alpha}$ be the truncations of $\widetilde \Delta$ and 
  $\delta^2$ to action less than $\alpha$. It is enough to show 
  that $\widetilde \Delta^\alpha=\delta^{2,\alpha}$ for a suitable choice of the 
  collection $\{f'_\gamma\}$ which depends on $\alpha$. By letting 
  $\alpha\to\infty$ we then get $\widetilde \Delta=\delta^2$.  
 
  The set $\cP_\lambda^{\le\alpha}$ is finite and, for each pair 
  $\og',\ug'\in \cP_\lambda^{\le\alpha}$, the moduli space of holomorphic 
  curves asymptotic to $\og'$,$\ug'$ is compact and therefore involves 
  only a finite number of homology classes. As a consequence, for each 
  $\gamma'\in\cP_\lambda^{\le\alpha}$ we can choose an open neighbourhood 
  $\cU_\gamma^{\prime,\alpha}$ of $P_{\gamma'}$ in $S'_\gamma$ such that  
  every collection $\{q_{\gamma'}\}\in\prod_{\gamma'\in\cP_\lambda^{\le\alpha}} 
  \cU_\gamma^{\prime,\alpha}$ is regular and the map 
  $\widehat\Delta^\alpha$ associated to $\{q_{\gamma'}\}$ is equal to 
  the map $\widetilde \Delta^\alpha$ associated to $\{P_{\gamma'}\}$.  
         
  By choosing a generic collection $\{q_{\gamma'}\}$, 
  $q_{\gamma'}\in\cU_\gamma^{\prime,\alpha}$, 
  $\gamma'\in\cP_\lambda^{\le\alpha}$ and small neighbourhoods  
  $\cV_\gamma^{\prime,\alpha}\subset 
  \cU_\gamma^{\prime,\alpha}$ of $q_{\gamma'}$,  
  we can further assume that the evaluation maps $\uev$, $\oev$ 
  defined on the spaces $\cM^A_c(q_{\og'},S'_\ug;J)$ and 
  $\cM^A_c(S'_\og,q_{\ug'};J)$ respectively, with 
  $|\og'|-|e^A\ug'|=1$, miss the neighbourhoods $\cV_\gamma^{\prime,\alpha}$.   
 
  We choose the Morse functions 
  $f'_\gamma$, $\gamma'\in\cP_\lambda^{\le\alpha}$ generically so that 
  both critical points $m'$ and $M'$ of $f'_\gamma$ lie inside 
  $\cV_\gamma^{\prime,\alpha}$ and so that the ``long arc'' 
  in $S'_\gamma$ running from $m'$ to $M'$ and exiting 
  $\cV_\gamma^{\prime,\alpha}$ has the same orientation as the chosen 
  orientation of $S'_\gamma$. Let us show that in this situation we 
  have $\widehat \Delta^\alpha=\delta^{2,\alpha}$.    
 
  We first note that $\delta^{2,\alpha}$ is built out of two kinds of moduli spaces 
of capped punctured $S^1$-parametrized broken $J_\infty$-holomorphic 
cylinders, namely having either one or two sublevels. Indeed, we cannot have 
more than two sublevels since the difference of indices at the 
extremities is equal to two, and each sublevel introduces a difference 
of index of at least one due to the fact that the moduli spaces of 
punctured $S^1$-parametrized $J_\infty$-holomorphic cylinders are 
regular and carry a one-dimensional symmetry given by the action of $S^1$.  
 
  Let us fix $\og',\ug'\in\cP_\lambda^{\le\alpha}$ and denote by $M'$ the maximum of 
$f'_\og$ and by $m'$ the minimum of $f'_\ug$. Since 
$M'\in\cV_\og^{\prime,\alpha}\subset \cU_\og^{\prime,\alpha}$ 
and $m'\in\cV_\ug^{\prime,\alpha}\subset \cU_\ug^{\prime,\alpha}$, the 
elements of $\cM^A_c(M',m';\{f'_\og\},J)$  
having one sublevel are in one-to-one correspondence with elements of 
$\widetilde \cM^A_{1,c}(q_{\og'},q_{\ug'};J)$ and their signs coincide.  
 
We claim now that the subset of $\cM^A_c(M',m';\{f'_\gamma\},J)$  
consisting of elements with two sublevels is in bijective 
correspondence with $\widetilde \cM^{A,+}_{2,c}(q_{\og'},q_{\ug'};J)$, with 
preservation of signs. Indeed, such elements are of the form  
$$ 
(\overline u,\underline 
u) \in \cM^{A_1}_c(M',S'_\gamma;J)\times 
\cM^{A-A_1}_c(S'_\gamma,m';J) 
$$ 
so that there is a gradient trajectory of $f'_\gamma$ running from 
$\uev(\overline u)$ to $\oev(\underline u)$. Since these evaluation 
maps miss the neighbourhood $\cV_\gamma^{\prime,\alpha}$ 
and by our choice of order for the two critical points of $f'_\gamma$, 
this is equivalent to saying that the cyclic order on the triple 
$(q_{\gamma'},\uev(\overline u),\oev(\underline u))$ is the same as the 
one induced by the chosen orientation of $S'_\gamma$. Since $M'$ is 
close to $q_{\og'}$ and $m'$ is close to $q_{\ug'}$ there is a unique 
element in $\widetilde \cM^{A,+}_{2,c}(q_{\og'},q_{\ug'};J)$ corresponding to such a 
pair $(\overline u,\underline u)$. Their signs coincide for continuity 
reasons and this proves our claim.  
 
We have shown that $\delta^{2,\alpha}=\widehat 
\Delta^\alpha$ on $(E^1_\delta,\bar \delta^1)$, as desired. 
\hfill{$\square$}

 
\section{Examples} \label{sec:examples} 
 
\subsection{Riemann surfaces} \label{sec:Riemann}  
 
We compute in this section the exact sequence~\eqref{eq:intro} for 
genus $g$ Riemann surfaces $\Sigma=\Sigma_{g,1}$ with one boundary 
component. We shall see that the cases $g=0$ and $g\ge 1$ are 
fundamentally different: although the boundary $M=\p \Sigma$ is the 
same, i.e. the circle $S^1$, the linearized contact homology groups 
$HC_*(M)$ differ, and so do the corresponding maps $D$. Note that 
regularity is automatic when the target manifold is a Riemann 
surface. 
 
Let us note that free 
homotopy classes of loops in $M=S^1$ are indexed by $\Z$ via the 
degree of the corresponding maps $S^1\to  
S^1$. Given a free homotopy class $b$ of loops we denote by the same 
symbol $b\in \Z$ its degree. There are no closed Reeb 
orbits in any class $b\le 0$, whereas each class $b\in\Z^+$ 
contains exactly one closed Reeb orbit $\gamma_b$. Since the contact 
distribution is zero-dimensional we need to use the special  
convention $\mu(\gamma_b)=2b$ for the Maslov index, corresponding to 
the index of the linearized Reeb flow in the symplectization. 
 
Let us first consider the case $g\ge 1$. The inclusion $i$ of free 
homotopy classes of loops from the boundary to $\Sigma$ is injective 
and we denote $i(b)$ by $b$. For each $b\in\Z^+$ we have 
$HC_*^b(M)=\Q$ if $*=2b+(1-3)=2b-2$ and  
$0$ otherwise. The Reeb orbit $\gamma_b$ gives rise to two Hamiltonian 
orbits of  
indices $2b$ and $2b+1$ and we have $SH_*^b(\Sigma)=\Q$ if $*=2b,2b+1$ 
and $0$ otherwise. The interesting portion of the exact sequence is 
therefore  
$${
\def\objectstyle{\scriptstyle} 
\def\labelstyle{\scriptstyle} 
\xymatrix 
@C=15pt 
@R=10pt 
{ 
0\ar[r]^-{D} & HC_{2b-2}^b \ar[r] \ar@{=}[d] & SH_{2b+1}^b(\Sigma) \ar[r] 
\ar@{=}[d] & 0 \ar[r]^-{D} & 
0 \ar[r] & SH_{2b}^b(\Sigma) \ar[r] \ar@{=}[d] & HC_{2b-2}^b(M) \ar[r]^-{D} 
\ar@{=}[d] & 0 \\ 
& \Q & \Q &  &  & \Q & \Q &  
} 
} 
$$ 
We see in particular that $D$ vanishes.  
 
We now consider the case $g=0$, so that we can assume without loss of 
generality that $\Sigma=D^2$, the unit disc in the complex plane. Note 
that all closed Reeb orbits are contractible in $D^2$. Since 
$SH^0_*(D^2)=0$ we have $SH^+_*(D^2)\simeq H_*(D^2,\p D^2) = \Q$ if 
$*=2$ and $0$ otherwise. The linearized contact 
complex is lacunary with generators of index $2b-2$, $b\ge 1$ and 
therefore $HC_*^{i^{-1}(0)}=\Q$ if $*=2b-2$, $b\ge 1$ and $0$ 
otherwise. The long exact sequence~\eqref{eq:intro} therefore splits 
into short exact sequences of which the nontrivial ones are  
$${ 
\def\objectstyle{\scriptstyle} 
\def\labelstyle{\scriptstyle} 
\xymatrix 
@C=15pt 
@R=10pt 
{0 \ar [r] & SH_2^+(D^2) \ar@{=}[d] \ar[r]^\sim & HC_0(S^1) \ar@{=}[d] \ar[r]^-{D} & 
  HC_{-2}(S^1) \ar@{=}[d] \ar[r] & 0 \\ 
& \Q & \Q & 0 &  
} 
} 
$$ 
and  
$$ { 
\def\objectstyle{\scriptstyle} 
\def\labelstyle{\scriptstyle} 
\xymatrix 
@C=15pt 
@R=10pt 
{0\ar [r] & SH_{2b}^+(D^2) \ar@{=}[d] \ar[r] & HC_{2b-2}(S^1)
  \ar@{=}[d] \ar[r]^-{D}_-{\sim} &  
  HC_{2b-4}(S^1) \ar@{=}[d] \ar[r] & 0 \\ 
& 0 & \Q & \Q &  
}, \quad b\ge 2. 
} 
$$  
We see in particular that the map $D:HC_{2b-2}(S^1) \stackrel \sim \longrightarrow 
HC_{2b-4}(S^1)$ does not vanish for $b\ge 2$. We can actually describe 
it explicitly as follows. The only contractible Reeb orbit of 
normalized index $0$ is $\gamma_1$ in the class $b=1$ and,  
up to reparametrization, there is a unique holomorphic plane in 
$\C=\widehat {D^2}$ asymptotic to it, namely $z\mapsto cz+d$, $c\in 
\C^*$, $d\in\C$. Since there are no rigid  
nontrivial holomorphic cylinders in $S^1\times \R$ the map $D$ is 
obtained by a count of punctured curves in $S^1\times \R$ with only 
one sublevel. These must necessarily have three punctures: a positive 
one asymptotic to $\gamma_b$ and two negative ones asymptotic to 
$\gamma_{b-1}$ and $\gamma_1$. Note that the puncture asymptotic to 
$\gamma_1$ corresponds to the augmentation of 
Remark~\ref{rmk:augmentation} and the count of these curves gives the 
coefficient of $\gamma_{b-1}$ in $D(\gamma_b)$. Such curves correspond to  
meromorphic functions on the Riemann sphere with one pole of order $b$ and two 
zeroes of order $b-1$ and $1$, respectively. Meromorphic functions with 
those properties are unique up to    
reparametrization and thus the sum defining $D$ reduces to 
$D(\gamma_b)=\gamma_{b-1}$, so that $D$ is the obvious isomorphism.

\subsection{Subcritical Stein manifolds} \label{sec:subcritical} 
 
A {\bf Stein manifold} $\widehat W$ is a triple $(\widehat W,J,\phi)$, 
where $J$ is a complex structure on $\widehat W$ and $\phi:\widehat 
W\to \R$ is an exhausting plurisubharmonic function. That $\phi$ is 
exhausting means that $\phi$ is proper and bounded from below. That 
$\phi$ is plurisubharmonic means that $\omega_\phi:= 
-dJ^*d\phi$ is a symplectic form and $J$ is compatible with 
$\omega_\phi$. We say that $\widehat W$ is of {\bf finite type} if we 
can choose $\phi$ such that the set of its critical points is 
compact. In this situation we can assume without loss of generality 
that $\phi$ is Morse~\cite[Theorem~8.1.C]{Bi}. All its 
critical points have index at most $\frac 1 2 
\dim_\R\widehat W$, and we call $\widehat W$ {\bf subcritical} if all 
critical points have index strictly smaller than   
$\frac 1 2 \dim_\R\widehat W$. We assume in this section that  
$\widehat W$ is a subcritical Stein manifold of finite type.  
 
A {\bf Stein domain}  
$W\subset \widehat W$ is a domain such that $W=\{\phi\le c\}$ for some 
$c\in \R$ large enough. In particular $c$ is a regular value of 
$\phi$ and $W$ is a smooth manifold with boundary $\partial 
W=\{\phi=c\}$. Actually $W$ is an exact symplectic manifold with 
boundary of contact type and Liouville vector 
field $\nabla \phi$, the gradient with respect to the metric 
$\omega_\phi(\cdot,J\cdot)$.   
Moreover, the symplectic completion of $W$ is symplectomorphic to 
$(\widehat W,\omega_\phi)$.  
 
The isotopy class of the contact structure $\xi_\phi$ induced on 
$\partial W$ does not depend on $\phi$ and $c$. A contact manifold 
$(M,\xi)$ is called {\bf subcritically Stein fillable} if it can 
be realized as $(\partial W,\xi_\phi)$ for some subcritical Stein 
manifold of finite type $(\widehat W,J,\phi)$. We call the Stein 
domain $W$ a {\bf filling} of $M$.  
 
Mei-Lin Yau has computed in~\cite{Y} the cylindrical contact homology 
groups with rational coefficients in the trivial homotopy class 
for a subcritically Stein fillable contact manifold $(M^{2n-1},\xi)$, 
$n\ge 2$ satisfying $c_1(\xi)=0$. A crucial ingredient in the proof is an  
estimate on the reduced Conley-Zehnder index of 
Reeb orbits running between different handles~\cite[Lemma~4.2]{Y}, 
which implies in particular that every Reeb orbit $\gamma$ which is 
contractible in $W$ satisfies $\bar\mu(\gamma)\ge 1$. Hence there are 
no rigid holomorphic planes in $\widehat W$ and cylindrical contact 
homology is tautologically isomorphic to linearized contact homology 
$HC_*^0(M)$.   
 
The next theorem states Mei-Lin Yau's result in a form which is 
equivalent to the original one of~\cite{Y}.

\begin{theorem}[M.-L.~Yau~\cite{Y}] \label{thm:MLYau} 
 Let $(M^{2n-1},\xi)$, $n\ge 2$ be a subcritically Stein 
fillable contact manifold with $c_1(\xi)=0$, 
and let $W$ be any subcritical Stein filling of $M$. Then, using rational 
coefficients, we have  
\begin{equation} \label{eq:MLYau} 
HC_k^0(M)\simeq \bigoplus _{m\ge 0} H_{k-2m+2}(W,\partial W), \qquad 
k\in \Z.  
\end{equation} 
\end{theorem}  
 
Denote by $i$ the map associating to the free homotopy class 
of a loop in $\partial W$ the free homotopy class of the same loop in 
$W$. We claim that subcriticality implies 
$i^{-1}(0)=0$, i.e. a loop in $\partial W$  
which is contractible in $W$ is actually contractible in $\partial 
W$. Indeed, the homotopy to a point defines a chain of dimension 
$2$. Since $n\ge 2$ and the isotropic skeleton of $W$ is of dimension 
at most $n-1$, we can perturb the homotopy so as to avoid it. Finally 
we push the homotopy to $\partial W$ by the Liouville vector field 
$\nabla \phi$. 
 
We recall now the discussion on regularity in 
Remark~\ref{rmk:transv_exples}. If $W$ is a stabilization  
of a subcritical Stein manifold, i.e. $W=W'\times D^2$ with $W'$ 
subcritical, then all necessary regularity assumptions are met. In 
fact, unpublished work of Cieliebak~\cite{Csplit} shows that  
$W$ is deformation equivalent to such a stabilization if and only if 
it has the homotopy type of a complex of dimension at most $n-2$. The 
exact sequence~\eqref{eq:intro} becomes   
\begin{equation*}  
{\scriptstyle 
\ldots \longrightarrow SH_{k-(n-3)}^+(W) \longrightarrow  
HC_k^0(M) \stackrel D \longrightarrow HC_{k-2}^0(M)\longrightarrow 
SH_{k-1-(n-3)}^+(W) \longrightarrow  
\ldots 
} 
\end{equation*}

\begin{proposition} \label{prop:newMLYau} 
 Let $W$ be a stabilization of a subcritical Stein 
 manifold, and denote $M=\partial W$. The exact 
sequence relating contact and symplectic homology in the trivial 
homotopy class is isomorphic to an exact sequence of the form 
\begin{eqnarray*} 
\lefteqn{\scriptstyle \ldots \stackrel 0 \longrightarrow 
  H_{k+2}(W,\partial W) \longrightarrow   
\bigoplus_{m\ge 0} H_{k-2m+2}(W,\partial W) \longrightarrow  
\bigoplus_{m\ge 1} H_{k-2m+2}(W,\partial W) \stackrel 0 \longrightarrow} \\ 
& \qquad \qquad \qquad \qquad \qquad \qquad \qquad\qquad \qquad \qquad  &  
\scriptstyle \stackrel 0 \longrightarrow H_{k+1}(W,\partial W) 
\longrightarrow \ldots   
\end{eqnarray*} 
\end{proposition} 
 
\proof The symplectic homology groups $SH_*(W)$ of the subcritical 
Stein domain $W$ vanish for the trivial homotopy 
class~\cite{C,Kunneth}, and the tautological long exact 
sequence~\eqref{eq:Vseq} implies  
\begin{equation} \label{eq:SH+} 
SH_{k-(n-3)}^{+}(W)\simeq H_{k+2}(W,\partial W), \qquad k\in \Z. 
\end{equation} 
 
 We prove by induction that the maps  
$HC_k^0(M)\to H_{k+3}(W,\partial W)$, $k\in \Z$ vanish.  
This holds for $k\le n-3$ or $k\ge 2n-2$ because 
$W$ is an oriented manifold of dimension $2n$ having  
the homotopy type of a CW-complex of dimension  
$\le \, n-1$, and therefore $H_{k+3}(W,\partial W)=0$.  
Assuming that the map $HC_k^0(M)\to 
H_{k+3}(W,\partial W)$ vanishes for some $k\in \Z$,  
the exactness of the sequence  
\begin{equation*}\scriptstyle 
 HC_{k+1}^0(M) \to H_{k+4}(W,\partial W) \stackrel i \to  
\bigoplus_{m\ge 0} H_{k-2m+4}(W,\partial W) \to 
\bigoplus_{m\ge 1} H_{k-2m+4}(W,\partial W) \to 0  
\end{equation*} 
shows, for dimensional reasons, that $i$ is necessarily 
injective. Hence the map $HC_{k+1}^0(M) \to H_{k+4}(W,\partial W)$  
vanishes as well and the induction step is completed. 
\hfill{$\square$}

\begin{remark} \label{rmk:subcriticalarrows} {\rm  
Proposition~\ref{prop:newMLYau} shows that our long exact sequence 
splits in the case under study into short exact sequences of 
the form 
$$\scriptstyle 
0 \longrightarrow H_{k+2}(W,\partial W) \stackrel i \longrightarrow  
\bigoplus_{m\ge 0} H_{k-2m+2}(W,\partial W) \stackrel p \longrightarrow  
\bigoplus_{m\ge 1} H_{k-2m+2}(W,\partial W) \longrightarrow 0. 
$$ 
It is actually the case that $i$ and $p$ are, respectively,   
the obvious injection and 
projection. This follows from the results in~\cite{BO}.  
 
Although for 
transversality reasons  we 
stated Proposition~\ref{prop:newMLYau} only for stabilizations of 
subcritical Stein manifolds, we expect it to hold for arbitrary 
subcritical Stein manifolds (see Remark~\ref{rmk:transv_exples}).  
} 
\end{remark}

\subsection{Negative disc bundles}   
 
Let $\cL \stackrel{\pi}{\to} B$ be a Hermitian line bundle over a 
compact, symplectically aspherical manifold $(B,\beta)$ with $c_1(\cL) 
= - [\beta]$. Let $W = \{ v \in \cL : | v | \le 1 \}$ be the 
corresponding disc bundle, with symplectic form $\omega = 
\pi^* \beta + d(r^2 \theta)$, where $r$ is the  
radial coordinate in the fibers and $\theta$ is the angular form. The 
boundary $M = \partial W$ is a contact manifold with contact form 
$\theta$. The closed Reeb orbits on $M$ are the fibers of the natural 
projection to $B$. 
 
Let us assume $\dim_\R B=2n-2$, so that $\dim_\R W=2n$. It follows 
from~\cite[Theorem~D]{O1} that $SH_*(W)=0$, and therefore  
$$ 
SH_*^+(W)\simeq H_{*+n-1}(W,\p W) \simeq H_{*+n-3}(B).  
$$ 
 
\begin{proposition}[\cite{B}] The (linearized) contact homology of the  
prequantization bundle $M=\p W$ is well-defined and equal to 
\begin{equation} \label{eq:HClinebundle}   
HC_*(M,\xi)=\bigoplus_{m=0}^\infty H_{*-2m}(B).  
\end{equation}  
\end{proposition}  
 
\proof As in Remark~\ref{rmk:transv_exples}.(ii), we choose compatible 
almost complex structures $J$ on $W$ and $J_B$ on $B$, and a generic 
Morse function $f:B\to \R$. The Morse-Bott contact 
complex~\cite[\S8.1]{B} is generated by closed Reeb orbits 
$\gamma_{p,k}$ of multiplicity $k$ above critical points $p \in 
\mathrm{Crit}(f)$. The grading of $\gamma_{p,k}$ with respect to the 
symplectic trivialization given by the fiber is $\ind(p;f) + 
2k -\frac12 (2n-2) + n-3 = \ind(p;f) + 2k -2$ \cite[Lemma~2.4]{B}.  
All closed orbits are good, because the parity of the grading of 
$\gamma_{p,k}$ does not depend on the multiplicity $k$. The 
differential in the Morse-Bott complex counts \emph{rigid} 
configurations consisting of $J$-holomorphic curves with gradient 
fragments of $f$ (for this reason, the underlying gluing theorem 
is   
similar to~\cite[Theorem~3.7]{BOauto}). Since 
$J$-holomorphic curves in the symplectization $M \times \R$  
are branched covers of vertical cylinders over closed Reeb orbits, and 
since branch points can always be displaced along the corresponding 
vertical cylinder, it follows that the only such rigid configurations 
are rigid gradient trajectories of $f$. Hence, for each multiplicity 
$k \ge 1$, the contact differential coincides with the Morse 
differential of $f$ on $B$, which proves~\eqref{eq:HClinebundle}.  
\hfill{$\square$}

\begin{proposition}   \label{prop:negdiscbdle}  
 Let $W$ be a disc bundle over a compact, symplectically 
 aspherical manifold $(B,\beta)$ with $c_1(\cL) = - [\beta]$. 
 The exact sequence relating contact and symplectic homology  
 is isomorphic to an exact sequence of the form 
 $$ 
 \scriptstyle 
  \ldots \stackrel 0 \longrightarrow 
   H_{k}(B) \longrightarrow   
 \bigoplus_{m\ge 0} H_{k-2m}(B) \longrightarrow  
 \bigoplus_{m\ge 1} H_{k-2m}(B) \stackrel 0 \longrightarrow H_{k-1}(B) 
 \longrightarrow \ldots   
 $$ 
 \end{proposition} 
 
 \proof Similarly to Proposition~\ref{prop:newMLYau}, we prove by descending 
 induction on $k\in\Z$ that the maps $HC_k(M)\to H_{k+1}(B)$ 
 vanish. The claim is true for $k\ge 2n-2$ for dimensional 
 reasons. Assuming the claim to be true for some $k\in \Z$, the exactness 
 of the sequence  
$${ 
\def\objectstyle{\scriptstyle} 
\def\labelstyle{\scriptstyle} 
\xymatrix 
@C=10pt 
@R=10pt 
{HC_k(M)\ar[r]^0 & H_{k+1}(B) \ar[r] & \bigoplus_{m\ge 0} 
  H_{k+1-2m}(B) \ar[r] & \bigoplus_{m\ge 1} H_{k+1-2m}(B) \ar[r] 
  \ar@{=}[d] & H_k(B)  \\ 
& & & HC_{k-1}(M) &  
} 
} 
$$  
implies, for dimensional reasons, that the last map vanishes, so that 
the claim is true for $k-1$. 
\hfill{$\square$} 
 
\medskip 
 
\begin{remark}  
 The similarity between Proposition~\ref{prop:negdiscbdle} and 
 Proposition~\ref{prop:newMLYau} is best explained via the 
 $S^1$-equivariant approach in~\cite{BO}. The spectral sequence 
 in~\cite{O1} admits an $S^1$-equivariant version which implies that 
 positive $S^1$-equivariant symplectic homology is isomorphic to 
 $H_*(B)\otimes H_*(\mathbb{C}P^\infty)$. The isomorphism between 
 contact homology and positive $S^1$-equivariant symplectic homology 
 therefore implies~\eqref{eq:HClinebundle}. Moreover, the exact 
 sequence in Proposition~\ref{prop:negdiscbdle} is the tensor product 
 of $H_*(B)$ with the Gysin exact sequence of the subcritical pair 
 $(D^2, S^1)$.   
\end{remark}

\subsection{Cotangent bundles} \label{sec:cotangent} 
 
Our next example are unit cotangent bundles  
$$ 
W=DT^*L:=\{p\in T^*L \, : \, |p| \le 1\} 
$$ 
of closed Riemannian manifolds $L$. We recall the  
transversality discussion in Remark~\ref{rmk:transv_exples}    
where we imposed conditions~$(A)$ and~$(B_a)$ in all statements 
involving linearized contact homology (we assume either $\dim\, L\ge 
5$ or $L$ has no contractible closed geodesics, and we work in a free 
homotopy class $a$ containing only simple closed geodesics). As 
mentioned in Remark~\ref{rmk:transv_exples}, we expect these two 
conditions to be completely removed in the future, and we do not 
mention them anymore in the discussion that follows.

The symplectic manifold $W=DT^*L$ is exact with boundary of restricted 
contact type  
$$ 
M=ST^*L:=\{p\in T^*L \, : \, |p|=1\}. 
$$ 
The Liouville form determines a contact structure on $M$ whose isotopy 
class does not depend on the choice of metric, since the 
space of Riemannian metrics is convex.

The first ingredient involved in our long exact sequence are the 
symplectic homology groups $SH_*^a(W)$ in a free homotopy class 
$a$. These have been computed by Viterbo~\cite{Vcotangent}, 
Salamon-Weber~\cite{SaWe} and Abbondandolo-Schwarz~\cite{AbboSch}: 
$$ 
SH_k^a(DT^*L)\simeq H_k(\Lambda^a L), \qquad k\in \Z, 
$$ 
where $\Lambda^aL\subset \Lambda L$ is\footnote{In \cite{SaWe}  
the symbol $\Lambda^aL$ stands for a sublevel set of 
the energy functional, and $a\in \R$.}  
the connected component $a$ in 
the free loop space $\Lambda L$, i.e. the space of continuous 
maps from $S^1$ to $L$.  
The space $\Lambda L$ is endowed with 
the canonical $S^1$-action 
$\big(\theta,\gamma(\cdot)\big)\mapsto\gamma(\cdot+\theta)$, 
 $\theta\in S^1=\R/\Z$, $\gamma\in\Lambda L$. 
 
 The above isomorphism actually holds 
in an improved version involving, on the left hand 
side, the symplectic homology groups truncated by the values of the 
action functional and, on the right hand side, the relative homology 
groups of sublevel sets of the energy functional on the loop 
space~\cite{AbboSch,SaWe}. When $a=0$, in particular, the 
tautological exact sequence~\eqref{eq:Vseq} is identified with the 
exact sequence of the pair $(\Lambda^0 L,L)$ and we have  
$$ 
SH_k^+(DT^*L) \simeq H_k(\Lambda^0 L,L).  
$$ 
 
The second ingredient in our long exact sequence are the linearized 
contact homology groups $HC_*^{i^{-1}(a)}(ST^*L)$, where $i$ is the map 
associating to the free homotopy class of a loop in $ST^*L$ its free 
homotopy class as a loop in $DT^*L$. An argument similar to that of  
Section~\ref{sec:subcritical} shows that, if $\dim\, L\ge 3$, the map $i$ is 
bijective and $i^{-1}(0)=0$.  
An ongoing project of K.~Cieliebak and J.~Latschev~\cite{CiLa} aims at 
computing the entire SFT of $ST^*L$.   
The theorem below is a particular case of their more general results.  
 
\begin{theorem}[Cieliebak-Latschev~\cite{CiLa}] \label{thm:CiLa} 
 Let $L$ be a closed oriented Riemannian manifold. Given a free homotopy 
class of loops $a$ in $DT^*L$ (hence in $L$), the  
following isomorphisms hold for $k\in\Z$: 
\begin{eqnarray}  
  HC_k^{i^{-1}(a)}(ST^*L) & \simeq & H_{k-(n-3)}(\Lambda^a L/S^1),  
\qquad a\neq 0, \nonumber \\ 
  HC_k^{i^{-1}(0)}(ST^*L) & \simeq & H_{k-(n-3)}(\Lambda ^0 L/S^1,L).  
\label{eq:ST*Lcontractible} 
\end{eqnarray} 
\end{theorem} 
 
\begin{remark} {\rm  
 One can rephrase Theorem~\ref{thm:CiLa} within the setting of 
 $S^1$-equivariant homology. Given a topological space $M$ endowed 
 with an $S^1$-action, the $S^1$-equivariant homology groups are 
 defined as  
$$ 
H_*^{S^1}(M):=H_*(M\times_{S^1} ES^1),  
$$ 
where $ES^1$ is a contractible 
space on which $S^1$ acts freely (for example $ES^1=S^\infty$, the 
infinite dimensional sphere). If we work with rational coefficients 
and $S^1$ acts with finite  -- hence cyclic -- isotropy groups, then 
$S^1$-equivariant homology is isomorphic to the  
homology of the quotient. The reason is that the map $M\times _{S^1} 
ES^1 \to M/S^1$ induced by the projection on the first factor behaves 
like a fibration with fibers $B\Z/k\Z$, $k\ge 1$ and the latter are 
$\Q$-acyclic. As a consequence we obtain that   
$$ 
H_*(\Lambda^0L/S^1,L) \simeq H_*^{S^1} (\Lambda^0L,L). 
$$ 
} 
\end{remark} 
 
The above discussion can be summarized as follows.  
 
\begin{proposition} \label{prop:newcotangent} 
Given a closed oriented Riemannian manifold $L$, the long exact 
sequence relating contact and symplectic homology for $DT^*L$ in the 
free homotopy class $a$ is isomorphic to an exact sequence of the form 
$$ 
\ldots \!\to\! H_k(\Lambda ^a L) \!\to\! H_k^{S^1}(\Lambda^a L) \!\to\!  
H_{k-2}^{S^1}(\Lambda^a L) \!\to\! H_{k-1}(\Lambda ^a L) \!\to\! \ldots  
$$ 
if $a\neq 0$, respectively 
$$ 
...\!\to\! H_k(\Lambda ^0\!L,L)\!\to\!  H_k^{S^1}(\Lambda^0\!L,L) \!\to\! 
H_{k-2}^{S^1}(\Lambda^0\!L,L)\!\to\!  H_{k-1}(\Lambda ^0\!L,L) \!\to\!...  
$$  
if $a=0$. \hfill{$\square$} 
\end{proposition} 
It follows from the results in~\cite{BO} that these are the classical 
Gysin exact sequences for $S^1$-equivariant homology.


\end{document}